\documentclass{cmslatex}
\usepackage[paperwidth=7in, paperheight=10in, margin=.875in]{geometry}
 \usepackage[backref,colorlinks,linkcolor=red,anchorcolor=green,citecolor=blue]{hyperref}
\usepackage{amsfonts,amssymb}
\usepackage{amsmath}
\usepackage{graphicx}
\usepackage{cite}
\usepackage{enumerate}
\usepackage{comment}
\usepackage{dsfont}
\usepackage{tikz,siunitx} % to draw domains 
\usepackage[caption=false]{subfig} % <- Preamble
\newcommand{\ueps}[0]{u^{\epsilon}}
\newcommand{\peps}[0]{p^{\epsilon}}

\newcommand{\pderiv}[3][]{% \pderiv[<order>]{<func>}{<var>} 
  \ensuremath{\frac{\partial^{#1} {#2}}{\partial {#3}^{#1}}}}

\newcommand{\bigoh}[1]{\mathcal{O}\left( #1 \right)}
\newcommand\norm[1]{\left\lVert#1\right\rVert}
\newcommand{\chevron}[1]{\left\langle #1 \right\rangle}
\newcommand{\Omegamac}[0]{\Omega^{\rm mac}}
\newcommand{\Omegamic}[0]{\Omega^{\rm mic}}
\newcommand{\U}{\mathcal{U}}

\usepackage{mathtools}

\DeclarePairedDelimiter\floor{\lfloor}{\rfloor}

\renewcommand{\theequation}{\arabic{section}.\arabic{equation}}

   \allowdisplaybreaks

\graphicspath{ {figs/} }
\begin{document}
 \title{Heterogeneous multiscale methods for rough-wall laminar viscous flow
\thanks{}}

          %For each author, make a block with the following macros:

\author{Sean P. Carney\thanks{Department of Mathematics, University of California, Los Angeles, CA, 
90095, USA (spcarney@math.ucla.edu). }
\and Bj\"{o}rn Engquist \thanks{Department of Mathematics and Oden Institute for Computational Engineering
and Sciences, The University of Texas at Austin, Austin, TX, 78712, USA (engquist@oden.utexas.edu).}
}

\pagestyle{myheadings} \markboth{MULTISCALE METHODS FOR VISCOUS LAMINAR FLOW}{SEAN P. CARNEY AND BJ\"{O}RN ENGQUIST} 
\maketitle

\begin{abstract}
We develop numerical multiscale methods for viscous fluid flow over a rough boundary. 
The goal is to derive effective boundary conditions, or wall laws, 
through high resolution simulations localized to the boundary 
coupled to a coarser simulation in the domain interior following the
framework of the heterogeneous multiscale method. 
Rigorous convergence of the coupled system is shown in a
simplified setting. 
%Under simplifying assumptions, the 
%wall law produced by the multiscale method is quantitatively 
%compared to that derived from periodic homogenization theory. 
Numerical experiments illustrate the utility of the 
method for more general roughness patterns and far field 
flow conditions. 

\end{abstract}
\begin{keywords} Rough boundaries; multiscale methods; fluid dynamics
\end{keywords}

 \begin{AMS} 65N22; 65N30; 76D05
\end{AMS}

%%%%%%%%%%%%%%%%%%%%%%%%%%%%%%%%%%%%%%%%%%%%%%%%%%%%%%%%%%%%%%%%%%%%%%%%%%%%%%%%%%%%%%%%%%%%%%%%%%%%%%%%%%

\section{Introduction}

Standard partial differential equations for viscous flow such as 
the Stokes and Navier-Stokes equations naturally have no-slip boundary conditions. 
The velocity vector $u = 0$ at the boundary. 

%There are cases where this is not accurate. At a slip line where two immiscible 
%fluids meet at a solid boundary is one example \cite{E_contact_lines,dussan_contact_lines}. The fluid
% molecules actually slide along the boundary near the slip line. Another example is gas at very low pressure. 
%
%Even when the no-slip boundary condition is valid, it 
The no-slip condition however    % this is the new line
can generate a large separation of  
scales in a problem. The analysis and simulation of 
fluid flows over rough surfaces, for example, 
is challenging because the no-slip condition generates
boundary layers in the vicinity of the roughness whose 
resolution can be computationally expensive. 

Surface roughness plays an important role in a variety of physical applications.
In geophysical fluid dynamics, 
meteorological flows are known to be affected by mountain ranges, city landscapes, 
and wavy seas, while ocean currents are impacted by variations in the 
ocean floor as well as the coastline \cite{pedlosky}.  %\cite{Bresch2005,pedlosky} 
Rough surfaces also can also effect a reduction 
in skin friction drag; the morphology of a swordfish's sword, shark dermal 
denticles, and riblets on the Stars and Stripes yacht in the 1987 America's
Cup finals are all examples 
\cite{drag_reduction_nature,mikelic2009}. In hypersonic flows, surface
roughness is relevant because space shuttles typically 
contain periodic gaps between 
covering tiles designed for heat control \cite{Carrau:1991}.
Outside of fluid mechanics, electromagnetic wave scattering
by an obstacle is known to depend on small-scale rugosity, 
or imperfections, along its surface \cite{achdou1}.

In such cases where the no-slip condition either (i) is inaccurate
or (ii) generates sharp boundary layers, it may be better in a computation to replace  
it with an effective boundary condition, 
or wall-law. Ideally, the wall-law captures the effect of the 
asymptotic small scales on the large scales. Computing in a 
domain without the small scale structure then results in a large
reduction in the degrees of freedom necessary in a simulation. 

Sometimes effective boundary conditions can be rigorously 
justified from first principles. Relevant to the present work is 
the Navier-slip law for viscous laminar flow. Other examples include the 
the Leontovitch boundary conditions for
electromagnetic wave scattering, and the Beavers-Joseph-Saffman
law for viscous, laminar flow over a porous bed.  
%wave scattering over rough surfaces, respectively, that can 
%be rigorously justified \cite{jager,achdou1,achdou2,artola}. 
%Another example is the Beavers-Joseph-Saffman law for viscous
%flow over a porous bed \cite{jager_porous}. 
There are other well-known wall-laws for which a theoretical 
justification is lacking, however. The logarithmic 
law-of-the-wall for wall-bounded turbulent flows \cite{Millikan:1938tv,piomelli:2002}
and the electro-osmotic slip velocity (and associated zeta potential)
for electrokinetic flows \cite{squires_bazant_2004,yariv:2012} are two examples. 

The focus here will be on viscous laminar flow over a rough boundary that varies with 
characteristic amplitude and period $0<\epsilon\ll 1$. The problem is well understood
mathematically and has been extensively analyzed with the tools of asymptotic 
homogenization theory. Rigorous estimates have been obtained for the linear, Stokes
case \cite{amirat}, and for Poiseuille \cite{jager} and Couette \cite{jager_couette}
flows governed by the stationary Navier-Stokes equations in channel geometries with
periodic roughness, as well as random, ergodic roughness \cite{basson2008,dalibard}.
In Ref.~\cite{mikelic2013} the authors obtain rigorous convergence results 
for more general, three-dimensional flows under the assumption that the solution to 
the corresponding problem with the 
no-slip boundary condition imposed along a smooth boundary possesses a smooth solution.
Related mathematical studies also include \cite{jager_porous} for channel flow over a porous bed,
\cite{Bresch2005,gerard_varet:2003} for meteorological flows, and 
\cite{achdou1,achdou2,artola} for Maxwell's equations. 

There additionally exists several numerical studies of rough-wall viscous flow. Computational 
techniques based on domain decomposition \cite{achdou:1998dd} and
asymptotic expansions \cite{pironneau} have been previously proposed 
for modeling the effect of surface roughness on the flow in the domain interior.  
Similar strategies have also been 
explored for compressible flows over rough surfaces \cite{deolmi__2017}
and for shape optimization with the purpose of minimizing drag for
both laminar and turbulent incompressible flows, assuming the roughness 
is within the viscous sublayer \cite{friedmann:2010,friedmann:2011}. 
More recently, Bottaro and Naqvi \cite{bottaro:2020} and L\={a}cis, Bagheri, and 
coauthors \cite{Lacis:2017,Sudhakar2019HigherOH}
used asymptotic expansions to develop high-order wall-laws for flows over 
rough surfaces and porous media, respectively. In the latter groups' work,
wall-laws for rough surfaces are included as a limiting case.

The goal in this work is to derive the wall-law that 
describes the average effect of surface roughness on 
the large scale flow by local, high-resolution 
`microscale' simulations. These are in turn coupled to a coarse-scale 
simulation in the interior following the framework of 
the heterogeneous multiscale method (HMM) 
\cite{abdulle:2012,hmm_review}.  This `macroscale'
simulation will use the effective boundary condition along a
smooth boundary. 

The proposed method shares some features of classical domain 
decomposition and adaptive mesh refinement, however it 
is fundamentally different in that the classical methods 
try to resolve the $\epsilon$-scale in the neighborhood of 
the rough boundary throughout the entire computational domain, 
while here local refinement is only used selectively to derive 
effective boundary conditions. The domain size of the local 
simulations scales with $\epsilon$, and hence the overall numerical 
degrees of freedom are independent of $\epsilon$. The overall 
computational cost is thus drastically reduced compared to a 
full discretization of the problem. 

%%%%%% MOVE TO HMM SECTION %%%%%%%%%%%%%%
%The HMM framework has previously been applied to the slip line problem mentioned 
%above \cite{hmm_microfluidics}.
% The local high-resolution model is then molecular dynamics. The outer 
%coarser scale model is the Navier-Stokes equations, which gets an 
%effective boundary condition at the slip line from the molecular dynamics 
%simulation. 

The existing mathematical homogenization theory for viscous, rough-wall 
flow justifies 
an effective boundary condition of Robin-type for the 
case of periodic or random, stationary ergodic roughness. 
Guided by the theory, the wall-law in the 
proposed method is also a Robin condition; however the 
method is designed to be applicable to general roughness
profiles, so long as there is a separation of scales 
between the roughness characteristics and the size of the 
full domain. 

For the case of Stokes flow in a channel with periodic
roughness, we exploit linearity to prove several 
properties of the proposed method. We first show that a 
back-and-forth iteration used to solve the coupled 
macro/microscale system converges to a fixed point, 
and that the resulting wall-law leads to a coercive
 macroscale problem. Moreover, we show the convergence
happens rapidly; the coefficients in the wall-law 
produced by the first and second iterations differ
by less than $\bigoh{\epsilon^2}$. 
Finally, we show that the macroscale approximation
produced by the HMM and the true rough-wall flow both
limit to the same quantity as $\epsilon$ vanishes.
We do not show, however, 
any rigorous estimates relating the HMM approximation
to the true rough-wall flow, as the creation of 
new results in homogenization theory is outside 
of the scope of the present work.

The structure of the paper is as follows. Section \ref{sec:background} 
provides background first for the problem of laminar, viscous 
flow over a rough boundary and then for the heterogeneous 
multiscale method. Included is a brief asymptotic analysis 
of the problem that slightly modifies that of Achdou et al.\ 
\cite{pironneau}. The asymptotics motivate the form of the 
effective boundary condition used in Section \ref{sec:algorithm_for_laminar_flow}, 
where the heterogeneous multiscale method is proposed and 
analyzed. Numerical results in Section \ref{sec:laminar_numerics} 
illustrate the method accurately and efficiently 
captures the average effect of surface roughness both where
the homogenization theory is applicable and in more general 
settings. Finally, Section \ref{sec:conclusions} contains some 
concluding remarks.

\section{Background}
\label{sec:background}

There exists a large amount of mathematical results available in
 the literature concerning the asymptotic behavior of laminar
 incompressible flow in the presence of a rough boundary of characteristic
height and wavelength $\epsilon$; see
\cite{achdou2,pironneau,achdou:1998dd,amirat,basson2008,dalibard,jager,jager_couette,Sudhakar2019HigherOH}, 
mentioned in the introduction, as well as references therein.
 %\cite{achdou1,achdou2,pironneau,artola,amirat,basson2008,dalibard,jager}, and references therein.
 Although there are differing physical assumptions and
 levels of mathematical rigour associated to each work,
 all justify the use of a Robin-type condition
 on a smooth boundary near the original rough boundary. 
If $u_{\tau}$ denotes the fluid velocity vector tangential to the smooth boundary,
 and $n$ denotes the unit vector normal 
to the boundary, 
the wall-law is of the form
\begin{equation}\label{eq:ex_wall_law}
u_{\tau} = \alpha \nabla u_{\tau} \cdot n. 
\end{equation}
This is typically coupled with the no-penetration condition $u \cdot n = 0$. 
In the available homogenization theory, either periodic roughness or 
random, stationary ergodic roughness is assumed. Furthermore, the smooth boundary 
is assumed to be flat; it aligns with a coordinate axis and exhibits no curvature.
In these settings, the slip amount $\alpha$ is a
constant given by the average of a local corrector that decays exponentially fast
in the wall-normal variable. In more general situations, however, it should be
noted that $\alpha$ could exhibit tangential variations.

Regardless of whether or not $\alpha$ is constant or varies 
tangentially, in two dimensions it is a scalar value. In three 
dimensions, however, it could more generally be a tensor. 
We focus in this work on the two-dimensional case for ease
of exposition and so that numerical examples of the full, rough-wall 
problem are feasible on a workstation. All that follows, 
however, can be readily generalized to three dimensions, and
Remark \ref{rmk:3d_case} in Section \ref{sec:algorithm_for_laminar_flow} 
briefly describes how our scheme generalizes to this case.

To see how a wall-law of the form \eqref{eq:ex_wall_law} can arise, 
we now briefly conduct formal asymptotic analysis of two-dimensional, 
viscous laminar flow in rough domain.  
The analysis is nearly identical to that of Achdou et al. \cite{pironneau}, 
and below we quote two of their theorems relevant to current work.
In particular the asymptotics lead to the same linear cell problem 
whose solution is the first order corrector in the asymptotic expansion. 
After scaling by the small-scale parameter $\epsilon$, the coefficient 
$\alpha$ in \eqref{eq:ex_wall_law} then is simply the average 
of this solution.

There are two differences. The first is that we consider a rough boundary 
parameterized by the product of a rapidly oscillating periodic function $\varphi^{\epsilon}$
and a slowly varying function $\beta$. This simply results in the modulation of
the slip-amount by $\beta$.
The second difference is that the authors in 
\cite{pironneau} consider high Reynolds number ($Re$) stationary flow, where $Re$ is proportional to 
$1/\epsilon$. Since such flows are expected to be sensitive to turbulent 
instabilities in the presence of rough boundaries on physical grounds, we 
consider the asymptotic regime of low Reynolds number flow. The low $Re$ case is 
additionally in line with the more rigorous analysis in \cite{jager,jager_couette,basson2008}.

We preface the asymptotics with some definitions to be used throughout the paper.

\label{subsec:analysis_of_laminar_flow}
\subsection{Preliminary definitions. }\label{subsec:definitions}
First, let $e_1$ and $e_2$ be the standard unit vectors in $\mathbb{R}^2$, and let 
$x_1$ and $x_2$ parameterize the horizontal and vertical directions as in 
Figure \ref{fig:channel_modulated_domain}. 
%%%%%%% cell prob domain %%%%%%%%%%%
\begin{figure}[h!]
\centering
\begin{tikzpicture}[xscale = 10.875, yscale = 6.3] %[xscale=6.5,yscale=3.8]
\draw[black, domain = 0.2:0.3] plot(\x, {0.05*(cos(20*pi*\x r) -1)});
\draw[->, black] (0.2,0) -- (0.2,0.35);
\draw[black] (0.2,0.25);
\draw[<-, black] (0.3,0.35) -- (0.3, 0);
\draw[black,densely dotted] (0.2,0.)node[left]{$y_2=H$} -- (0.3, 0.);
\draw[black,densely dotted] (0.2,-0.1)node[left]{$y_2=0$} -- (0.3, -0.1);
\draw (0.29,-0.05)node[right]{$\partial Y$}; 
\end{tikzpicture}
\caption{Semi-infinite domain $Y$ with lower boundary $\partial Y = \{(y_1,y_2)\in [0,1) \times \mathbb{R}\,\,\vert \,\,y_2=\varphi(y_1)\}$ for 
some periodic function $\varphi$. }
\label{fig:per_cell}
\end{figure}
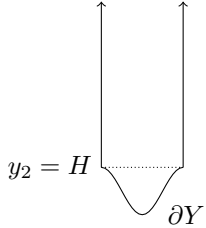

\begin{figure}[h!]
\centering
\begin{tikzpicture}[xscale=3.75,yscale=1.9]
%%%%%%% rough, full domain %%%%%%%%%%%
\draw (0.5,-0.075) node[below]{$\Gamma^{\epsilon}$}; 
\draw (0.5,0.75) node[below]{$\Omega^{\epsilon}$}; 
\draw[variable=\y, domain=0:1] plot({0.4*\y*\y-0.4*\y},\y);
\draw[variable=\y, domain=0:1] plot({-0.4*\y*\y+0.4*\y+1},\y);
\draw[black, domain = 0:1] plot(\x, {-1.6*\x*\x*\x*\x+3.2*\x*\x*\x-2.4*\x*\x+0.8*\x+1});
\draw[black, domain = 0:0.1] plot(\x, {0.025*(cos(40*pi*\x r) -1)});
\draw[black, domain = 0.1:0.2] plot(\x, {0.025*(cos(40*pi*\x r) -1)});
\draw[black, domain = 0.2:0.3] plot(\x, {0.025*(cos(40*pi*\x r) -1)});
\draw[black, domain = 0.3:0.4] plot(\x, {0.025*(cos(40*pi*\x r) -1)});
\draw[black, domain = 0.4:0.5] plot(\x, {0.025*(cos(40*pi*\x r) -1)});
\draw[black, domain = 0.5:0.6] plot(\x, {0.025*(cos(40*pi*\x r) -1)});
\draw[black, domain = 0.6:0.7] plot(\x, {0.025*(cos(40*pi*\x r) -1)});
\draw[black, domain = 0.7:0.8] plot(\x, {0.025*(cos(40*pi*\x r) -1)});
\draw[black, domain = 0.8:0.9] plot(\x, {0.025*(cos(40*pi*\x r) -1)});
\draw[black, domain = 0.9:1.0] plot(\x, {0.025*(cos(40*pi*\x r) -1)});
%%%%%%% axes %%%%%%%%%%%
\draw[->] (-.4,-0.075)--(-0.4,0.4) node[right]{$x_2$}; 
\draw[->] (-.4,-0.075)--(-0.2,-0.075) node[right]{$x_1$}; 
%%%%%%% second, smooth domain %%%%%%%%%%%
\draw (1.4,-0.075) -- (2.4,-0.075);
\draw (1.9,-0.075) node[below]{$\Gamma^{0}$}; 
\draw (1.9,0.75) node[below]{$\Omega^{0}$}; 
\draw[black, domain = 1.4:2.4] plot(\x, {-1.6*(\x-1.4)*(\x-1.4)*(\x-1.4)*(\x-1.4)+3.2*(\x-1.4)*(\x-1.4)*(\x-1.4)-2.4*(\x-1.4)*(\x-1.4)+0.8*(\x-1.4)+1});
\draw[variable=\y, domain=-0.075:1] plot({0.26087*\y*\y-0.241304*\y+1.38043},\y);
\draw[variable=\y, domain=-0.075:1] plot({-0.26087*\y*\y+0.241304*\y+2.41957},\y);
\draw[variable=\y, domain=0:1] plot({-0.4*\y*\y+0.4*\y+1},\y);
\end{tikzpicture}
\caption{Domain $\Omega^{\epsilon}$ with periodic, sinusoidal roughness and corresponding 
$\Omega^0$ with flat boundary.}
\label{fig:channel_modulated_domain}
\end{figure}

Let $\varphi:\mathbb{R}\to\mathbb{R}$ be a bounded, Lipschitz continuous, periodic
function with maximum value $H:= \norm{\varphi}_{\infty}$ that satisfies $\varphi(N)=H$ for every $N\in\mathbb{Z}$, and 
$\varphi(t+1) =\varphi(t) \ge 0$  $\forall t\in \mathbb{R}$. 
Let $\epsilon$ be some fixed small parameter,  $0<\epsilon \ll 1$, and define $\varphi^{\epsilon}(x_1) := \epsilon \varphi(x_1/\epsilon)$. 
For smooth, bounded function $\beta:\mathbb{R}\to\mathbb{R}$ that is independent
of $\epsilon$, define $\zeta^{\epsilon}(x_1) := \beta(x_1)\varphi^{\epsilon}(x_1)$ to be the
function that parameterizes the rough boundary. Without loss of generality, assume
that $\norm{\beta}_{\infty} = 1$ for ease of exposition below. Further assume that $\beta$ is 
bounded below by a positive constant, so that $\beta \ge \beta_{\ast} > 0$ and hence 
$\zeta^{\epsilon}(x_1) \ge 0$ $\forall x_1 \in \mathbb{R}$. 

Let $Y = \left\{ (y_1,y_2) \in [0,1)\times \mathbb{R}\vert \,\, y_2 \ge \varphi(y_1)\right\}$ be a
domain containing a periodic ``cell" of $\varphi$, semi-infinite in the vertical direction with lower boundary 
$\partial Y = \{y\in[0,1)\times\mathbb{R}\vert y_2=\varphi(y_1)\}$; for example, see Figure \ref{fig:per_cell}.
Let $L^2_{\rm per}(Y)$ be the space of square integrable functions in $Y$ that are $1$-periodic in the 
$y_1$ variable, and let $H^1_{\rm per}(Y) \subset L^2_{\rm per}(Y)$ be the subspace whose first derivatives 
also belong to $L^2_{\rm per}(Y)$. Also, let $\mathcal{S}_{\rm per}(Y) \subset L^2_{\rm per}(Y)$ be the 
subspace of functions that decay exponentially fast in $y_2$, as well as all of their derivatives.  

Let $\Theta^{\epsilon} = \left\{(x_1,x_2)\in\mathbb{R}^2\vert \,\, x_2 \ge \zeta^{\epsilon}(x_1)\right\}$ be the
 semi-infinite domain contained in the upper half plane 
$x_2 \ge 0$ in $\mathbb{R}^2$, and let $\Omega$ be a bounded domain in $\mathbb{R}^2$ made
of one piece that intersects the line $\{x_2=0\}$. Take $\Omega^0 = \Omega \cap \{x_2 \ge0\}$, and
let $\Gamma^0 = \partial \Omega^0 \cap \{x_2 = 0\}$. Finally, take 
$\Omega^{\epsilon} := \Theta^{\epsilon} \cap \Omega^0$, so that $\Omega^{\epsilon}$ has 
a rough boundary $\Gamma^{\epsilon}$ with characteristic amplitude and wavelength $\epsilon$; 
for example, see Figure \ref{fig:channel_modulated_domain}. 
 Note that $\overline{\Omega^{\epsilon}} \to 
\overline{\Omega^0}$ as $\epsilon \to 0$.

\subsection{Asymptotic analysis.}\label{subsection:asymptotic_analysis}
Given $\Omega^{\epsilon}$, consider the following stationary Navier-Stokes problem
\begin{align}\label{eq:NS_eps}
\mathcal{L}(\ueps,p^{\epsilon}):=  
- \nu \Delta \ueps + \ueps\nabla \ueps  + \nabla p^{\epsilon} &= f   \qquad \text{in } \Omega^{\epsilon} \nonumber \\
\nabla \cdot \ueps   &= 0  \qquad \text{in } \Omega^{\epsilon} \nonumber \\
\ueps &= 0 \qquad \text{on } \partial \Omega^{\epsilon}
\end{align}
where $\nu = \bigoh{1}$. 
Other combinations of well-posed boundary conditions are possible, 
so long as the no-slip condition $u^{\epsilon}=0$ is imposed on the rough wall 
$\Gamma^{\epsilon}$. 

First consider the approximation $(u^0, p^0)$ that satisfies 
\begin{align*}%\label{eq:NS_0}
\mathcal{L}(u^0,p^{0})   &= f   \qquad \text{in } \Omega^{0} \nonumber \\
\nabla \cdot u^0   &= 0  \qquad \text{in } \Omega^{0} \nonumber \\
u^0 &= 0 \qquad \text{on } \partial \Omega^{0}.
\end{align*}
A simple Taylor expansion shows that the error along the rough boundary 
$\Gamma^{\epsilon}$ is $\bigoh{\epsilon}$; indeed let 
$x^0 \in \Gamma^0$, $x^0 + \zeta^{\epsilon}(x_1^0) e_2 \in \Gamma^{\epsilon}$. Then 
\begin{align}
u^0\left(x^0 +\zeta^{\epsilon}(x^0_1)e_2 \right) &= u^0(x^0) + \zeta^{\epsilon}(x_1^0) \frac{\partial u^0}{\partial x_2} (x^0) 
+ \frac{1}{2}\left(\zeta^{\epsilon}(x_1^0)\right)^2 \frac{\partial^2 u^0}{\partial x_2^2}(\xi(x^0)) \nonumber  \\
&=  \epsilon \beta(x_1^0) \varphi(x_1^0/\epsilon)  \frac{\partial u^0}{\partial x_2} (x^0) + 
\mathcal{O}(\epsilon^2). \label{boundary_err}
\end{align}
The error can be improved by considering a higher order approximation that 
accounts for the geometry of the rough boundary. 

Consider next the approximations 
\begin{align} 
\label{eq:expansions}
\ueps(x) &\approx u^1(x) + \epsilon \, u_{BL}^1(x,x/\epsilon)
=  u^1(x) + \epsilon \, \beta(x_1) \pderiv{u^1_1}{x_2}(x_1,0) \left(\chi(x/\epsilon) - \overline{\chi}e_1\right) \nonumber \\
\peps(x) &\approx p^1(x) + p_{BL}^1(x,x/\epsilon)
= p^1(x) + \beta(x_1) \pderiv{u^1_1}{x_2}(x_1,0) \pi(x/\epsilon),
\end{align}
where $(u^1, p^1)$ satisfy
\begin{align}\label{eff_prob_p}
\mathcal{L}(u^1, p^1) &= f  \qquad \text{in } \Omega^0 \nonumber \\
\nabla \cdot u^1 &= 0 \qquad \text{in } \Omega^0 \nonumber \\
u^1(x)-  \epsilon \beta(x_1) \overline{\chi_1}\, \pderiv{u^1_1}{x_2}(x) e_1 
&=0 \qquad x\in \Gamma^0  \nonumber \\
u^1 &= 0 \qquad \text{on } \partial \Omega^0\setminus \Gamma^0. 
\end{align}
The constant $\overline{\chi_1}$ in the slip boundary condition 
posed along $\Gamma^0$ is defined to be the horizontal average of the horizontal component of 
the vector $\chi$ at $y=H$, that is, just above the roughness $\varphi$:
\begin{equation*}
\overline{\chi_1} = \int_0^1 \chi_1(y_1, H) \, dy_1.
\end{equation*} 
The pair $(\chi, \pi)$ satisfy  
\begin{align}\label{eq:cell_prob_mod}
-\mu \Delta_y\chi + \nabla_y \pi &= 0 \qquad\qquad\,\,\,\, \text{in }Y \nonumber \\
\nabla_y \cdot \chi &= 0 \qquad\qquad \,\,\,\,\text{in } Y \nonumber \\
\chi(y) &= -\varphi(y_1) \qquad y\in \partial Y \nonumber \\ 
\chi - \overline{\chi}  &\in H^1_{\rm per}(Y) \nonumber \\
\pi &\in L^2_{\rm per}(Y) 
\end{align}
and hence only depend on the roughness geometry.  

Another Taylor expansion shows that indeed $u^1 + \epsilon u^1_{BL}$ is $\bigoh{\epsilon^2}$ 
along $\Gamma^{\epsilon}$, i.e.\ one order higher than $u^0$ as desired. Using the formal 
differentiation rule 
\begin{equation*}
\nabla \Phi(x,x/\epsilon) = \nabla_x \Phi(x,y) + \frac{1}{\epsilon} \nabla_y \Phi(x,y), 
\end{equation*}
inserting the approximations \eqref{eq:expansions} into $\mathcal{L}$, and using 
\eqref{eff_prob_p} and \eqref{eq:cell_prob_mod} give that the approximations
\eqref{eq:expansions} are $\bigoh{\epsilon}$ pointwise in $\Omega^0$.

The first theorem we quote from \cite{pironneau} asserts that the so-called
cell-problem \eqref{eq:cell_prob_mod} has a unique solution that decays 
exponentially in $y_2$, and that $\overline{\chi_2} = 0$. This means  
the influence of the correctors $u^1_{BL}$ and $p_{BL}$ in \eqref{eq:expansions}
are only asymptotically felt in a $\bigoh{\epsilon \log(1/\epsilon)}$ neighborhood 
of the rough boundary $\Gamma^{\epsilon}$. 
%The following two theorems are due to Achdou et al. The first gives 
%rigorous backing to the assertion that the local 
%correctors $(u^1_{BL}, p^1_{BL})$ decay exponentially fast in the fast variable $x/\epsilon$, 
%so that the term \eqref{eq:closer_look_term} resulting from the 
%insertion of the asymptotic expansions \eqref{asym_exp_p1} and \eqref{asym_exp_p2}
%into $\mathcal{L}$ is indeed $\bigoh{\epsilon}$ throughout $\Omega^0$.  
%%%
%%% Pironneau theorem on exponential decay 
%%%

\smallskip
\begin{theorem}\label{thrm1_p}
[Achdou et al.]\cite{pironneau} 
There exists a unique pair of functions $(\chi, \pi)$ and a unique 
vector $\overline{\chi} \in \mathbb{R}^2$ such that 
$\chi-\overline{\chi} \in \left(H^1_{\rm per}(Y)\right)^2 \cap 
\left(\mathcal{S}_{\rm per}(Y)\right)^2, \pi \in L^2_{\rm per}(Y) \cap 
\mathcal{S}_{\rm per}(Y)$ and \eqref{eq:cell_prob_mod} is satisfied in a weak sense. 
Furthermore, $\overline{\chi}$ is horizontal, 
\begin{equation*}
\overline{\chi} = \overline{\chi_1} e_1 . 
\end{equation*}
\end{theorem}

The second quoted theorem provides a bound on the size of the constant $-\overline{\chi_1}$
and is crucial for the well posedness of the effective problem \eqref{eff_prob_p}
and hence for its numerical approximation as well. 
\smallskip
\begin{theorem}\label{thm:thrm2_p}
[Achdou et al.]\cite{pironneau} Let $H:= \max_{y \in \partial Y} y = \norm{\varphi}_{\infty}$. 
Then the constant $-\overline{\chi_1}$ satisfies the bound
\[
0 \le -\overline{\chi_1} \le H.
\]
\end{theorem}
As a result, the problem \eqref{eff_prob_p} is generally ill-posed; its variational form
contains the term 
\begin{equation*}
\frac{\mu}{\overline{\chi}} \int_{\Gamma^0} u^1_1 v_1/\beta  \, ds 
\end{equation*}
where $v$ is some test function, which is not coercive when $\overline{\chi}< 0$. For this reason,
the effective boundary condition is posed just \emph{above} the roughness, as oppposed
to at $y=0$. This is consistent with location of the effective boundary condition in 
the more rigorous studies \cite{jager,jager_couette,basson2008}. 
For any $\delta > \epsilon H$, the problem \eqref{eff_prob_p} becomes
\begin{align*}\label{eff_prob_p_mod}
\mathcal{L}(u^1, p^1) &= f  \qquad \text{in } \Omega_{\delta}^0 \nonumber \\
\nabla \cdot u^1 &= 0 \qquad \text{in } \Omega_{\delta}^0 \nonumber \\
u^1(x)-  \beta(x_1) (\epsilon \overline{\chi}+\delta)\, \pderiv{u^1_1}{x_2}(x) &=0 \qquad x\in \Gamma_{\delta}^0 \nonumber  \\
u^1 &= 0 \qquad \text{on } \partial \Omega_{\delta}^0\setminus \Gamma_{\delta}^0,
\end{align*}
where 
\begin{equation}\label{eq:Omega_epsH}
\Omega^0_{\delta} = \Omega^0 \cap \{x_2\ge \delta\}, \qquad
\Gamma^0_{\delta} = \{x+(0,\delta), x \in \Gamma^0 \}.
\end{equation}
%so that $\Omega^0_{\delta} \subset \Omega^{\epsilon} \subset \Omega^0$ and the 
%wall law is imposed at the height $x_2 = \epsilon H$, 
%that is, \emph{just above} the rough boundary $\Gamma^{\epsilon}$. 
%A Taylor expansion of the effective boundary condition \eqref{eq:ebc_pironneau}
%implies that the approximation
%\begin{equation}\label{eq:final_expansion}
%u^1(x) + \epsilon \beta(x_1) \pderiv{u^1_1}{x_2}(x_1,\epsilon H) 
%\left(\chi(x/\epsilon)-\overline{\chi}\right) = \bigoh{\epsilon^2}, \qquad x\in \Gamma^0_{\epsilon H} 
%\end{equation}
%still holds.
%
For this reason, the wall-law determined by the multiscale method 
described in Section \ref{sec:algorithm_for_laminar_flow} is likewise 
imposed just above the surface roughness. 

%%%%%%%%%%%%%%%%%%%%%%%%%%%%%%%%%%%%%%%%%%%%%%%%%%%%%
%\subsection{\red{Convergence results of \cite{jager}}}
%\red{Worth citing the theorems from J\"{a}ger, Mikeli\'{c} here?}

\subsection{Heterogeneous multiscale method.}
\label{subsec:hmm}
The heterogeneous multiscale method (HMM) is a general framework for
 designing multiscale algorithms that aims to capture the macroscopic
 behavior of a system without resolving the microscopic details in their
 entirety. Under the assumption of scale separation in the underlying
 physical system, HMM couples macroscopic simulations to local, 
microscopic simulations so that the simulation has an overall computational complexity
 independent of the fine scale. Comprehensive introductions to and reviews of HMM
 can be found in \cite{abdulle:2012,hmm_intro,hmm_review}; below we briefly describe 
the main idea of the method and its applicability to designing effective
 boundary conditions for fluid simulations. 

Suppose there is a general model for the macroscopic state of a physical system that
 can be expressed as $M(\Psi,D) = 0$, where $D$ represents the macroscopic data 
necessary for the model to be complete. Then the main goal of HMM is to
 approximate $D$ by solving microscale problems locally in space and/or 
time that are constrained
 by the macroscopic solution. If the microscale problem is denoted by 
$m(\psi, d) = 0$, where the data $d$ represents the input from
 the macroscopic system, then the HMM can be succinctly expressed as 
\begin{align}\label{eq:hmm_model}
M(\Psi, D) &= 0, \qquad D = D(\psi) \nonumber \\
m(\psi, d) &= 0, \qquad d = d(\Psi).
\end{align}
With a macroscopic solver in hand, the procedure is to first
 constrain the micro simulation to be consistent with local macro data:
 $d = d(\Psi)$. After solving for $\psi$ in the micro domain,
 the missing macro data is estimated using the results from the micro simulation: $D = D(\psi)$. 

%As mentioned in the introduction, 
The HMM framework has been utilized to compute effective boundary conditions for    % new stuff
fluid simulation problems before. For instance, in \cite{hmm_microfluidics,E_contact_lines},
 the authors model fluid-fluid and fluid-solid interactions in which 
the standard no-slip boundary conditions for a continuum fluid are no longer accurate
 and must be inferred from microscopic models, such as molecular dynamics (MD).
 Using such a microscopic model throughout the entire computational domain
 is prohibitively expensive, due to the disparate spatial and temporal
 scales between the continuum and molecular dynamics involved. Instead, local
 molecular dynamic simulations are computed only along the interfaces for
 which a boundary condition is needed. In the language of \eqref{eq:hmm_model} above, 
the macro-scale model $M$ is Navier-Stokes equations, but the usual 
stress tensor is the missing data $D$ to be replaced by a more accurate
 model coming from MD simulation where it is needed. The MD simulation 
is initialized to be consistent with the local values of the continuum velocity. 

While not initially proposed as an example of HMM, the method of
 Superparameterization proposed by Grabowski \cite{grabowski_1}
 and developed by Majda and others \cite{majda_1,majda_2,Grooms4464,MAJDA201460} 
is a multiscale method for the simulation of atmospheric flows that fits
 into the framework of HMM. The original idea of the method is to couple
 local computations for the turbulent transport quantities to a global 
macroscopic model for the atmosphere. The local computations impose
 artificial scale separation in both space and time between the 
large scale energetic motions and the small scale fluctuations and hence 
allow for a reduced computational cost. 

%We mention also recent work for determining the effective boundary 
%condition at the interface between a free fluid and a porous medium 
%\cite{Lacis:2017}. The method fits in the HMM framework and is 
%quite similar to the one proposed
%in this work; the coefficients in a generalized
%Beavers-Joseph law are determined by solving Stokes problems in a 
%microscale domain containing a unit cell of the porous media. 

In the present setting of viscous laminar flow over a 
rough boundary, the macroscopic model $M$ consists of 
the Navier-Stokes equations posed in domain with a smooth 
boundary. The missing data $D$ necessary for the model 
to be complete is the coefficient in the wall-law \eqref{eq:ex_wall_law} 
coming from the homogenization theory. 
The microscopic model $m$ consists again of the Navier-Stokes 
equations, this time posed on a single ``element" of 
roughness whose size is finite in the wall-normal direction. 
The constraint $d$ is that the values of the microscopic 
solution variables $\psi$ at the computational boundaries 
(those that are \textit{not} the rough wall, where the no-slip 
condition is prescribed) must be consistent with the 
local values of the macroscopic flow variables $\Psi$ 
at those locations. Once the microscopic problem is 
solved, the solution $\psi$ is suitably averaged to 
estimate the slip amount $D$; in this way the models 
are formally coupled. This coupling is more fully detailed 
and analyzed below in Section \ref{sec:algorithm_for_laminar_flow}.

\section{HMM for viscous laminar flow over a rough boundary}
\label{sec:algorithm_for_laminar_flow}
%%%%%%%%%%%%%%%% 
%%%%%%%
%%%%%%%
%%%%%%%
%%%%%%%%%%%%%%%% 

After some preliminary definitions, 
we now describe a heterogeneous multiscale method (HMM) 
for the efficient computation of the effective boundary condition, or wall-law, 
for the case of laminar flow over a rough surface. 
We then discuss the details of its practical implementation and 
analyze some of its properties. 
 
%\blue{
%Under the same assumptions as the asymptotic analysis described in \cref{subsec:analysis_of_laminar_flow}, 
%the macroscopic method is shown to converge to the homogenized 
%solution from \cite{pironneau} in \cref{subsec:convergence_theory}. 
%\red{[update later]} As discussed below, in practice the HMM technique iteratively produces 
%slip amounts for the macroscopic flow problem. In \cref{subsec:contraction_map} a contraction map argument 
%for the case of Stokes flow in a channel configuration
%demonstrates the convergence of the this iterative procedure to 
%a positive slip amount, which ensures the macroscopic problem is coercive.\red{[update later]}}
% Assuming the 
%Analysis of the algorithm in the setting of periodic roughness 
%then demonstrates the method's convergence to the homogenized solution
%from \cite{pironneau} described in \cref{subsec:analysis_of_laminar_flow}. 

\subsection{Preliminary definitions.}
Consider a translation of the domains
$\Omega^{\epsilon}$ and $\Omega^{0}_{\epsilon H}$ defined
in Section \ref{subsec:definitions} and \eqref{eq:Omega_epsH} by $\epsilon H$ units 
in the negative $x_2$ direction, where $H = \norm{\varphi}_{\infty}$ as before, 
so that 
\begin{equation*}
(x_1, x_2) \mapsto (x_1, x_2 -\epsilon H);
\end{equation*}
note that $\Omega^{0}_{\epsilon H} \subset \Omega^{\epsilon}$ still of course holds
after the translation.
Define $\Omegamac$ to be the resulting translation of $\Omega_{\epsilon H}^0$, 
and for simplicity continue to refer to the translation of $\Omega^{\epsilon}$ as 
$\Omega^{\epsilon}$ (and similarly for the rough boundary $\Gamma^{\epsilon}$). 
In addition, rename $\Gamma_{\epsilon H}^0$--the flat part of the boundary 
of $\Omegamac$ defined by \eqref{eq:Omega_epsH}--to be simply $\Gamma$.

Consider also a collection of points $\{s_1,s_2,\ldots,s_J\}$, each
$s_j\in\mathbb{R}$, and assume $|s_j - s_{j+1}| \ge \epsilon$ for each $j$. 
Let $\{ L_1^{\rm mic}, L_2^{\rm mic}, \ldots , L_J^{\rm mic}\}$ be a collection of
positive real values. 
Define the micro-domains $\Omegamic_j$ to be the domains bounded by the curves 
$x_1 = s_j$ on the left, $x_1 = s_j+L_j^{\rm mic}$ on the right, $x_2 = \gamma>0$ above, and
$\{(x_1,x_2)|x_2 = \zeta^{\epsilon}(x_1)-\epsilon H\}$ below. The lower
curve is simply the portion of $\Gamma^{\epsilon}$ from $x_1=s_j$ to $s_j + L_j^{\rm mic}$. 
Note that of course $\zeta^{\epsilon}$ need not agree at the locations $s_j$ and $s_j + L_j^{\rm mic}$.
Denote this portion of the micro-domains $\partial \Omegamic_{j,\rm noslip}$, as this is where
the physical wall is located. Denote the remaining portion of the boundary 
$\partial \Omegamic_{j,\rm D} = \partial \Omegamic_j \setminus \partial \Omegamic_{j,\rm noslip}$, 
since the flow will generally satisfy some Dirichlet condition there. 
Lastly, assume $\gamma$ and each $L^{\rm mic}_j$ are $\bigoh{\epsilon}$. See \ref{fig:macro_micro_cartoon} for 
an example of such a configuration. 

\smallskip
\begin{remark}
\emph{
In general, the HMM algorithm described below is not limited to rough domains with boundaries parameterized
by oscillatory functions of the form $\zeta^{\epsilon}(x_1) = \beta(x_1) \varphi(x_1/\epsilon)$, i.e.\
rapidly oscillating periodic functions modulated by a slowly varying smooth function, 
as assumed in Section \ref{subsec:definitions}. }
\end{remark}

\smallskip
\begin{remark}
\emph{
For the micro-domain problems defined below to be well-posed, 
the corners of the $\Omegamic_j$ domains should be mollified; such technical details 
are not considered here.  }
\end{remark}

\begin{figure}[tbhp]
\centering
\subfloat[
]
{
%\begin{tikzpicture}[xscale=6.5,yscale=3.4]
\begin{tikzpicture}[xscale=6.0,yscale=3.0]
\draw[dotted]  (0,1) -- (0,0);
%\draw[dotted] (0,0) -- (1,0) node[right]{$\Gamma$} ;
\draw[dotted] (1,0) -- (1,1);
\draw[dotted] (0,1) -- (1,1);
\draw[densely dotted, domain = 0:0.05] plot(\x, {0.025*(cos(40*pi*\x r) -1)*(0.5+(sin(sqrt(2)*2*pi*\x r))*(sin(sqrt(2)*2*pi*\x r)))});

\draw[densely dotted, domain = 0.05:0.1] plot(\x, {0.025*(cos(40*pi*\x r) -1)*(0.5+(sin(sqrt(2)*2*pi*\x r))*(sin(sqrt(2)*2*pi*\x r)))});
%\draw[black, domain = 0.05:0.1] plot(\x, {0.025*(cos(40*pi*\x r) -1)*(0.5+(sin(sqrt(2)*2*pi*\x r))*(sin(sqrt(2)*2*pi*\x r)))});
%\draw[black](0.05, 0) -- (0.05, 0.125); 
%\draw[black](0.05,0.125) -- (0.1, 0.125);
%\draw[black](0.1, 0.125) -- (0.1, 0);
%\draw (0.05,0.125)node[above]{$s_1$}; 

\draw[densely dotted, domain = 0.1:0.2] plot(\x, {0.025*(cos(40*pi*\x r) -1)*(0.5+(sin(sqrt(2)*2*pi*\x r))*(sin(sqrt(2)*2*pi*\x r)))});
%\draw[densely dotted, domain = 0.2:0.25] plot(\x, {0.025*(cos(40*pi*\x r) -1)*(0.5+(sin(sqrt(2)*2*pi*\x r))*(sin(sqrt(2)*2*pi*\x r)))});
\draw[black, domain = 0.2:0.25] plot(\x, {0.025*(cos(40*pi*\x r) -1)*(0.5+(sin(sqrt(2)*2*pi*\x r))*(sin(sqrt(2)*2*pi*\x r)))});
\draw[black, domain = 0.25:0.3] plot(\x, {0.025*(cos(40*pi*\x r) -1)*(0.5+(sin(sqrt(2)*2*pi*\x r))*(sin(sqrt(2)*2*pi*\x r)))});
%\draw[black](0.25, 0) -- (0.25, 0.125); 
\draw[black](0.20, 0) -- (0.20, 0.125); 
%\draw[black](0.25,0.125) -- (0.3, 0.125);
\draw[black](0.20,0.125) -- (0.35, 0.125);
\draw[black](0.35, 0.125) -- (0.35, 0);
%\draw[black](0.3, 0.125) -- (0.3, 0);
\draw (0.275,0.125)node[above]{$s_1$}; 
%\draw (0.25,0.125)node[above]{$s_2$}; 

\draw[black, domain = 0.3:0.35] plot(\x, {0.025*(cos(40*pi*\x r) -1)*(0.5+(sin(sqrt(2)*2*pi*\x r))*(sin(sqrt(2)*2*pi*\x r)))});
\draw[densely dotted, domain = 0.35:0.4] plot(\x, {0.025*(cos(40*pi*\x r) -1)*(0.5+(sin(sqrt(2)*2*pi*\x r))*(sin(sqrt(2)*2*pi*\x r)))});
%\draw[densely dotted, domain = 0.3:0.4] plot(\x, {0.025*(cos(40*pi*\x r) -1)*(0.5+(sin(sqrt(2)*2*pi*\x r))*(sin(sqrt(2)*2*pi*\x r)))});
\draw[densely dotted, domain = 0.4:0.45] plot(\x, {0.025*(cos(40*pi*\x r) -1)*(0.5+(sin(sqrt(2)*2*pi*\x r))*(sin(sqrt(2)*2*pi*\x r)))});
\draw[densely dotted, domain = 0.45:0.5] plot(\x, {0.025*(cos(40*pi*\x r) -1)*(0.5+(sin(sqrt(2)*2*pi*\x r))*(sin(sqrt(2)*2*pi*\x r)))});
%\draw[black, domain = 0.45:0.5] plot(\x, {0.025*(cos(40*pi*\x r) -1)*(0.5+(sin(sqrt(2)*2*pi*\x r))*(sin(sqrt(2)*2*pi*\x r)))});
%\draw[black](0.45, 0) -- (0.45, 0.125); 
%\draw[black](0.45,0.125) -- (0.5, 0.125);
%\draw[black](0.5, 0.125) -- (0.5, 0);
%\draw (0.45,0.125)node[above]{$s_3$}; 

\draw[densely dotted, domain = 0.5:0.6] plot(\x, {0.025*(cos(40*pi*\x r) -1)*(0.5+(sin(sqrt(2)*2*pi*\x r))*(sin(sqrt(2)*2*pi*\x r)))});
\draw[densely dotted, domain = 0.6:0.65] plot(\x, {0.025*(cos(40*pi*\x r) -1)*(0.5+(sin(sqrt(2)*2*pi*\x r))*(sin(sqrt(2)*2*pi*\x r)))});
\draw[densely dotted, domain = 0.65:0.7] plot(\x, {0.025*(cos(40*pi*\x r) -1)*(0.5+(sin(sqrt(2)*2*pi*\x r))*(sin(sqrt(2)*2*pi*\x r)))});
%\draw[black, domain = 0.65:0.7] plot(\x, {0.025*(cos(40*pi*\x r) -1)*(0.5+(sin(sqrt(2)*2*pi*\x r))*(sin(sqrt(2)*2*pi*\x r)))});
%\draw[black](0.65, 0) -- (0.65, 0.125); 
%\draw[black](0.65,0.125) -- (0.7, 0.125);
%\draw[black](0.7, 0.125) -- (0.7, 0);
%\draw (0.65,0.125)node[above]{$s_4$}; 

\draw[densely dotted, domain = 0.7:0.75] plot(\x, {0.025*(cos(40*pi*\x r) -1)*(0.5+(sin(sqrt(2)*2*pi*\x r))*(sin(sqrt(2)*2*pi*\x r)))});
%\draw[densely dotted, domain = 0.75:0.8] plot(\x, {0.025*(cos(40*pi*\x r) -1)*(0.5+(sin(sqrt(2)*2*pi*\x r))*(sin(sqrt(2)*2*pi*\x r)))});
%\draw[densely dotted, domain = 0.8:0.85] plot(\x, {0.025*(cos(40*pi*\x r) -1)*(0.5+(sin(sqrt(2)*2*pi*\x r))*(sin(sqrt(2)*2*pi*\x r)))});
\draw[black, domain = 0.75:0.8] plot(\x, {0.025*(cos(40*pi*\x r) -1)*(0.5+(sin(sqrt(2)*2*pi*\x r))*(sin(sqrt(2)*2*pi*\x r)))});
\draw[black, domain = 0.8:0.85] plot(\x, {0.025*(cos(40*pi*\x r) -1)*(0.5+(sin(sqrt(2)*2*pi*\x r))*(sin(sqrt(2)*2*pi*\x r)))});
\draw[black, domain = 0.85:0.9] plot(\x, {0.025*(cos(40*pi*\x r) -1)*(0.5+(sin(sqrt(2)*2*pi*\x r))*(sin(sqrt(2)*2*pi*\x r)))});
\draw[black](0.75, 0) -- (0.75, 0.125); 
%\draw[black](0.85, 0) -- (0.85, 0.125); 
\draw[black](0.75,0.125) -- (0.9, 0.125);
%\draw[black](0.85,0.125) -- (0.9, 0.125);
\draw[black](0.9, 0.125) -- (0.9, 0);
\draw (0.825,0.125)node[above]{$s_2$}; 
%\draw (0.85,0.125)node[above]{$s_5$}; 

\draw[densely dotted, domain = 0.9:1.0] plot(\x, {0.025*(cos(40*pi*\x r) -1)*(0.5+(sin(sqrt(2)*2*pi*\x r))*(sin(sqrt(2)*2*pi*\x r)))});
\end{tikzpicture}
}
\subfloat[
]
{
%\begin{tikzpicture}[xscale = 10.875, yscale = 6.3] %[xscale=6.5,yscale=3.8]
\begin{tikzpicture}[xscale = 9.875, yscale = 5.5] %[xscale=6.5,yscale=3.8]
\draw[black, domain = 0.2:0.3] plot(\x, {0.05*(cos(20*pi*\x r) -1)});
\draw[black, domain = 0.3:0.4] plot(\x, {0.07*(cos(20*pi*\x r) -1)});
\draw[black, domain = 0.4:0.5] plot(\x, {0.03*(cos(20*pi*\x r) -1)});
\draw[black] (0.2,0) -- (0.2,0.35);
\draw[black] (0.2,0.25);
%\draw[black] (0.3,0.35) -- (0.3, 0);
\draw[black] (0.5,0.35) -- (0.5, 0);
\draw[black] (0.2,0.35)node[left]{$(s_j,\gamma)$} -- (0.5,0.35); 
\draw[fill=black] (0.2,0.35) ellipse (0.1pt and 0.2pt); 
%\draw[black] (0.2,0.35)node[left]{$(s_j,\gamma)$} -- (0.3,0.35); 
%\draw[black,densely dotted] (0.2,0.00) -- (0.3, 0)node[right]{$(s_j+\epsilon,0)$};
\draw[black,densely dotted] (0.2,0.00) -- (0.5, 0)node[right]{$(s_j+L^{\rm mic}_j,0)$};
%\draw (0.5,0.0) circle (2.0pt); 
\draw[fill=black] (0.5,0.0) ellipse (0.1pt and 0.2pt); 
\draw (0.25,-0.1)node[left]{$\partial \Omegamic_{j,\rm noslip}$}; 
\draw (0.5,0.3)node[right]{$\partial \Omegamic_{j,\rm D}$}; 
\end{tikzpicture}
}
\caption{ 
(a) Example domain $\Omega^{\epsilon}$ 
(\protect\tikz[baseline]{\protect\draw[line width=0.2mm,densely dotted] (0,.4ex)--++(0.25,0) ;})
containing several $\Omegamic_j$
(\protect\tikz[baseline]{\protect\draw[line width=0.2mm] (0,.4ex)--++(0.25,0) ;}).
(b) One instance of a micro domain $\Omegamic_j$ whose boundary consists of 
two pieces, $\partial \Omegamic_{j, \rm D}$ and $\partial \Omegamic_{j, \rm noslip}$. 
}
\label{fig:macro_micro_cartoon}
\end{figure}

\subsection{Multiscale method.}
Let $\ueps$ be a solution to the stationary Navier-Stokes equations in a
domain $\Omega^{\epsilon}$ with a rough-boundary.
The purpose of the multiscale method is to efficiently produce an approximation $U$ 
to the true, oscillatory flow $\ueps$ by enforcing that $U$ satisfy a wall-law
of the form from the homogenization theory
\begin{equation}\label{eq:mac_wall_law}
U = \alpha \pderiv{U_1}{x_2} \, e_1 
\end{equation}
on $\Gamma$, the boundary of the smooth domain $\Omegamac$. 
The coefficient $\alpha$ in the wall-law parameterizes the 
approximation and, in the language of HMM from Section \ref{subsec:hmm}, 
is the missing data $D$ needed to complete the macro-scale model
of $\ueps.$

The strategy utilized in \cite{pironneau} to determine $\alpha$ 
in the case of periodic roughness results from the asymptotic
analysis from Section \ref{subsec:analysis_of_laminar_flow}.
One simply precomputes
the solution $\chi$ to the (truncated) cell problem 
\eqref{eq:cell_prob_mod} and then 
takes the average of the horizontal component 
$\overline{\chi_1}$. After scaling by $\epsilon$, this constant (plus some
amount $\delta > \epsilon H$, in light of Theorem \ref{thm:thrm2_p}) 
is then taken to be the missing data $\alpha$. 
 The precomputing step is possible because the cell
problem depends only on the geometry of the roughness.

In contrast, the multiscale method defined below to estimate $\alpha$
generally involves coupling a Navier-Stokes system posed in the 
macroscale domain $\Omegamac$ to $J$ separate Navier-Stokes 
systems posed in microscale domains $\Omegamic_j$, $1\le j\le J$. 
Similar to the cell problem from the homogenization theory, 
the microscale systems account for the geometry of the rough surface, 
however, the current method is not restricted to domains with 
periodic roughness. The micro-systems are additionally constrained 
to match the averaged local flow values of the macroscale system, 
which could allow for a more accurate representation of the 
effect of surface roughness on the local macroscopic flow. 
In each $\Omegamic_j$, the ratio
of the average flow and the average flow gradient (the shear) is computed. 
These values are interpolated along $\Gamma$, and the resulting function is then 
used for the slip amount $\alpha$ in \eqref{eq:mac_wall_law}.
%Furthermore, the inclusion of the nonlinear term allows for                      %%%% KEEP THIS COMMENT ??? 
%convective effects not captured by the cell problem's linear Stokes system. 

%The coupled system is thus designed to reproduce the effective
%boundary condition from the homogenization theory whenever it is 
%applicable, but also to perform favorably 
% in more general situations to which the theory does not apply.

%%%
% define macro solver
%%%
We define first the macro and microscale systems, as well as
a projection, smoothing, and interpolation operator before 
fully formulating the coupled multiscale model. 

\smallskip
\begin{definition}\label{dfn:macro_prob}
[Macroscale system] Let $M\left(U,P,\alpha\right)$ define the following PDE system
posed in $\Omegamac$ and 
parameterized by the slip amount $\alpha$:
\begin{align*}\label{eq:macro_prob}
-\nu \Delta U + U \nabla U + \nabla P - f &= 0,  \qquad \text{in } \Omegamac \nonumber \\
\nabla \cdot U &= 0, \qquad \text{in } \Omegamac \nonumber \\
U - \alpha\,\, \pderiv{U_1}{x_2} \, e_1 &=0 , \qquad \text{on } \Gamma \nonumber \\
U &= 0, \qquad \text{on } \partial \Omegamac \setminus \Gamma 
\end{align*}
\end{definition}
Note that in general $\alpha$ can vary along $\Gamma$ so that $\alpha = \alpha(x_1)$.
%%%
% define micro solver
%%%
\smallskip
\begin{definition}\label{dfn:micro_prob}
[Microscale system] Let $m_j\left(u^j,p^j,\Upsilon_j\right)$ define the PDE system posed 
in $\Omegamic_j$ and parameterized by the Dirichlet boundary condition 
$\Upsilon_j: \partial \Omegamic_{j,\rm D}\to \mathbb{R}^2$:
\begin{align*} %\label{eq:micro_prob}
-\nu \Delta u^j + u^j \nabla u^j + \nabla p^j - f &= 0, \qquad \text{in } \Omegamic_j \nonumber \\
\nabla \cdot u^j &=0, \qquad \text{in } \Omegamic_j \nonumber \\
u^j &= \Upsilon_j, \,\,\,\,\,\,\,\, \text{on } \partial \Omegamic_{j,\rm D} \nonumber \\
u^j &= 0 \,\,\,\,\,\,\,\,\,\,\,\,\,\, \text{on } \partial \Omegamic_{j, \rm no slip} 
%u \text{ and its derivatives} \,\,\, &x_1 \text{ periodic in } \Omega_j^m \text{ for } x_2\ge 0
\end{align*}
\end{definition}
For well-posedness, i.e.\ conservation of mass, $\Upsilon_j$ should satisfy 
\begin{equation}\label{eq:cons_of_mass}
\int_{\partial \Omegamic_{j,\rm D}} \Upsilon_j \cdot n \, ds = 0. 
\end{equation}
Furthermore, let $\{s_n\}_{n=1}^{\infty} \in \partial \Omegamic_{j,\rm D}$ be 
some convergent sequence with limiting point $\sigma = (s_j,0)$ or $\sigma = (s_j +L_j^{\rm mic},0)$. 
Then $\Upsilon_j$ should also satisfy 
\begin{equation}\label{eq:micro_cont}
\lim_{n\to \infty} \Upsilon_j(s_n) = 0
\end{equation}
for consistency with the no slip condition posed along $\partial \Omegamic_{j, \rm noslip}$.
%%%
% define projection operator
%%%

\smallskip
\begin{definition}\label{dfn:projection_op}
[Projection operator] 
Let $\mathcal{P}_j$ be the collection of continuous maps
from $C(\Omegamac,\mathbb{R}^2)$ to $C\left(\partial \Omegamic_{j, \rm D},\mathbb{R}^2\right)$, 
where $C(X,Y)$ denotes the set of continuous functions from $X$ to $Y$.
Then we say $\pi_j$ is a projection operator if 
$\pi_j \in \mathcal{P}_j$ and $\pi_j(f)$ satisfies the properties 
\eqref{eq:cons_of_mass} and \eqref{eq:micro_cont} for any 
$f \in C(\Omegamac,\mathbb{R}^2)$. 
\end{definition}

This projection operator is the mechanism by which the micro-problems $m_j$ are 
constrained to match the macroscopic solution $U$. 
Suppose that $U$ is a continuous solution to $M(U,P,\alpha)$ for $\alpha\ne 0$. 
Note then that simply taking
the trace of $U$ along $\partial \Omegamic_{j, \rm D}$ is not sufficient to
be a projection operator in the above sense;  
even though the conservation of mass property \eqref{eq:cons_of_mass} holds, 
the constraint \eqref{eq:micro_cont} will not because of the slip condition on $\Gamma$. 
A specific example of a projection map is mentioned in Remark \ref{rmk:bc_micro} and 
more fully detailed in Appendix \ref{appendix:micro_bc}.
%%%
% define smoothing operator
%%%

\smallskip
\begin{definition}\label{dfn:smoothing_op}
[Smoothing operator] For continuous and integrable $\psi:\mathbb{R}^2\to \mathbb{R}$ and 
$y \in \mathbb{R}$, define the operator
\begin{equation*}
\chevron{\psi}_L(x,y) := \int_{x}^{x+L} \psi(s,y) \, ds
\end{equation*}
which integrates $\psi$ in the horizontal direction along a $L$-sized 
strip $[x,x+L]$ at fixed height $y$. 
\end{definition}
%%%
% define interpolation operator
%%%

\smallskip
\begin{definition}\label{dfn:interp_op}
[Interpolation operator] For $\{(s_j,\alpha_j)\}_{j=1}^J \in \mathbb{R}^{2\times J}$, define
\begin{equation*}
\mathcal{I}((s_1,\alpha_1),\ldots,(s_J,\alpha_J)):\mathbb{R}^{2\times J} \to C\left(\Gamma\right)
\end{equation*}
denote a piecewise continuous polynomial interpolant based on the given points.
\end{definition}
%%%
% formally define the coupled model
%%%

Using the above definitions, the HMM is formally defined as follows; 
given a roughness profile $\zeta^{\epsilon}$, macroscale domain $\Omegamac$,
and a collection of points $\{s_j,L^{\rm mic}_j\}_{j=1}^J$ with associated 
micro-domains $\Omegamic_j$, find 
$\left((U,P),(u^1,p^1), \ldots, (u^J, p^J)\right)$ satisfying the 
coupled system of equations:
\begin{align}
M\left( U,P, \alpha \right)  = 0& \nonumber \\
m_j\left( u^j,p^j, \Upsilon_j \right) = 0& \qquad 1\le j\le J  \label{eq:hmm_system}
\end{align} 
where
\begin{align*}
\Upsilon_j = \,\,\pi_j(U)& \qquad 1\le j\le J,  \\
%\alpha_j &= \frac{\chevron{u^j_1}(s_j,0)}{  \chevron{\partial u^j_1/\partial x_2}(s_j,0)} \label{eq:hmm4} \\
%\alpha_j = \chevron{u^j_1}_{L_j^{\rm mic}}(s_j,0) \Big(
%\chevron{\partial u^j_1/\partial x_2}_{L_j^{\rm mic}}(s_j,0)\Big)^{-1}&, \qquad 1\le j\le J \label{eq:hmm4} \\
\alpha = \mathcal{I}((s_1,\alpha_1), \ldots,(s_J,\alpha_J))&, 
\end{align*}
and
\begin{equation}
\alpha_j = \frac{ \langle u^j_1 \rangle_{L_j^{\rm mic}}(s_j,0)}{ 
\langle \partial u^j_1/\partial x_2 \rangle_{L_j^{\rm mic}}(s_j,0)}, \qquad 1\le j\le J. \label{eq:hmm5} 
\end{equation}
The HMM thus consists of a stationary Navier-Stokes equation in $\Omegamac$ with 
slip amount $\alpha$ and $J$ stationary Navier-Stokes equations posed in the
domains $\Omegamic_j$, each of which depend on the projection of $U$ onto the 
boundaries $\Omegamic_{j,D}$. Given a boundary condition $\Upsilon_j$ and corresponding
micro-solution $u^j$, the slip amounts $\alpha_j$ are defined to be the ratio of 
the average horizontal flow velocity to the average flow derivative in the 
vertical direction. The average is taken across the length $L_j$ of the micro-domain
at $x_2 =0$, i.e.\ just above the roughness. In our numerical tests, we also considered 
computing the $\alpha_j$ values at some $x_2 \in (0,\gamma)$ and extrapolating the result
back to $x_2= 0$ for use in the macro-solver; see the discussion in Section \ref{subsec:discussion} for the
details. The values $\{\alpha_j\}_{j=1}^J$ are patched together with some interpolation
scheme $\mathcal{I}$ and in turn utilized by the macro-solver. 

%%%%%%%%%%%%%%%%%%%%%%%%%%%%% remark on three 3 case
\smallskip
\begin{rem}[Three-dimensional case]
\label{rmk:3d_case}
\emph{
As mentioned in Section \ref{sec:background}, the slip amount $\alpha$ in 
the wall-law \eqref{eq:ex_wall_law}
can generally be a tensor in the case of three-dimensional flow. Suppose
$(x_1,x_3)$ and $x_2$ denote the horizontal and vertical directions, 
respectively. Then the wall-law \eqref{eq:ex_wall_law} for the macro-solver would be
\begin{equation*}
U_i = \sum_{k \in \{1,3\}} \alpha_{ik} \pderiv{U_k}{x_2}, \qquad i = 1,3.
\end{equation*}
Each individual slip amount then is
\begin{equation*}
 (\alpha_j)_{ik} = \frac{ \langle u^j_i \rangle_{L_j^{\rm mic}}(s_j,0)}{ 
\langle \partial u^j_k/\partial x_2 \rangle_{L_j^{\rm mic}}(s_j,0)}, 
\qquad i,k \in \{1,3\},
\end{equation*}
where the index $j$ runs over each micro-domain. Each entry is interpolated
across the two-dimensional surface $\Gamma$ to produce a slip amount 
$\alpha$ for the macro-solver.}
\end{rem}
%%%%%%%%%%%%%%%%%%%%%%%%%%%%% END remark on three 3 case

%%%%%%%%%%%%%%%%%%%%%%%%%%%%% remark on BCs 
\smallskip
\begin{remark}[Boundary conditions for the micro-systems]
\label{rmk:bc_micro}
\emph{
Given some macroscopic flow $U$, each of the microscopic problems
depends on the boundary condition from the projection operator:
\begin{equation*}
u^j = \Upsilon_j = \pi_j(U) \qquad \text{on } \partial \Omegamic_{j,\rm D}.
\end{equation*}
In the case when $U$ is horizontal, i.\,e.\  the vertical component
of the velocity vector is zero, or at least negligible compared
with the horizontal component, then the boundary conditions for the micro-systems
 can be simplied to a ``free stream" condition along the upper computational boundary $x_2 = \gamma$:}
\begin{equation}\label{eq:free_stream}
u^j %=  \frac{1}{L_j^{\rm mic}}\int_{s_j}^{s_j + L_j^{\rm mic}} U(x, \gamma) dx 
= \frac{1}{L_j^{\rm mic}} \chevron{U_1}_{L_j^{\rm mic}}(s_j,\gamma)\, \, e_1
\end{equation}
\emph{and periodic boundary conditions at $x_1=s_j$ and $x_1 = s_j+L_j^{\rm mic}$. 
In this case the projection operator $\pi_j$ simply maps from}
$C(\Omegamac,\mathbb{R}^2)$ \emph{to} $\mathbb{R}$. \emph{Note that this approach 
is only practically feasible if the roughness function $\zeta^{\epsilon}$ 
satisfies $\zeta^{\epsilon}(s_j) \approx \zeta^{\epsilon}(s_j + L_j^{\rm mic})$
so that a periodic mesh for the micro domain $\Omegamic_j$ can be constructed without much error.}
\emph{In situations for which this is not true or when the macroscopic 
flow has a nontrivial vertical component, a more general approach is 
needed. We propose quadratic polynomial Dirichlet conditions for
both the horizontal and vertical components of the velocity along
each of the three faces $\partial \Omegamic_{j,D}$. Enforcing the
no-slip condition, some interpolation constraints, and that 
the quadratic profiles preserve the macro-scale mass fluxes
guarantees a unique, well-defined projection operator. The details
can be found in Appendix \ref{appendix:micro_bc}. }
\end{remark}
%%%%%%%%%%%%%%%%%%%%%%%%%%%%% END remark on BCs

%%%%%%%%%%%%%%%%%%%%%%%%%%%%% remark on locations of s_js 
\smallskip
One key feature of the method clearly is the specification of the locations $\{s_j\}_{j=1}^J$ of the micro-domains
$\Omegamic_j$ and the domain lengths $\{L_j^{\rm mic}\}_{j=1}^J$. 
%A general strategy is to select a micro-domain location $s_j$ where either the 
%roughness varies nontrivially or the macroscopic flow is qualitatively different, or both. 
For Poiseuille 
type channel flow with periodic roughness, for instance, only one micro-domain covering a single
periodic roughness element is necessary. In more realistic settings for which the surface roughness
is nonperiodic and additionally varies over macroscopic length scales, the micro-domain lengths 
$L_j^{\rm mic}$ should be chosen large enough to cover a few of the estimated correlation lengths, 
or approximate periods. 
The $s_j$ locations should also be placed frequently enough along $\Gamma^{\epsilon}$
 to capture its large-scale, macroscopic variations. 
The numerical examples in Section \ref{sec:laminar_numerics} are chosen to approximate such situations. 
%%%%%%%%%%%%%%%%%%%%%%%%%%%%% END remark on loc of s_js

\subsection{Simulation algorithm.}\label{subsec:algo}
In Section \ref{subsec:analysis}, we use the contraction mapping principle 
to prove that the coupled, stationary system 
\eqref{eq:hmm_system} has a unique solution 
for the simplified case of linear, Stokes flow in a 
channel with periodic roughness. In practice, the 
coupled system is also solved iteratively starting from 
some initial guess for $\alpha$. A simple choice is
to use the no-slip condition along $\Gamma$, i.e.\ 
use $\alpha = 0$. One then solves the macroscale
problem, estimates the boundary conditions for the 
microscale problems, computes their solution, and 
estimates an updated slip amount. This process is 
repeated until the difference between subsequent 
slip amounts is smaller than some prescribed tolerance.

More precisely, let $\tau>0$ be some fixed tolerance,
and let the macro/microscale domains $\Omegamac$ and 
$\{\Omegamic_j\}_{j=1}^J$ be given. Then: 
\begin{enumerate}
\item[(1)] Let $\underline{\alpha}=0$ and solve the macroscale
problem 
$$M(U, P, \underline{\alpha}) = 0 \qquad \text{for } (U,P).$$
 
\item[(2)] For each $j= 1, \ldots, J$ determine boundary 
conditions for the micro-domain $\Upsilon_j = \pi_j(U)$
and then solve the microscale problems 
$$ m(u^j, p^j, \Upsilon_j) = 0 \qquad \text{for } (u^j,p^j).$$

\item[(3)] For each $j = 1, \ldots, J$, estimate the local 
slip amounts
$$
\alpha_j = \frac{ \langle u^j_1 \rangle_{L_j^{\rm mic}}(s_j,0)}{ 
\langle \partial u^j_1/\partial x_2 \rangle_{L_j^{\rm mic}}(s_j,0)}
$$
and then interpolate the values to define the a new slip 
amount
$$
\alpha_{(0)} = \mathcal{I}((s_1,\alpha_1), \ldots,(s_J,\alpha_J)).
$$

\item[(4)] If $\norm{\alpha_{(0)} - \underline{\alpha}}_{\infty} < \tau$, compute
a final macroscale solution $(U,P)$ by solving
\begin{equation}\label{eq:macro_alpha1}
 M(U,P, \alpha_{(0)}) = 0.
\end{equation}
Otherwise, use $(U,P)$ from \eqref{eq:macro_alpha1}, repeat steps (2) 
and (3) to produce a new $\alpha_{(1)}$, and compare $\alpha_{(1)}$ 
and $\alpha_{(0)}$. Repeat until successive slip amounts differ
by less than the tolerance $\tau$. 
\end{enumerate}

\smallskip
If the tolerance is set to $\bigoh{\epsilon^2}$, 
we show in Theorem \ref{thm:rapid_convergence} that for the case of linear, 
Stokes flow in a channel with periodic roughness,
the back-and-forth procedure just outlined will 
terminate in a single iteration. 
More precisely, 
we show $(\alpha_{(1)}-\alpha_{(0)})/\epsilon^2$
vanishes as $\epsilon \searrow 0$. For this case,
the HMM hence requires only one solution 
of the microscale systems in practice; computing
an updated slip amount $\alpha_{(2)}$ is not necessary. 
Although we have no proof for more general situations, we
observe the same behavior in all the numerical 
experiments presented in Section \ref{sec:laminar_numerics}.
Despite the rapid convergence, the need to have some 
initial macro-flow $U$ in hand from
which to estimate boundary conditions 
$\Upsilon = \pi(U)$ for the micro-solvers
is a disadvantage of the proposed HMM. In practice,
the overhead cost to produce some initial $\Upsilon$ 
could for example be reduced by first computing 
$U$ on a mesh that is relatively coarse. 
For time-dependent problems, the iterative procedure 
could more naturally be incorporated into a time-marching scheme, 
such as a predictor-corrector method. 

Finally, we note that since the microscopic systems $m_j(u^j,p^j,\Upsilon_j)$ are all 
independent of one another, they are trivially parallelizable.

\subsection{Analysis of the the coupled system.}
\label{subsec:analysis}
We now prove some convergence results for our HMM 
 in a simplified setting. 
Let $0 < \gamma_1 < \gamma_2 < 1$ and consider a 
macroscale channel configuration, so that 
$\Omegamac = [0,l]\times [\gamma_1,1]$, where 
$0 < \gamma_1 < 1$ and $l >0$. 
We assume the full channel with rough boundary 
 that contains $\Omegamac$ has periodic 
roughness with a maximum at $x_2 =0$ as before. 
For technical reasons however we instead measure 
the slip amount in the micro-domain at $x_2 = \gamma_1 >0$, 
and hence impose the wall-law at this location in 
the macro-domain as well. This is consistent with 
Theorem \ref{thm:thrm2_p} 
from Achdou et al.\ \cite{pironneau}. The micro-domain
has height $\gamma_2$.

Let the fluid viscosity $\nu = 1$ and assume both 
the macro and micro problems are driven by a 
constant forcing $f = -2 e_1$ representing the 
macroscale pressure gradient forcing the flow. 

%%%%%%%% define macro --> micro map
For some fixed slip amount $\alpha > 0$, a modified
Poiseuille flow solves the stationary Navier-Stokes 
equations posed in $\Omegamac$ with the no-slip 
condition $U = 0$ at $y = 1$ and the slip-condition 
$U = \alpha\,\, \partial U_1/\partial x_2 \, e_1$ 
at $y=\gamma_1$: 
\begin{equation*}\label{eq:outersoln}
U(x_2)  = \left(-x_2^2 
+ \left(\frac{\gamma_1^2 -1 - 2\gamma_1 \alpha}{\gamma_1 - 1-\alpha}\right)(x_2 - 1) 
+ 1\right)  e_1.
\end{equation*}
Since the macro-flow is horizontal and the roughness
is assumed to be periodic, the boundary condition for
the micro-domain can be simplified to a free-stream 
condition at $x_2 = \gamma_2$. We can thus define the 
map from the macro to microscale flow. 

\smallskip
\begin{definition}\label{dfn:macro-micro-map}
[Macro $\longrightarrow$ micro map]
Let $\alpha \ge 0$ and define 
\begin{equation} \label{eq:macromap}
T_2(\alpha) := -\gamma_2^2 
+ \left(\frac{\gamma_1^2 -1 - 2\gamma_1 \alpha}{\gamma_1 - 1-\alpha}\right)(\gamma_2 - 1)
+ 1 .
\end{equation}
\end{definition}

\smallskip
The derivative is 
\begin{equation}\label{eq:dUdalpha}
\pderiv{T_2}{\alpha} = \left(1-\gamma_2 \right)\left(\frac{(\gamma_1-1)^2}{(\gamma_1-1-\alpha)^2} \right), 
\end{equation}
which is positive $\forall \alpha \ge 0$ and is guaranteed to be finite whenever $\alpha >0$. 

%%%%%%%% define micro --> macro map
Instead of the full nonlinear problem, we now consider a Stokes 
problem in the micro-domain $\Omegamic$
\begin{align}\label{eq:micro_problem}
-\Delta u + \nabla p = 2 e_1 & \nonumber \\ 
\nabla \cdot u = 0 &
\end{align}
with the boundary condtions $u = 0$ 
along $\Gamma_{\epsilon}$ and $u = \U e_1$ at $x_2 = \gamma_2$ for some 
$\U\in \mathbb{R}$. 
At the side-walls $\{x_1=0\}\times[0,\gamma_2]$ and $\{x_1=\epsilon\}\times [0,\gamma_2]$, 
$(u,p)$ is periodic. 
Since all of the problem data is smooth, it is known
\cite{temam} that $u$ and $p$ are smooth fields. 
\begin{figure}[h!]  %%%%%%%%%%%%%%%%%%% BEGIN FIGURE FOR MICRO DECOMP
\centering
\begin{tikzpicture}[xscale = 10., yscale = 10.] %[xscale=6.5,yscale=3.8]

%u1 domain 
\draw[dash dot, domain = 0.0:0.2] plot(\x, {0.05*(cos(10*pi*\x r) -1)});
\draw[black] (0.0,0) -- (0.0,0.35)node[left]{$x_2 = \gamma_2$}; %node[left]{$\vec{u}_1 = 0$};
\draw[black] (0.0,-0.10)node[left]{$x_2 = -M$} -- (0.2,-0.10); %node[right]{$y=-M$};
\draw[black] (0.0,-0.10) -- (0.0,0.00);
\draw[black] (0.2,-0.10) -- (0.2,0.00);
\draw[black,densely dotted] (0.0,0.18)node[left]{$x_2=\gamma_1$}-- (0.2, 0.18) ;
\draw (0.1,0.35)node[above]{$u_{1}=0$};
\draw (0.1,-0.10)node[below]{$u_{1}=0$};
%\draw (0.1,-0.10)node[below]{$\Gamma_{\epsilon}$};
\draw[black] (0.2,0.35) -- (0.2, 0);
\draw[black] (0.0,0.35) -- (0.2, 0.35); %node[right]{$y=\gamma_2$};
%\draw[black,densely dotted] (0.0,0.00)-- (0.2, 0)node[right]{$y=0$} ;
%\draw[black,densely dotted] (0.0,0.10)-- (0.2, 0.1)node[right]{$y=\gamma_1$} ;

%u2 domain 
\draw[black, domain = 0.4:0.6] plot(\x, {0.05*(cos(10*pi*(\x -0.4 )  r) -1)});
\draw[black] (0.4,0) -- (0.4,0.35); %node[left]{$\vec{u}_2 = 0$};
%\draw[dash dot] (0.4,-0.10) -- (0.6,-0.10)node[right]{$y=-M$};
%\draw[dash dot] (0.4,-0.10) -- (0.4,0.00);
%\draw[dash dot] (0.6,-0.10) -- (0.6,0.00);
\draw (0.5,0.35)node[above]{$u_{2}=0$};
\draw (0.5,0.00)node[above]{$B_{2}(x)$};
\draw (0.5,-0.10)node[below]{$u_2 = -u_1 \Big\vert_{\Gamma_{\epsilon}}$};
\draw[black] (0.6,0.35) -- (0.6, 0);
\draw[black] (0.4,0.35) -- (0.6, 0.35); %node[right]{$y=\gamma_2$};
\draw[black,densely dotted] (0.4,0.00)-- (0.6, 0)node[right]{$x_2=0$} ;
%\draw[black,densely dotted] (0.4,0.10)-- (0.6, 0.1)node[right]{$y=\gamma_1$} ;
%\draw (0.27,-0.1)node[right]{$\partial Y$}; 

%u3 domain 
\draw[black, domain = 0.8:1.0] plot(\x, {0.05*(cos(10*pi*(\x -0.8 )  r) -1)});
\draw[black] (0.8,0) -- (0.8,0.35); %node[left]{$\vec{u}_3 = \U \vec{e}_1$};
%\draw[dash dot] (0.4,-0.10) -- (0.6,-0.10)node[right]{$y=-M$};
%\draw[dash dot] (0.4,-0.10) -- (0.4,0.00);
%\draw[dash dot] (0.6,-0.10) -- (0.6,0.00);
\draw (0.9,0.35)node[above]{$u_{3}=\U e_1$};
\draw (0.9,0.00)node[above]{$B_{3}(x)$};
\draw (0.9,-0.10)node[below]{$u_3 = 0 $};
\draw[black] (1.0,0.35) -- (1.0, 0);
\draw[black] (0.8,0.35) -- (1.0, 0.35); %node[right]{$y=\gamma_2$};
\draw[black,densely dotted] (0.8,0.00)-- (1.0, 0)node[right]{$x_2=0$} ;

\end{tikzpicture}
\caption{Boundary conditions for each of the three pieces making up the 
micro solution $u$. Each piece is additionally periodic at $x_1=0$ and $x_1=\epsilon$.
The interior traces of the second and third piece at the line $x_2 = 0$
are marked as $B_2(x)$ and $B_3(x)$. 
}
\label{fig:three_pieces}
\end{figure}        %%%%%%%%%%%%%%%%%%% END FIGURE FOR MICRO DECOMP

The linearity of the Stokes problem \eqref{eq:micro_problem} allows 
us to decompose the micro solution into three pieces
\begin{equation}\label{eq:micro_decomp}
u = u_1\big\vert_{\Omega_{\rm mic}}+ u_2 +  u_3, 
\end{equation}
which lets us systematically quantify the HMM slip amount 
$\alpha$. The boundary conditions and domains for each piece are illustrated 
in Figure \ref{fig:three_pieces}. We now describe each 
solution $u_i$, $i=1,2,3$ in more detail.

The first piece accounts for the constant forcing $-2 e_1$ and is 
posed in the larger domain $[0, \epsilon]\times [-M, \gamma_2] 
\supset \Omega_{\rm mic}$, where $M = \norm{\varphi^{\epsilon}}_{\infty} = \epsilon H$ is the amplitude of 
the roughness. With no-slip at $x_2 = -M$ and $x_2 = \gamma_2$, it 
has analytic solution $u_1(x_2) = (x_2 + M)(\gamma_2-y) e_1$ and 
$p_1 = 0$. 

The second piece $u_2$ corrects for the fact that the restriction 
of $u_1$ to $\Omegamic$ does not satisfy the no-slip condition
along $\Gamma_{\epsilon}$. It satisfies the homogeneous Stokes
equations (i.e.\ without forcing), the no-slip condition at 
$x_2 = \gamma_2$ and 
\begin{equation*}
u_2 = -\left(u_1 \cdot e_1\right)\big\vert_{\Gamma_{\epsilon}} e_1
\end{equation*}
along $\Gamma_{\epsilon}$.
At this point and in the rest of the subsection, we 
abuse notation and refer to the \emph{horizontal} component
of the flow vector $u_i \cdot e_1$ simply as $u_i$ for each 
$i = 1,2,3$.  

Let $\chevron{\cdot}(x_2)$ be shorthand for the integral operator 
$\chevron{\cdot}_{\epsilon}(0,x_2)$ from Definition \ref{dfn:smoothing_op}, 
and let $\overline{B}_2$ denote $\chevron{u_2}(0)$. Since 
$\chevron{u_2}(x_2)$ is independent of $x_1$, it satisfies a 
simple boundary value problem in $[0, \gamma_2]$ with solution 
\begin{equation*}
\chevron{u_2}(x_2) = -(\overline{B}_2/\gamma_2) x_2 + \overline{B}_2.
\end{equation*}
In Appendix \ref{appendix:estimates} we show that in particular 
$|\overline{B}_2|/\epsilon \to 0$ as $\epsilon \searrow 0$. 

The third piece $u_3$ will account for how the
macro-solution $U$ influences the micro-solution. It satisfies the 
no-slip condition along $\Gamma_{\epsilon}$ and the free stream 
condition $u_3 = \U$ at $x_2 = \gamma_2$. Like $u_2$ it is not
externally forced, in contast to $u_1$. By linearity of the Stokes
equation, $\chevron{u_3}(x_2)$ can be written as 
$\U \chevron{\tilde u_3}(x_2)$, where $\tilde u_3$ solves the 
same Stokes problem as $u_3$ but satisfies $\tilde u_3 = 1$ at
$x_2 = \gamma_2$. Similar to $\chevron{u_2}(x_2)$, 
$\chevron{u_3}(x_2)$ can be written as 
\begin{equation*}
\chevron{u_3}(x_2) = \U \left( \frac{(1 - \tilde B_3)}{\gamma_2}\, x_2 + \tilde B_3\right)
\end{equation*}
where $\tilde B_3 = \chevron{\tilde u_3}(0)$. In Appendix \ref{appendix:estimates}
we show that $\tilde B_3 = \bigoh{\epsilon}$. 

Using the decomposition \eqref{eq:micro_decomp}, we are now able 
to define the slip amount $\alpha$ given by \eqref{eq:hmm5} in the 
present setting. The slip amount is a map from the micro to 
macroscale flow. 

\smallskip
\begin{definition}
[Micro $\longrightarrow$ macro map] Let $\U \in \mathbb{R}$ and 
define 
\begin{equation}\label{eq:micromap}
T_1(\U) :=  \frac{C_1 + C_2 \U}{C_3 + (1-\tilde B_3) \U} 
\end{equation}
where 
\begin{align}\label{eq:T1_Cs}
C_1 &= \epsilon\gamma_2(\gamma_1+ M)(\gamma_2 - \gamma_1) 
+ \overline{B}_2(\gamma_2-\gamma_1)  \nonumber \\
C_2 &= \tilde B_3(\gamma_2-\gamma_1) + \gamma_1 \nonumber \\
C_3 &= \epsilon \gamma_2 (\gamma_2-2\gamma_1-M)- \overline{B}_2, 
\end{align}
which is the ratio of the average of the horizontal component of the 
micro-flow $u$ to 
the average of its derivative in the $x_2$ direction.
\end{definition}

\smallskip
The derivative is 
\begin{equation*}\label{eq:dadU}
\pderiv{T_1}{\U} = \frac{\tilde B_{3} C_{1} + C_{2} C_{3} - C_{1}}
{{\left((1 - \tilde B_{3})\U + C_{3}\right)}^{2}}.
\end{equation*}

The idea now is to prove the convergence of our multiscale method 
in the current simplified setting by first showing that the 
back-and-forth map $T_2$ composed with $T_1$ is a contraction for 
$\epsilon$ sufficiently small. Starting then from 
$\U_0 = T_2(\alpha)\vert_{\alpha=0}$, we'll show not only that sequence $(\U_n)_n$ 
generated by the back-and-forth iteration converges, but also 
that the slip-amount resuling from the limiting point is positive. 

\smallskip
\begin{definition}
Let $T: \mathbb{R} \to \mathbb{R}$ be the nonlinear map given by the composition of 
the two maps defined by \eqref{eq:macromap} and \eqref{eq:micromap}:
\begin{equation*} \label{eq:nonlinear_map_defn}
T(\U) = T_2(T_1(\U)). 
\end{equation*}
\end{definition}

\smallskip
From the chain rule  
\begin{equation*}
\pderiv{T}{\U}(\U)   = \pderiv{T_2}{\alpha}(\alpha(\U)) \pderiv{T_1}{\U}(\U),
\end{equation*}
and after some calculation one can show 
\begin{equation}\label{eq:dT_succinct}
\pderiv{T}{\U}(\U) = \frac{K_1 (\gamma_1-1)^2 (1-\gamma_2)}{K_2 + (1+K_3)\U}
\end{equation}
where $K_1 = \tilde B_3 C_1 + C_2 C_3 - C_1$, $K_2 = C_1 + (1-\gamma_1)C_3$, and 
$K_3 = C_2 - \gamma_1 + \tilde B_3 (\gamma_1 - 1)$. 

We now prove $T$ is a contraction map for all $\epsilon$ sufficiently 
small by showing $\partial T/\partial \U$ vanishes in the limit  
$\epsilon \searrow 0$ for all $\U$ larger than $T_2(\alpha)\vert_{\alpha=0}$. 
The result relies on an assumption that $\gamma_2$ is asymptotically 
larger than both $\epsilon$ and $\gamma_1$

\smallskip
\begin{lemma}\label{lem:contraction} %%%%%%%%%%%%% LEMMA that dT/dU vanishes as eps vanishes
Let $\U_0 = T_2(0)$, and assume that $\gamma_1$ and $\gamma_2$ monotonically 
tend to zero as $\epsilon \searrow 0$ as well as the quantities 
$\gamma_1/\gamma_2$ and $\epsilon/\gamma_2$. 
Then $\forall \U \in [\U_0, \infty)$, 
\begin{equation*}
\lim_{\epsilon \downarrow 0} \Big \vert \pderiv{T}{\U} \Big \vert =  0.
\end{equation*}
\end{lemma}
\begin{proof}
\eqref{eq:dT_succinct} can be bounded as
\begin{align*}
\left|\pderiv{T}{\U}\right| &\le  \frac{|K_1|}{|K_2 + (1+K_3)\U|} \\
&= \frac{|K_1/\U|}{|K_2/\U +K_3+1|}.
\end{align*}
Since $K_3 = \tilde B_3(\gamma_2 - 1) \to 0$ as $\epsilon \searrow 0$
by \eqref{eq:B3_bound}, then the result will follow if  
\begin{equation*}
\lim_{\epsilon \downarrow 0} K_i/\U = 0 \qquad \forall \U \in [\U_0, \infty)
\end{equation*}
for $i = 1, 2$. 
Inserting the values of $C_1$, $C_2$, and $C_3$ into the expressions for each $K_i$, we see: 
\begin{equation*}
K_1/\U_0 = -{\left[\tilde B_{3} \epsilon \gamma_{1}^{2} - 2 \, \tilde B_{3} \epsilon \gamma_{1} \gamma_{2} 
+ \tilde B_{3} \epsilon \gamma_{2}^{2} - \epsilon \gamma_{1}^{2} - M \epsilon \gamma_{2} + \overline{B}_{2}\right]} \gamma_{2}/\U_0  
\end{equation*}
and
\begin{equation*}
K_2/\U_0 = \left[-\epsilon \gamma_{2}^{2} - {\left(\epsilon \gamma_{2}^{2} - \epsilon \gamma_{2}\right)} M 
+ \overline{B}_{2} {\left(\gamma_{2} - 1\right)} - {\left(\epsilon \gamma_{1}^{2} - 2 \, \epsilon \gamma_{1}\right)} \gamma_{2} \right]/\U_0.
\end{equation*}
Note that $\U_0 = \gamma_2 (1 + \gamma_1 - \gamma_2 - \gamma_1/\gamma_2)$, $M = \epsilon H$, and 
from the estimates in Appendix \ref{appendix:estimates} we know both $\overline{B}_2/\gamma_2$  
and $\tilde B_3$ vanish as $\epsilon$ tends to zero. Hence both $K_1$ and $K_2$ limit to 0 as 
$\epsilon \searrow 0$ based on the stated assumptions. Ergo, $\partial T/\partial \U$ 
also vanishes as $\epsilon$ tends to zero $\forall \U \in [\U_0, \infty)$. 
\end{proof}

\smallskip
%%%%%%%%%%%%%%%%%%%%% define sequences (U_n)_n and (alpha_n)_n
\begin{definition}\label{dfn:sequence_Un_and_alpha_n}
Let 
\begin{equation*}\label{eq:U_n_sequence}
\left(\U_0, \U_1, \ldots, \U_n, \ldots\right)
\end{equation*}
be a sequence defined by $\U_0 := T_2(0)$ and $\U_{n+1} := T(\U_n)$ for $n\ge0$. Additionally 
define the sequence
\begin{equation*}\label{eq:alpha_n_sequence}
\left(\alpha_0, \alpha_1, \ldots, \alpha_n, \ldots\right), 
\end{equation*}
where $\alpha_n := T_1(\U_n)$. 
\end{definition}

We can now prove the convergence of the HMM. We need an 
additional assumption, however, which places a 
 limit on how close to $x_2=0$ one can theoretically take 
$\gamma_1$. 

%%%%%%%%%%%%%%%%%%%% theorem that U_n converges using the lemma
\smallskip
\begin{theorem}\label{thm:convergence}
Keeping the assumptions of Lemma \ref{lem:contraction}, additionally 
assume 
\begin{equation}\label{eq:extra_assumption}
\lim_{\epsilon \downarrow 0}\epsilon \gamma_2/\gamma_1 = 0. 
\end{equation}
Then $\exists \epsilon_0$ such that $\forall \epsilon$ satisfying 
$0 <\epsilon < \epsilon_0$, 
 $\exists \, \U^{\ast}$ to which $(\U_n)_n$ converges as $n \to \infty$. 
Moreover, $\alpha^{\ast}:= T_1(\U^{\ast}) > 0$. 
\end{theorem}

\smallskip
The positivity of $\alpha^{\ast}$ is important to ensure the macroscopic problem 
with slip boundary condition remains coercive and can be compared to 
Theorem \ref{thm:thrm2_p} from Achdou et al. 
\begin{proof}(Theorem \ref{thm:convergence})
By Lemma \ref{lem:contraction}, $T$ is a contraction map on $X:= [\U_0, \infty)$. 
So, the fixed point iteration will converge so long as $\U_n \in X$ for each $n$. 
By the mean value theorem,
\begin{equation}\label{eq:mvt_Un}
\U_{n+1} = T(\U_n) =  T_2(\alpha_n) = \U_0 + \alpha_n \pderiv{T_2}{\alpha}(\xi_n), \qquad \forall n\ge0
\end{equation}
for some $\xi_n \in (0, \alpha_n)$. Since $\partial T_2/\partial \alpha$ (given by 
\eqref{eq:dUdalpha}) is strictly positive and finite for $\alpha>0$,  
the desired result will follow if it can be shown that $\alpha_n > 0$ for all $n$. 

We first show $\alpha_0 >0$, which implies $\U_1 > \U_0$ by \eqref{eq:mvt_Un}. We 
then inductively assume $\U_n > \U_0$ for $n \ge 1$ and show that 
$\alpha_n > 0$ follows. Computing $\alpha_0 = T_1(\U_0) 
= T_1\big(\gamma_2(1+\gamma_1-\gamma_2 - \gamma_1/\gamma_2)\big)$: 
\begin{equation}\label{eq:alpha0}
\alpha_0 =  \frac{C_1 + C_2 \U_0}{C_3 + (1 - \tilde B_3)\U_0} 
\end{equation}
$$
 = \frac{{\left(M \epsilon \gamma_{2} - \tilde B_{3} \gamma_{1} \gamma_{2} + \epsilon \gamma_{1} \gamma_{2} 
+ \tilde B_{3} \gamma_{2}^{2} + \tilde B_{3} \gamma_{1} - \tilde B_{3} \gamma_{2} + \gamma_{1} \gamma_{2} 
- \overline{B}_{2} - \gamma_{1}\right)} {\left(\gamma_{1} - \gamma_{2}\right)}}
{M \epsilon \gamma_{2} - \tilde B_{3} \gamma_{1} \gamma_{2} 
+ 2 \, \epsilon \gamma_{1} \gamma_{2} + \tilde B_{3} \gamma_{2}^{2} - \epsilon \gamma_{2}^{2} + \tilde B_{3} \gamma_{1} 
- \tilde B_{3} \gamma_{2} + \gamma_{1} \gamma_{2} - \gamma_{2}^{2} - \overline{B}_{2} - \gamma_{1} + \gamma_{2}} .
$$
Since $\gamma_1 < \gamma_2$, the numerator will be positive if 
\begin{equation*}
M \epsilon \gamma_{2}/\gamma_1 - \tilde B_{3} \gamma_{2} + \epsilon \gamma_{2} 
+ \tilde B_{3} \gamma_{2}^{2}/\gamma_1 + \tilde B_{3} 
- \tilde B_{3} \gamma_{2}/\gamma_1 + \gamma_{2} - \overline{B}_{2}/\gamma_1  < 1.\label{eq:alpha0_num} 
\end{equation*}
Since $\tilde B = \bigoh{\epsilon}$, the LHS of the inequality vanishes in the limit $\epsilon \searrow 0$ by 
the assumption \eqref{eq:extra_assumption}. Thus the numerator is positive for $\epsilon$ 
sufficiently small. The denominator is positive if
\begin{equation*} \label{eq:alpha0_denom}
1 > -M \epsilon + \tilde B_3 \gamma_1 - 2 \epsilon \gamma_1 - \tilde B_3 \gamma_2 + \epsilon \gamma_2
-\tilde B_3 \gamma_1/\gamma_2 + \tilde B_3 - \gamma_1 + \gamma_2 + \overline{B}_2/\gamma_2 + \gamma_1/\gamma_2.
\end{equation*}
Again by \eqref{eq:extra_assumption}, the RHS vanishes as $\epsilon$ vanishes. Hence
\begin{equation*}
\alpha_0 > 0
\end{equation*}
for $\epsilon$ sufficiently small. 

Next, assume that $\U_n > \U_0$ for $n \ge 1$. 
We just showed that both the numerator $C_1 + C_2 \U_0$ 
and denominator $C_3 + (1-\tilde B_3) \U_0$ in 
\eqref{eq:alpha0} are positive for $\epsilon$ sufficiently 
small. From the estimate \eqref{eq:B3_bound} we
know that both $1-\tilde B_3 > 0$ and 
\begin{equation*}
C_2 > 0 \iff \tilde B_3 \gamma_2 /\gamma_1 - \tilde B_3 + 1 > 0,
\end{equation*}
are true for vanishing $\epsilon$, from which we conclude
\begin{equation*}
C_1 + C_2 \U_n > C_1 + C_2 \U_0 > 0
\end{equation*}
and 
\begin{equation*}
C_3 + (1-\tilde B_3) \U_n> C_3 + (1-\tilde B_3) \U_0 > 0. 
\end{equation*}
This implies $\alpha_n > 0$. 

By the method of induction, we conclude that for all 
$\epsilon$ sufficiently small, $\U_{n} \in X$ and 
$\alpha_n >0$  $\forall n\ge0$. Lemma \ref{lem:contraction} 
then guarantees the existence of a fixed point 
$\U^{\ast}$ such that $\U_n \to \U^{\ast}$ and 
$\U^{\ast} = T(\U^{\ast})$.  
Furthermore, since $T_1(\U)$ is continuous for 
$\U \in X$, we define the slip amount 
\begin{equation*}
\alpha^{\ast} = T_1(\U^{\ast}) = \lim_{n\to\infty} T_1(\U_n) .
\end{equation*}

Since $(\U_n)_n$ is a convergent sequence, it must also be bounded. 
The inductive argument additionally showed the sequence is strictly 
increasing, 
meaning $\U_n < \U^{\ast}$ and hence
$0 < C_3 + (1-\tilde B_3)\U_n < C_3 + (1-\tilde B_3)\U^{\ast}$
for all $n$. As just argued, we also have $C_1 + C_2 \U_n > C_1 +C_2 \U_0 > 0$. 
Combining these inequalities shows that $(\alpha_n)_n$ is uniformly 
bounded below: 
\begin{equation*}
\alpha_n = \frac{C_1 + C_2 \U_n}{C_3 + (1-\tilde B_3) \U_n} 
> \frac{C_1 + C_2 \U_0}{C_3 + (1-\tilde B_3) \U^{\ast}} > 0.
\end{equation*}
Ergo $\alpha^{\ast} >0$ as desired. 
\end{proof} %%%%%%%%%%%%% END proof of main theorem 

\smallskip
Next we show that the slip-amount resulting from the fixed point
iteration vanishes as $\epsilon$ limits to zero, ensuring
the macroscale flow $U$ will limit to a simple Poiseuille flow 
satisfying the no-slip condition along $x_2=0$. This is consistent
with the limiting behavior of the true, rough-wall flow $\ueps$
satisfying no-slip along $\Gamma_{\epsilon}$. 

\smallskip
\begin{corollary} %%%%%%%%%%%% proof alpha --> 0 as eps --> 0  
\label{cor:alphavanish}

Under the same assumptions as Theorem \ref{thm:convergence}, we have
\begin{equation*}
\lim_{\epsilon \downarrow 0} \,\, \alpha^{\ast} = 0.
\end{equation*}
\end{corollary}
\begin{proof}   %% proof of corollary
From the proof of Theorem \ref{thm:convergence} we know $\U_0 < \U^{\ast}$, so that 
\begin{equation*}
\alpha^{\ast} = T_1(\U^{\ast}) = \frac{C_1 + C_2 \U^{\ast}}{C_3 + (1-\tilde B_3) \U^{\ast}}
< \frac{C_1 + C_2 \U^{\ast}}{C_3 + (1-\tilde B_3) \U_0}. 
\end{equation*}
Using $\U_0 = \gamma_2 (1 + \gamma_1 - \gamma_2 - \gamma_1/\gamma_2)$ the RHS of the inequality 
can be written 
\begin{equation}\label{eq:bound_on_alpha}
\frac{C_1/\gamma_2 + C_2/\gamma_2 \, \U^{\ast}}{C_3/\gamma_2 +  (1- \tilde B_3)(1 + \gamma_1 - \gamma_2 -\gamma_1/\gamma_2)}; 
\end{equation}
from the definitions \eqref{eq:T1_Cs} and bounds \eqref{eq:B2_bound} and \eqref{eq:B3_bound}, 
it is clear \eqref{eq:bound_on_alpha} will vanish as 
$\epsilon \to 0$ so long as $\U^{\ast}$ is finite in the limit.

From \eqref{eq:mvt_Un} we have
\begin{equation*}
\U^{\ast} = \lim_{n \to \infty} T(\U_n) = \U_0 + \lim_{n\to\infty}\left( \alpha_n \pderiv{T_2}{\alpha}(\xi_n)\right)
\end{equation*}
for some $\xi_n \in (0,\alpha_n)$. Since $\xi_n$ is positive, we can bound $\partial T_2/\partial \alpha$ 
independently of $n$: 
\begin{equation*}
\pderiv{T_2}{\alpha}(\xi_n) < 1-\gamma_2, 
\end{equation*}
which implies 
\begin{equation*}
\U^{\ast} < \U_0 + (1-\gamma_2) \lim_{n\to\infty} \alpha_n  = \U_0 + (1-\gamma_2) \alpha^{\ast} < \infty
\end{equation*}
as desired. 
\end{proof} 

\smallskip
\begin{corollary} %%%%%%%%%%%% simple corollary where gamma1 = eps^s, gamma2 = eps^t
\label{cor:convergence}
Suppose $\gamma_2= \epsilon^t$ for $0 < t < 1$ and $\gamma_1 = \epsilon^s$ 
for $t < s < 1+t$. Then the assumptions of Theorem \ref{thm:convergence} are met, 
and hence the conclusions follow. 
\end{corollary}

Under an additional technical assumption on the locations $\gamma_1$ and 
$\gamma_2$, we next show the algorithm from Section \ref{subsec:algo} 
will converge rapidly
in the present simplified setting.
So long as the tolerance $\tau = \bigoh{\epsilon^2}$, 
we show 
the algorithm will convergence after just one iteration. 
\begin{theorem}\label{thm:rapid_convergence}  %%%%%%%%%%%%%%%% theorem on rapid convergence
Let $\alpha_0$ and $\alpha_1$ be the first two 
entries of the sequence from Definition \ref{dfn:sequence_Un_and_alpha_n}.
As in Corollary \ref{cor:convergence}, assume $\gamma_2= \epsilon^t$ for $1/2 < t < 1$
and $\gamma_1 = \epsilon^s$ for $t < s < 1+t$, and additionally 
assume that $s > (t+1)/3$. Then
\begin{equation*}
\lim_{\epsilon \downarrow 0} \frac{\alpha_1 - \alpha_0}{\epsilon^2} = 0 .
\end{equation*}
\end{theorem}
\begin{proof}
By definition $\alpha_1 - \alpha_0 = T_1(\U_1) - T_0(\U_0)$; 
from the mean value theorem 
\begin{equation*}
\alpha_1 - \alpha_0 = (\U_1-\U_0)\pderiv{T_1}{\U}(\xi)
\end{equation*}
for some $\xi \in (\U_0, \U_1)$. Let 
\begin{equation*}
b(z) := \frac{\gamma_1^2 - 1 -2 \gamma_1 z}{\gamma_1 - 1- z} 
\end{equation*}
(compare to \eqref{eq:macromap}). Then 
\begin{equation*}
\alpha_1 - \alpha_0 = [b(\alpha_0) - b(0)] (\gamma_2 -1)\pderiv{T_1}{\U}(\xi)
\end{equation*}
Since $(\gamma_2-1) \to 1$ as $\epsilon \to 0$, we analyze only 
$(b(\alpha_0) - b(0)) \partial T_1/\partial\U (\xi)$; 
with the aid of the computer algebra system Sage \cite{sage}, we first 
write its denominator as:  
\begin{align*}
&\big(\epsilon \gamma_{1}^{2} \gamma_{2} + M \epsilon \gamma_{2}^{2} + \tilde B_{3} \gamma_{1} \gamma_{2}^{2} 
- \tilde B_{3} \gamma_{2}^{3} - M \epsilon \gamma_{2} - 2 \tilde B_{3} \gamma_{1} \gamma_{2} - 2 \, \epsilon \gamma_{1} \gamma_{2} 
+ 2 \tilde B_{3} \gamma_{2}^{2} \nonumber \\
&+ \epsilon \gamma_{2}^{2} + \tilde B_{3} \gamma_{1} + \overline{B}_{2} \gamma_{2} 
- \tilde B_{3} \gamma_{2} + \gamma_{1} \gamma_{2} - \gamma_{2}^{2} - \overline{B}_{2} - \gamma_{1} + \gamma_{2}\big)\times \nonumber \\
&{\big(M \epsilon \gamma_{2} + 2 \, \epsilon \gamma_{1} \gamma_{2} - \epsilon \gamma_{2}^{2} + B_{3} \xi + B_{2} - \xi\big)}^{2}. 
\end{align*}
Of the quantities in the first set of parentheses, $\gamma_2$ is asymptotically the largest for vanishing 
$\epsilon$. Since $\xi \in (\U_0, \U_1)$ and $\U_0 = \gamma_2(1+\gamma_1-\gamma_2 - \gamma_1/\gamma_2)$, 
$\xi$ is at least asymptotically similar to $\gamma_2$. Thus we can say that in total the denominator is 
at least asymptotically  
similar to $\gamma_2^3$. 

The numerator of $(b(\alpha_0) - b(0)) \partial T_1/\partial\U (\xi)$ equals
\begin{align} \label{eq:foo_num}
\big(-\tilde B_{3} \epsilon \gamma_{1}^{2} &+ 2 \tilde B_{3} \epsilon \gamma_{1} \gamma_{2} 
- \tilde B_{3} \epsilon \gamma_{2}^{2} + \epsilon \gamma_{1}^{2} + M \epsilon \gamma_{2} + \overline{B}_{2}\big)\times \nonumber  \\
&\big(M \epsilon \gamma_{2} + \tilde B_{3} \gamma_{1} \gamma_{2} + \epsilon \gamma_{1} \gamma_{2} 
- \tilde B_{3} \gamma_{2}^{2} - \tilde B_{3} \gamma_{1} + \tilde B_{3} \gamma_{2} \nonumber \\
&- \gamma_{1} \gamma_{2} + \overline{B}_{2} + \gamma_{1}\big) {\left(\gamma_{1} - \gamma_{2}\right)} {\left(\gamma_{1} - 1\right)} \gamma_{2}.
\end{align}
To prove the desired result, it suffices to show that \eqref{eq:foo_num} divided 
by $\gamma_2^3 \epsilon^2$ limits to zero for vanishing $\epsilon$, meaning 
\begin{equation*}
\lim_{\epsilon \downarrow 0} 
\frac{q_1^{\epsilon} \, q_2^{\epsilon} \, (1-\gamma_2) (\gamma_1/\gamma_2 - 1)}{\gamma_2 \epsilon^2}  = 0 
\end{equation*}
where
\begin{align*}
q_1^{\epsilon} &= -\tilde B_{3} \epsilon \gamma_{1}^{2} + 2 \tilde B_{3} \epsilon \gamma_{1} \gamma_{2} 
- \tilde B_{3} \epsilon \gamma_{2}^{2} + M \epsilon \gamma_{2} + \epsilon \gamma_{1}^{2}  + \overline{B}_{2}\\
q_2^{\epsilon} &= \big(M \epsilon \gamma_{2} + \tilde B_{3} \gamma_{1} \gamma_{2} + \epsilon \gamma_{1} \gamma_{2}
- \tilde B_{3} \gamma_{2}^{2} - \tilde B_{3} \gamma_{1} + \tilde B_{3} \gamma_{2}
- \gamma_{1} \gamma_{2} + \overline{B}_{2} + \gamma_{1}\big)
\end{align*}
Since $\gamma_1$ is the asymptotically largest term in $q_2^{\epsilon}$, the result will 
follow if $ q_1^{\epsilon} \,  \gamma_1/(\gamma_2\epsilon^2) \to 0$ as $\epsilon \searrow 0$. 
Since $M = H \epsilon$ and the bound \ref{eq:B3_bound} shows that $\tilde B_3$ is asymptotically 
similar to $\epsilon$, the first four terms in $q_1^{\epsilon}$ vanish faster than $\epsilon^2$. 
The assumption in the theorem statement that $s > (t+1)/3$ ensures the fifth term in 
$q_1^{\epsilon}$ additionally vanishes. Finally, observe that $\sigma_{\epsilon}$ in 
\eqref{eq:sigma_defn} is asymptotically similar to $\epsilon^{1/2 + 3t/2}$, so that by 
\eqref{eq:B2_bound} $\overline{B}_2$ must vanish faster than $\epsilon^{1+2t}$. Since
$t > 1/2$, the result follows. 
\end{proof} %%%%%%%%%%%%% END proof of rapid convergence 

\section{Numerical results}
\label{sec:laminar_numerics}

We now present two-dimensional numerical tests of the HMM scheme in situations 
both where the periodic homogenization theory 
is applicable and where it is not.  All computations are 
performed using the open source finite element package FEniCS 
\cite{LoggMardalEtAl2012a,LoggWellsEtAl2012a}, and all 
meshes are generated using Gmsh \cite{gmsh}. 
In the first four cases considered, $\epsilon = 0.025$
and $|\Omegamac| = \bigoh{1}$, and $\nu = 1$. 
%\red{[keep?]}Strictly speaking, this set of parameters is not in the asymptotic
%regime analyzed in \cref{subsec:analysis_of_laminar_flow}; however, 
%the assumption \eqref{eq:contained_cond} certainly holds.\red{[keep?]} l
Different values are 
prescribed for the final example of a backwards facing step and
are detailed below. In all cases, the parameter defining the upper boundary 
of the microscopic domain $\gamma = 4\epsilon$. 

All discretizations are performed with the Taylor-Hood elements, 
i.e.\ $P_2$ and $P_1$ basis functions for the velocity and pressure fields, 
respectively \cite{gunzburger,Layton:2008}, and the resulting discrete nonlinear system is 
solved with Newton's method, using the solution to the corresponding Stokes problem as the initial guess. 

The direct numerical simulation (DNS) of the full problem \eqref{eq:NS_eps} 
is computed with a large number of elements and is compared with (i) the 
 1st order approximation satisfying the no-slip condition along $\Gamma$ (so-called
because of \eqref{boundary_err})
 and (ii) the HMM approximation satisfying the coupled system 
\eqref{eq:hmm_system}; the 1st order approximation and macroscale
HMM function are both computed on the 
same mesh. The coupled HMM system is solved iteratively, using 
the algorithm from Section \ref{subsec:algo}. The tolerance $\tau$ for the error
between successive slip amounts $\alpha$ is set to be 
$\tau = \epsilon^2$. % = (0.025)^2 = 0.000625$. 

%%%%%%%%%%%%%%%%%%%%%%%%%
%%%
%%% Big figure w all the cartoon rough domains
%%%
\begin{figure}[tbhp]
\centering
\subfloat[Periodic, sinusoidal roughness]
{
\label{fig:domain_a}
%\begin{tikzpicture}[xscale=6.5,yscale=3.8]
\begin{tikzpicture}[xscale=4.875,yscale=2.85]
\draw[loosely dotted]  (0,1) -- (0,0)node[left]{$(0, 0)$};
\draw[dotted] (0,0) -- (1,0) node[right]{$\Gamma$} ;
\draw[loosely dotted] (1,0) -- (1,1)node[above]{$(1,1)$};
\draw (0,1) -- (1,1);
\draw[black, domain = 0:0.1] plot(\x, {0.025*(cos(40*pi*\x r) -1)});
\draw[black, domain = 0.1:0.2] plot(\x, {0.025*(cos(40*pi*\x r) -1)});
\draw[black, domain = 0.2:0.3] plot(\x, {0.025*(cos(40*pi*\x r) -1)});
\draw[black, domain = 0.3:0.4] plot(\x, {0.025*(cos(40*pi*\x r) -1)});
\draw[black, domain = 0.4:0.5] plot(\x, {0.025*(cos(40*pi*\x r) -1)});
\draw[black, domain = 0.5:0.6] plot(\x, {0.025*(cos(40*pi*\x r) -1)});
\draw[black, domain = 0.6:0.7] plot(\x, {0.025*(cos(40*pi*\x r) -1)});
\draw[black, domain = 0.7:0.8] plot(\x, {0.025*(cos(40*pi*\x r) -1)});
\draw[black, domain = 0.8:0.9] plot(\x, {0.025*(cos(40*pi*\x r) -1)});
\draw[black, domain = 0.9:1.0] plot(\x, {0.025*(cos(40*pi*\x r) -1)});
\end{tikzpicture}
}
\subfloat[Periodic, ``sawtooth" roughness]
{
\label{fig:domain_b}
%\begin{tikzpicture}[xscale=6.5,yscale=3.8]
\begin{tikzpicture}[xscale=4.875,yscale=2.85]
\draw[loosely dotted]  (0,0.5) -- (0,0);
\draw[dotted] (0,0)node[left]{$(0,0)$} -- (1,0) node[right]{$\Gamma$} ;
\draw[loosely dotted] (1,0) -- (1,0.5)node[above]{$(1,0.5)$};
\draw[black, domain=0:1.0] plot (\x, {0.5+0.125*sin(2.0*pi*\x r)});
%\node at (1.3,0) {$\Gamma_0$};
%\draw [<-] (1.05,0) -- (1.2,0);
\foreach  \k in {1,...,20}
	\draw (\k*0.05-0.05,0) -- (\k*1/20, -1/20) -- (\k*0.05,0);
%\node at (0.5, -0.2) {$\Gamma^{\epsilon}$};
\end{tikzpicture}
}
\newline
\subfloat[Periodic roughness modulated by a smooth function]
{
\label{fig:domain_c}
%\begin{tikzpicture}[xscale=6.5,yscale=3.8]
\begin{tikzpicture}[xscale=4.875,yscale=2.85]
\draw[loosely dotted]  (0,1) -- (0,0)node[left]{$(0,0)$};
\draw[dotted] (0,0) -- (1,0) node[right]{$\Gamma$} ;
\draw[loosely dotted] (1,0) -- (1,1)node[above]{$(1,1)$};
\draw (0,1) -- (1,1);
\draw[black, domain = 0:0.1] plot(\x, {0.025*(cos(40*pi*\x r) -1)*(0.5+(sin(sqrt(2)*2*pi*\x r))*(sin(sqrt(2)*2*pi*\x r)))});
\draw[black, domain = 0.1:0.2] plot(\x, {0.025*(cos(40*pi*\x r) -1)*(0.5+(sin(sqrt(2)*2*pi*\x r))*(sin(sqrt(2)*2*pi*\x r)))});
\draw[black, domain = 0.2:0.3] plot(\x, {0.025*(cos(40*pi*\x r) -1)*(0.5+(sin(sqrt(2)*2*pi*\x r))*(sin(sqrt(2)*2*pi*\x r)))});
\draw[black, domain = 0.3:0.4] plot(\x, {0.025*(cos(40*pi*\x r) -1)*(0.5+(sin(sqrt(2)*2*pi*\x r))*(sin(sqrt(2)*2*pi*\x r)))});
\draw[black, domain = 0.4:0.5] plot(\x, {0.025*(cos(40*pi*\x r) -1)*(0.5+(sin(sqrt(2)*2*pi*\x r))*(sin(sqrt(2)*2*pi*\x r)))});
\draw[black, domain = 0.5:0.6] plot(\x, {0.025*(cos(40*pi*\x r) -1)*(0.5+(sin(sqrt(2)*2*pi*\x r))*(sin(sqrt(2)*2*pi*\x r)))});
\draw[black, domain = 0.6:0.7] plot(\x, {0.025*(cos(40*pi*\x r) -1)*(0.5+(sin(sqrt(2)*2*pi*\x r))*(sin(sqrt(2)*2*pi*\x r)))});
\draw[black, domain = 0.7:0.8] plot(\x, {0.025*(cos(40*pi*\x r) -1)*(0.5+(sin(sqrt(2)*2*pi*\x r))*(sin(sqrt(2)*2*pi*\x r)))});
\draw[black, domain = 0.8:0.9] plot(\x, {0.025*(cos(40*pi*\x r) -1)*(0.5+(sin(sqrt(2)*2*pi*\x r))*(sin(sqrt(2)*2*pi*\x r)))});
\draw[black, domain = 0.9:1.0] plot(\x, {0.025*(cos(40*pi*\x r) -1)*(0.5+(sin(sqrt(2)*2*pi*\x r))*(sin(sqrt(2)*2*pi*\x r)))});
\end{tikzpicture}
}
\subfloat[Quasi-periodic roughness]
{
\label{fig:domain_d}
%\begin{tikzpicture}[xscale=6.5,yscale=3.8]
\begin{tikzpicture}[xscale=4.875,yscale=2.85]
\draw[loosely dotted]  (0,1) -- (0,0)node[left]{$(0,0)$};
\draw[dotted] (0,0) -- (1,0) node[right]{$\Gamma$} ;
\draw[loosely dotted] (1,0) -- (1,1)node[above]{$(1,1)$};
\draw (0,1) -- (1,1);
\draw[black, domain = 0:0.1] plot(\x, {0.025*(sin(sqrt(2)*40*pi*\x r) +sin(2/0.025*pi*\x r)-2.25)});
\draw[black, domain = 0.1:0.2] plot(\x, {0.025*(sin(sqrt(2)*40*pi*\x r) +sin(2/0.025*pi*\x r)-2.25)});
\draw[black, domain = 0.2:0.3] plot(\x, {0.025*(sin(sqrt(2)*40*pi*\x r) +sin(2/0.025*pi*\x r)-2.25)});
\draw[black, domain = 0.3:0.4] plot(\x, {0.025*(sin(sqrt(2)*40*pi*\x r) +sin(2/0.025*pi*\x r)-2.25)});
\draw[black, domain = 0.4:0.5] plot(\x, {0.025*(sin(sqrt(2)*40*pi*\x r) +sin(2/0.025*pi*\x r)-2.25)});
\draw[black, domain = 0.5:0.6] plot(\x, {0.025*(sin(sqrt(2)*40*pi*\x r) +sin(2/0.025*pi*\x r)-2.25)});
\draw[black, domain = 0.6:0.7] plot(\x, {0.025*(sin(sqrt(2)*40*pi*\x r) +sin(2/0.025*pi*\x r)-2.25)});
\draw[black, domain = 0.7:0.8] plot(\x, {0.025*(sin(sqrt(2)*40*pi*\x r) +sin(2/0.025*pi*\x r)-2.25)});
\draw[black, domain = 0.8:0.9] plot(\x, {0.025*(sin(sqrt(2)*40*pi*\x r) +sin(2/0.025*pi*\x r)-2.25)} );
\draw[black, domain = 0.9:1.0] plot(\x, {0.025*(sin(sqrt(2)*40*pi*\x r) +sin(2/0.025*pi*\x r)-2.25)});
\end{tikzpicture}
}
\newline
\subfloat[Backwards facing step with roughness in the recirculation region. The two
dots mark the locations of the micro-domains discussed in Section \ref{subsec:bfs}.]
{
\label{fig:domain_e}
\begin{tikzpicture}[xscale=0.5,yscale=0.5]
\draw[dotted]  (0,2)node[above]{$(0,2)$} -- (0,1);
\draw (0,1)node[below]{$(0,1)$} -- (5, 1); 
\draw (5,1) -- (5,0);
\draw (5,0) node[below]{$(5,0)$} -- (6,0);
%\draw (0,0) node[below]{$x_1 = 0$}; 
%\draw[dotted] (6,0) -- (16,0); 
\draw (16,0) -- (23,0) node[below]{$(23,0)$}; 
\draw (23,2) -- (0,2); 
\draw[dotted] (23,0) -- (23,2);
%\draw[green, dotted] (5,0.25) -- (23,0.25); % draw where we'll show the plot of u1 vs x1
\draw[black, domain = 6:6.5] plot( \x, {0.1*(cos(2*2*pi*(\x) r) -1.0)} );
\draw[black, domain = 6.5:7.0] plot( \x, {0.1*(cos(2*2*pi*(\x) r) -1.0)} );
\draw[black, domain = 7.0:7.5] plot( \x, {0.1*(cos(2*2*pi*(\x) r) -1.0)} );
\draw[black, domain = 7.5:8.0] plot( \x, {0.1*(cos(2*2*pi*(\x) r) -1.0)} );
\draw[black, domain = 8:8.5] plot( \x, {0.1*(cos(2*2*pi*(\x) r) -1.0)} );
\draw[black, domain = 8.5:9.0] plot( \x, {0.1*(cos(2*2*pi*(\x) r) -1.0)} );
\draw[black, domain = 9.0:9.5] plot( \x, {0.1*(cos(2*2*pi*(\x) r) -1.0)} );
\draw[black, domain = 9.5:10.0] plot( \x, {0.1*(cos(2*2*pi*(\x) r) -1.0)} );
\draw[black, domain = 10:10.5] plot( \x, {0.1*(cos(2*2*pi*(\x) r) -1.0)} );
\draw[black, domain = 10.5:11.0] plot( \x, {0.1*(cos(2*2*pi*(\x) r) -1.0)} );
\draw[black, domain = 11.0:11.5] plot( \x, {0.1*(cos(2*2*pi*(\x) r) -1.0)} );
\draw[black, domain = 11.5:12.0] plot( \x, {0.1*(cos(2*2*pi*(\x) r) -1.0)} );
\draw[black, domain = 12:12.5] plot( \x, {0.1*(cos(2*2*pi*(\x) r) -1.0)} );
\draw[black, domain = 12.5:13.0] plot( \x, {0.1*(cos(2*2*pi*(\x) r) -1.0)} );
\draw[black, domain = 13.0:13.5] plot( \x, {0.1*(cos(2*2*pi*(\x) r) -1.0)} );
\draw[black, domain = 13.5:14.0] plot( \x, {0.1*(cos(2*2*pi*(\x) r) -1.0)} );
\draw[black, domain = 14:14.5] plot( \x, {0.1*(cos(2*2*pi*(\x) r) -1.0)} );
\draw[black, domain = 14.5:15.0] plot( \x, {0.1*(cos(2*2*pi*(\x) r) -1.0)} );
\draw[black, domain = 15.0:15.5] plot( \x, {0.1*(cos(2*2*pi*(\x) r) -1.0)} );
\draw[black, domain = 15.5:16.0] plot( \x, {0.1*(cos(2*2*pi*(\x) r) -1.0)} );
\draw (16,0)node[below]{$(16,0)$}; 
\draw [black, fill=black] (7.5,0.0) circle [radius=0.10];
\draw [black, fill=black] (13.5,0.0) circle [radius=0.10];
\end{tikzpicture}
}
\caption{Sketches of the rough domains used to test the numerical method.}
\label{fig:domain_cartoons}
\end{figure}
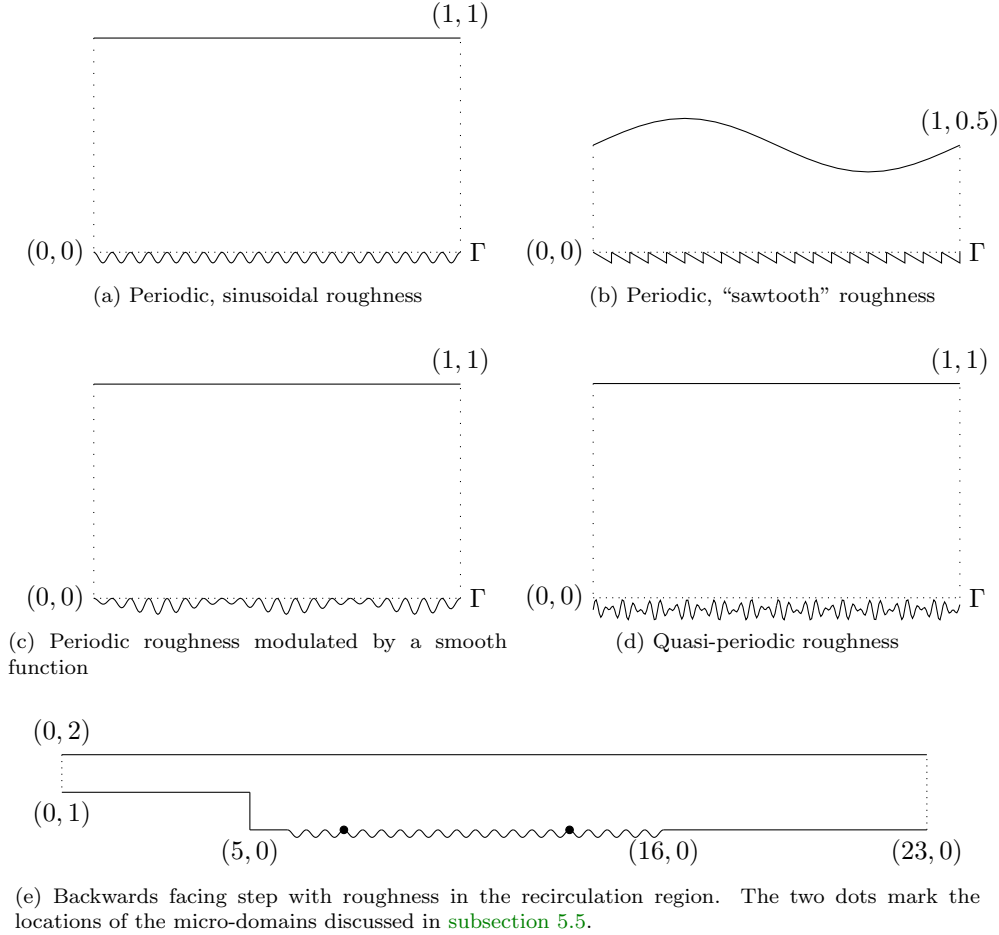

%%%%%%%%%%%%%%
%% Channel flow with periodic sinusoidal periodic roughness
%%%%%%%%%%%%%%
\subsection{Flow in a channel with periodic roughness.}\label{subsection_channel_comp_1}
First consider a channel domain with periodic roughness, as in Figure \ref{fig:domain_a}.
The macroscopic domain is simply $\Omegamac = [0,1]^2$, and the roughness is 
parameterized by the function 
$\varphi^{\epsilon}(x_1) = \epsilon/2( \cos(2\pi x_1/\epsilon)-1)$.
The no-slip condition is applied at the upper boundary $x_2 = 1$, and periodic
boundary conditions are applied on the left/right boundaries. A constant pressure 
gradient $-\nabla p = (1 \,\,\,\, 0)^T$ drives the flow from left to right. 
%\red{[keep?]}$\Omega^{\epsilon}$ is discretized with 50554 cells, while
%the discretization of $\Omegamac$ and $\Omegamic$ contain only 3200 and 927 cells, 
%totaling $6.3\%+1.8\%=8.1\%$ the amount of rough domain cells.\red{[keep?]}

In the setting just described, the macroscopic solution $U$ is one dimensional. Only the horizontal
component of the flow is nonzero, and it depends only on the wall-normal variable $x_2$. 
In this case, only one micro-domain is needed, and periodic boundary 
conditions can be prescribed along the left/right computational 
boundaries of the micro-domain ($x_1 = 0$ and $x_1 = \epsilon$) for simplicity, 
as discussed in Remark \ref{rmk:bc_micro}. The free stream condition 
\eqref{eq:free_stream} is then applied along the upper computational boundary $x_2 =\gamma = 4\epsilon$. 

In this setting of a domain with periodic side walls, there is a small difficulty using the FEniCS 
software package whose workaround we briefly describe. As mentioned, the mixed $P_2$--$P_1$ finite elements
are used to numerically solve the stationary Navier-Stokes system in $\Omegamic$.  

The discrete equations posed with Dirichlet condtions for the velocity vector $u$ 
along $\partial \Omegamic_{\rm noslip}$ and $x_2 = \gamma$ and periodic conditions
at the side-walls have a unique solution from the standard finite element theory. 
The pressure field exists to (weakly) enforce incompressibility and does not 
formally satisfy any boundary conditions. 
For mixed finite element methods however FEniCS requires that 
either \emph{both} function spaces or \emph{neither} have periodic boundaries. 
Since we want the velocity $u$ to be periodic for $x_2 \in [0,\gamma]$, we
consequently assume the microscale pressure can be written as $p = p_{\rm per} + \tilde p$, where $p_{\rm per}$ is periodic and
$-\nabla \tilde p = (1 \,\,\,\, 0)^T$ is a body force that is consistent with the macroscopic forcing 
driving the flow throughout the channel. We then compute the solution pair $(u,p_{\rm per})$. 
Whenever periodic boundary conditions 
are imposed at $x_1=s$ and $x_1 = s+L^{\rm mic}$ in the descriptions below, this decomposition is used. 

\begin{figure}[h!]
\centering
\includegraphics[width = 5in, height = 3.33in]{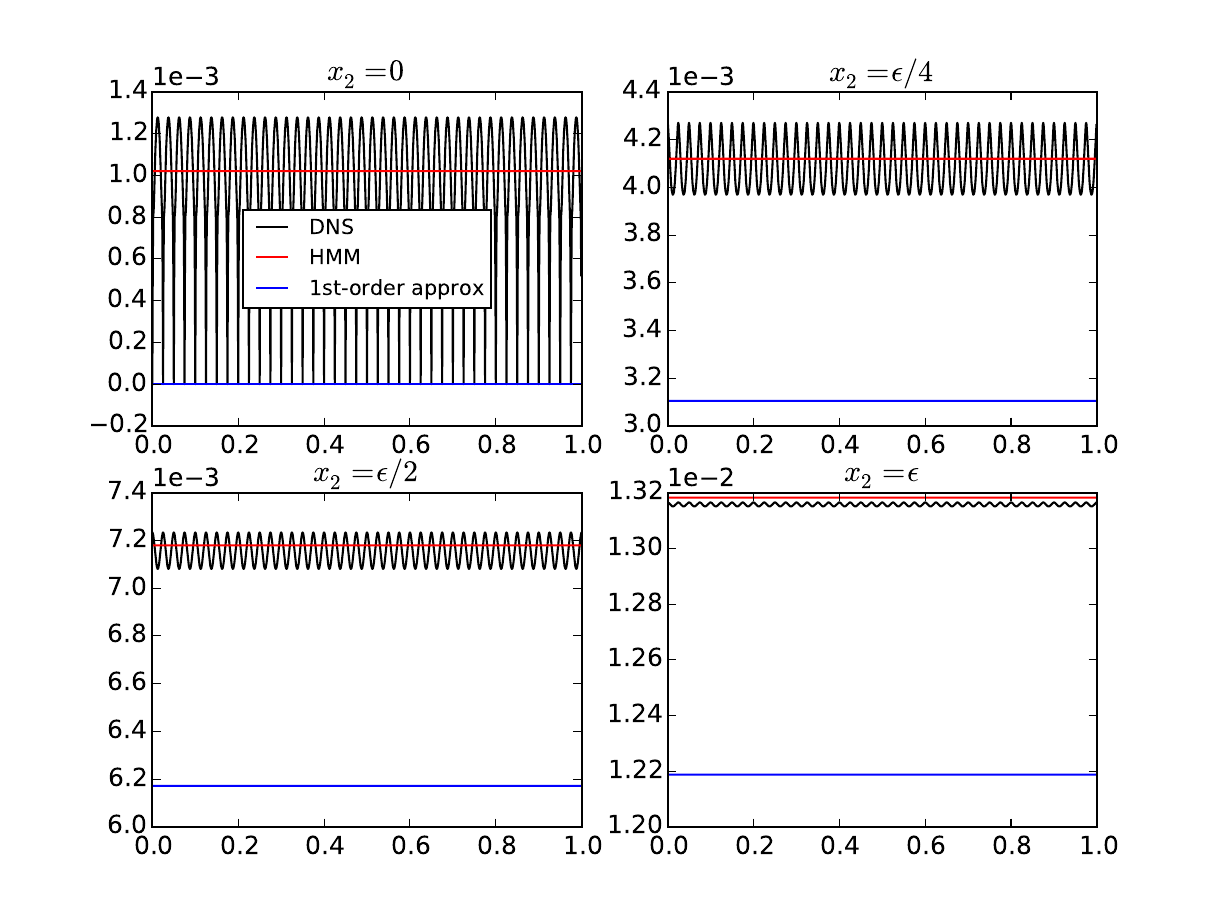}
\caption{Horizontal component of the flow $u_1$ versus $x_1$ plotted at various heights $x_2$ 
for the domain shown in Figure \ref{fig:domain_a}.}
\label{fig:u1channel}
\end{figure}

\begin{figure}[h!]
\centering
\includegraphics[width = 5in, height = 3.33in]{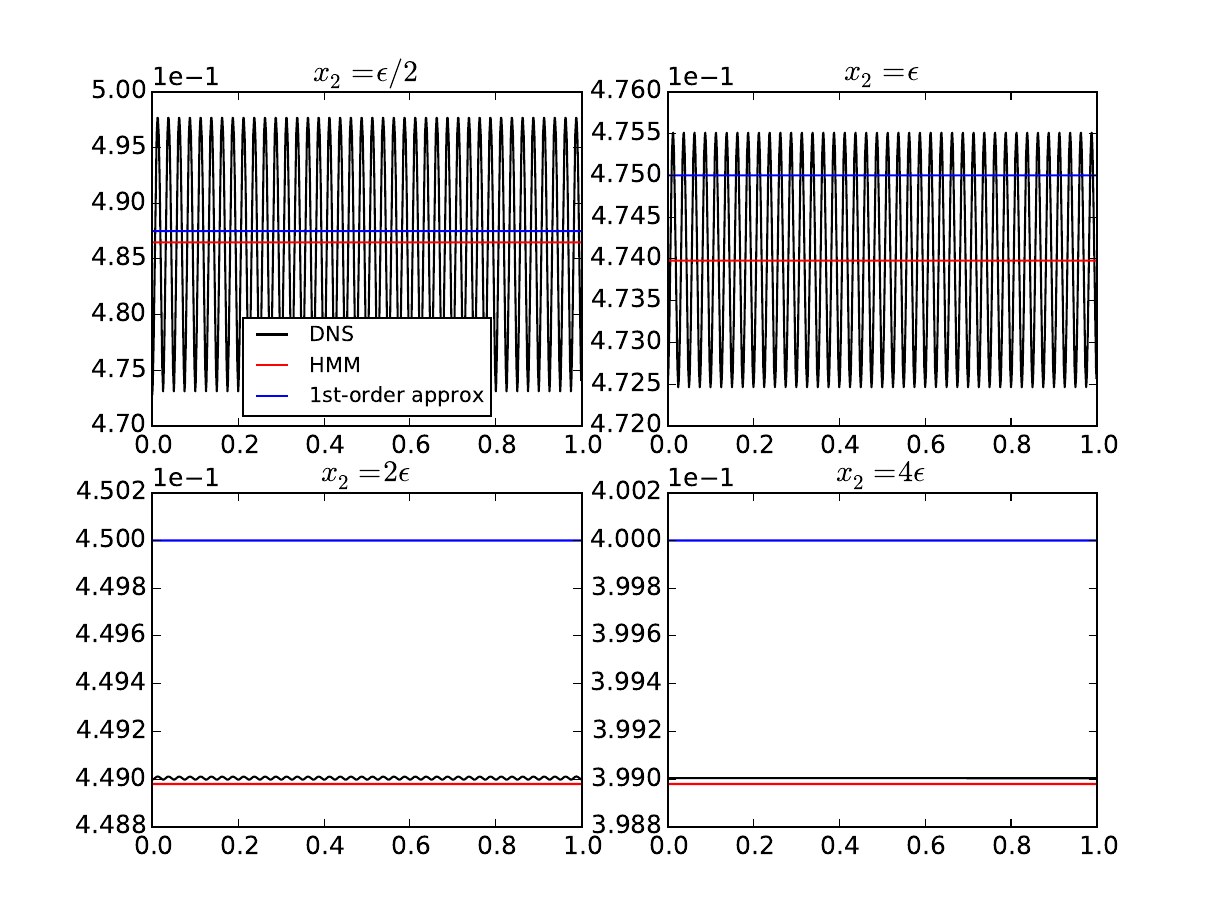}
\caption{Shear $\partial u_1/\partial x_2$ versus $x_1$ plotted at various heights $x_2$
for the domain shown in Figure \ref{fig:domain_a}.}
\label{fig:shearchannel}
\end{figure}

Figures \ref{fig:u1channel} and \ref{fig:shearchannel} plot $u_1$ and $\partial u_1/\partial x_2$, 
respectively, 
as functions of $x_1$ for several values of $x_2$ near the wall. 
Also computed was an HMM solution using the more general strategy for the projection 
operator $\pi$ defined by the constraints detailed in Appendix 
\ref{appendix:micro_bc} (not pictured). 
The resulting slip amount differed from 
the one computed with periodic boundary conditions in the micro domain only by 0.9\%. 
%It is interesting to
%note, however, that the micro-domain flow $u_1$ along the 
%line $x_2=0$ over which the average is taken has a lower value in the Dirichlet case than in the periodic 
%case. Since the slip amount is nearly identical, the shear $\partial u_1/\partial x_2$ accordingly is lower as well.


%%%%%%%%%%%%%%
%% Flow in a wavy tube  with periodic sawtooth roughness
%%%%%%%%%%%%%%
\subsection{Nonsquare domain with periodic roughness.}
Next, consider a nonsquare macroscopic domain with periodic, ``sawtooth" roughness
as shown in Figure \ref{fig:domain_b}. Let $h(x_1) := 0.5 - 0.125\sin(2\pi x_1)$. Then 
\begin{equation*}
\Omegamac = \left\{(x_1,x_2) \vert 0\le x_1\le 1, \,\, 0\le x_2 \le h(x_1)\right\}.  
\end{equation*}
The roughness is parameterized by the periodic function 
$\varphi^{\epsilon}(x_1) = -3\epsilon/4(x_1/\epsilon -\floor*{x_1/\epsilon}) $,
where $\floor*{\cdot}$ denotes the floor function. 
The no-slip condition is applied on the domain's upper, curved boundary, and periodic boundary
conditions are applied on the left/right boundaries. A constant body force $f = (1 \,\,\,\, 0)^T$ 
drives the flow from left to right. 
%\red{[keep?]}$\Omega^{\epsilon}$ is discretized with 41898 cells, while
%$\Omegamac$ is discretized with only 1854 cells and 
%each of the five discretized $\Omegamic$ domains below contain 903 cells, totaling
% $4.4\%+ 11.0\% = 15.4\%$ the amount of rough domain cells.\red{[keep?]}

To compute the HMM approximations, we use the algorithm from Section \ref{subsec:algo}
and set 
\begin{equation*}
\{s_1,s_2,s_3,s_4,s_5\} = \{0, 0.25, 0.5, 0.75, 1\}
\end{equation*}
chosen to capture influence of the macroscopic curvature of $\Omega^{\rm mac}$.
However,
the percent difference between the largest and smallest resulting values of
slip amounts is a negligible 0.3\%, indicating that simply performing one micro-solve
at a single $s_j$ is sufficient in this case. 

\begin{figure}[h!]
\centering
\includegraphics[width = 5in, height = 3.33in]{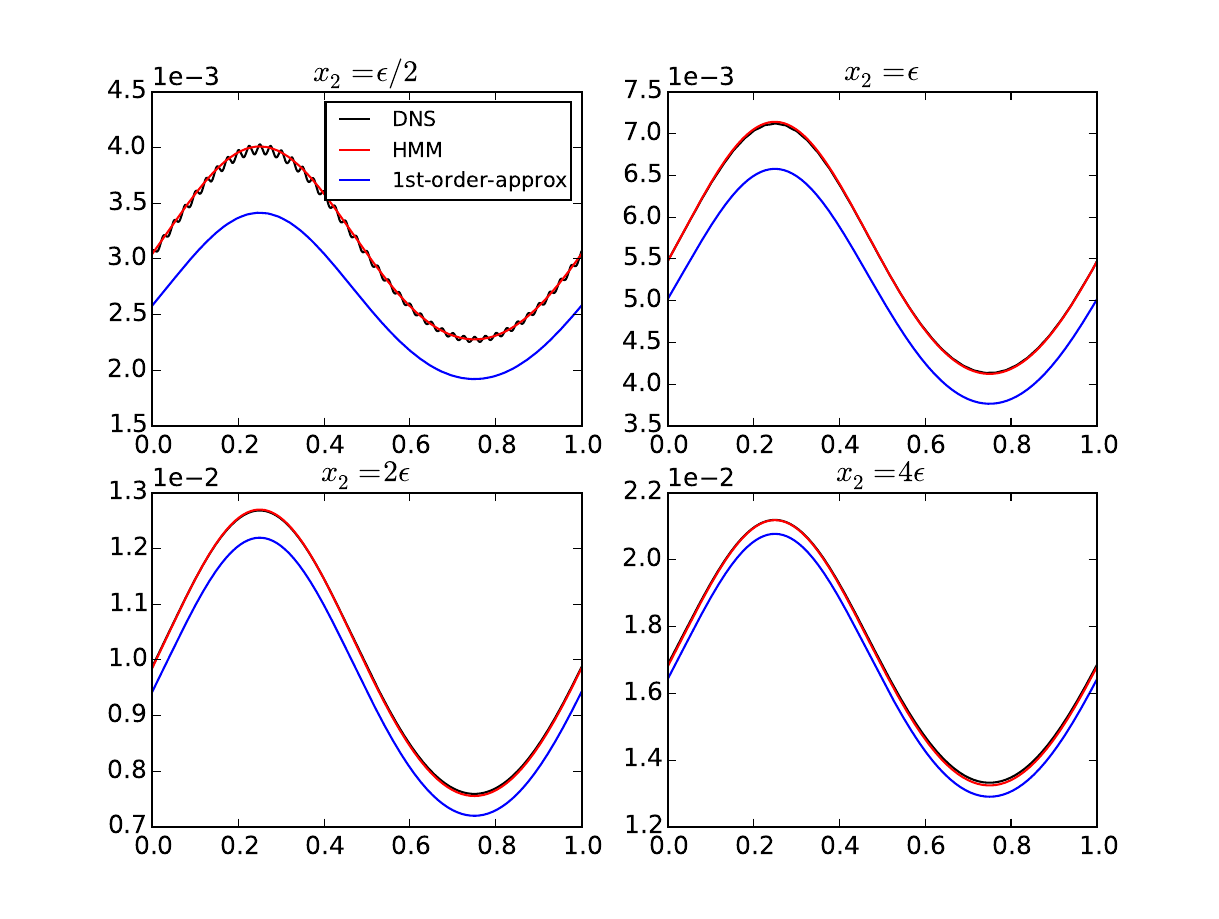}
\caption{Horizontal component of the flow $u_1$ versus $x_1$ plotted at various heights $x_2$ 
for the domain shown in Figure \ref{fig:domain_b}.}
\label{fig:u1wavy}
\end{figure}

\begin{figure}[h!]
\centering
\includegraphics[width = 5in, height = 3.33in]{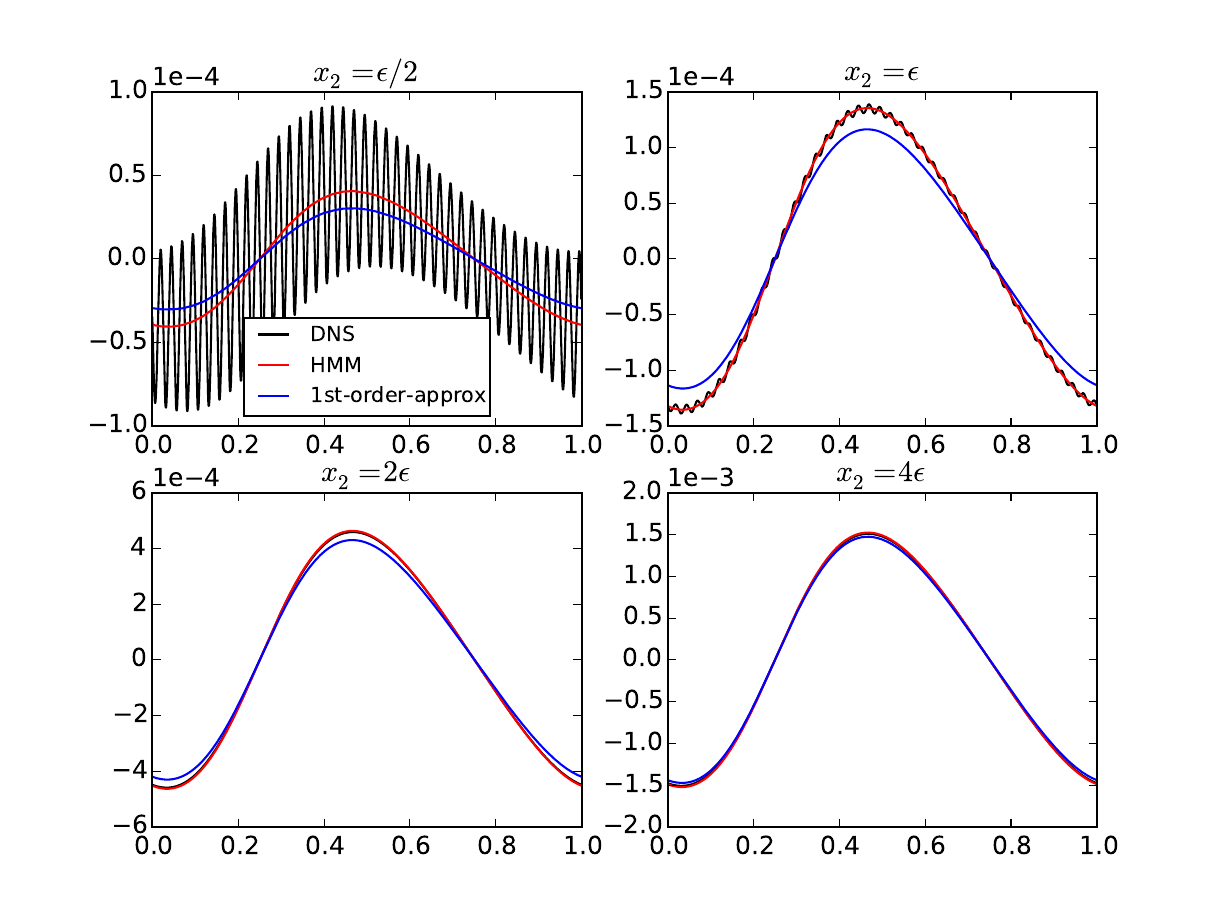}
\caption{Vertical component of the flow $u_2$ versus $x_1$ plotted at various heights $x_2$ 
for the domain shown in Figure \ref{fig:domain_b}; compared to $u_1$ in Figure \ref{fig:u1wavy}, 
$u_2$ is roughly one order of magnitude smaller.}
\label{fig:u2wavy}
\end{figure}

In contrast to the previous example, the macroscopic flow $U$ is not one dimensional, i.\,e.\ both 
$U_2$ and $\partial U_1/\partial x_1$ are nonzero, as can be seen from the DNS curve 
in Figures \ref{fig:u1wavy} and \ref{fig:u2wavy}. However, since at $x_2 = 4\epsilon$ the horizontal component of 
the flow is approximately one order of magnitude larger than the vertical component for a
given $x_1$, it is reasonable to attempt to approximate the vertical component as being zero
and compute with periodic boundary conditions and the free-stream condition \eqref{eq:free_stream}. 
Similar to the previous numerical example in Section \ref{subsection_channel_comp_1},
a difference of about one percent is observed between the slip amount computed this way and
the slip amount using the more general projection $\pi_j(U)$ from Appendix \ref{appendix:micro_bc}. 

\begin{figure}[h!]
\centering
\includegraphics[width = 5in, height = 3.33in]{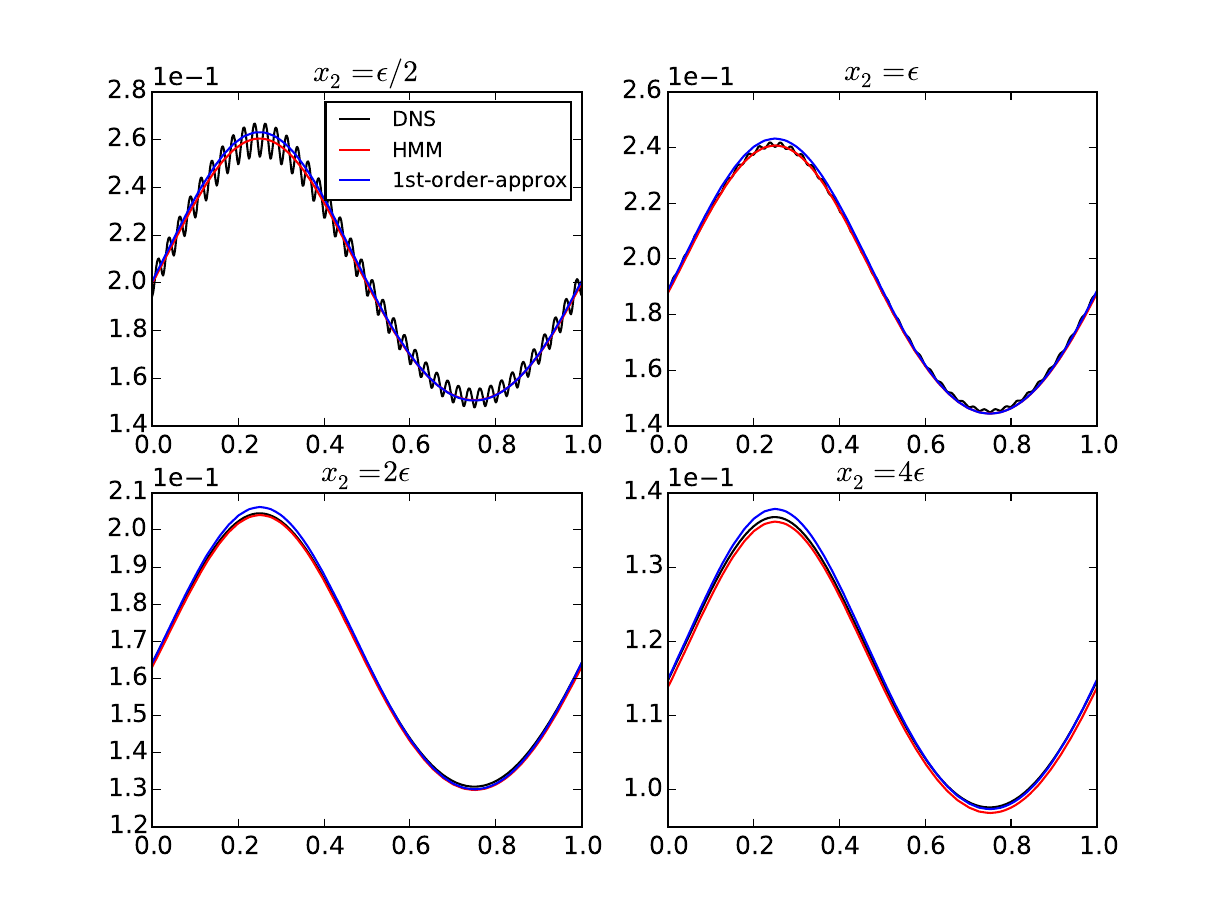}
\caption{Shear $\partial u_1/\partial x_2$ versus $x_1$ plotted at various heights $x_2$
for the domain shown in Figure \ref{fig:domain_b}.}
\label{fig:shearwavy}
\end{figure}
The results in Figures \ref{fig:u1wavy} and \ref{fig:shearwavy} again show $u_1$ and $\partial u_1/\partial x_2$ 
versus $x_1$ for various values of $x_2$. 

%%%%%%%%%%%%%%
%% Channel flow with periodic sinusoidal periodic roughness modulated by smooth function
%%%%%%%%%%%%%%
\subsection{Flow in a channel with non-periodic roughness.}\label{subsection_channel_modulated}
Consider again $\Omegamac = [0,1]^2$, but now let the roughness be parameterized by
\begin{align*}
\zeta^{\epsilon}(x_1) &= \beta(x_1) \varphi^{\epsilon}(x_1) \nonumber \\
\beta(x_1) &= \sin^2\left( \sqrt{2} \, 2\pi x_1\right) + 0.5 \nonumber \\
\varphi^{\epsilon}(x) &= \epsilon/2( \cos(2\pi x_1/\epsilon)-1), 
\end{align*}
so that the periodic roughness is modulated by a smooth function as shown in Figure \ref{fig:domain_c}. 
The no-slip condition is applied at $x_2 = 1$, periodic boundary conditions are enforced 
at $x_1 = 0$ and $x_1 = 1$, and a uniform pressure gradient 
$-\nabla p = (1 \,\,\,\, 0)^T$ again drives the flow from left to right. 
%\red{[keep?]}$\Omega^{\epsilon}$ is discretized with 51210 cells, while $\Omegamac$  
%is discretized with only 3200 cells and the sum of the discrete cells in the 
%7 separate micro domains below totals 9779, amounting to $6.3\% + 12.8\% = 19.1\%$
%of the rough domain cells. \red{[keep?]}

The algorithm from Section \ref{subsec:algo}
is used with 
\begin{equation*}
\{s_1,s_2,s_3,s_4,s_5,s_6,s_7\} = \{0, 0.15, 0.35, 0.525, 0.675, 0.875, 0.975\}
\end{equation*}
chosen to capture the large scale curvature of $\beta$. The asymptotic analysis
presented in Section \ref{subsection:asymptotic_analysis} suggests it is sufficient
to simply compute in a single microscopic domain with roughness parameterized only by $\varphi^{\epsilon}$
and then multiply the resulting slip amount by $\beta(x_1)$ in the effective boundary condition
\eqref{eq:mac_wall_law}. However, we chose to apply the general HMM algorithm to mimic the situation for
which an analytic formula for $\beta$ is not known. The total slip amount used in the macroscale domain
is a linear interpolation of the $\alpha_j$, $1 \le j \le 7$.

\begin{figure}[h!]
\centering
\includegraphics[width = 5in, height = 3.33in]{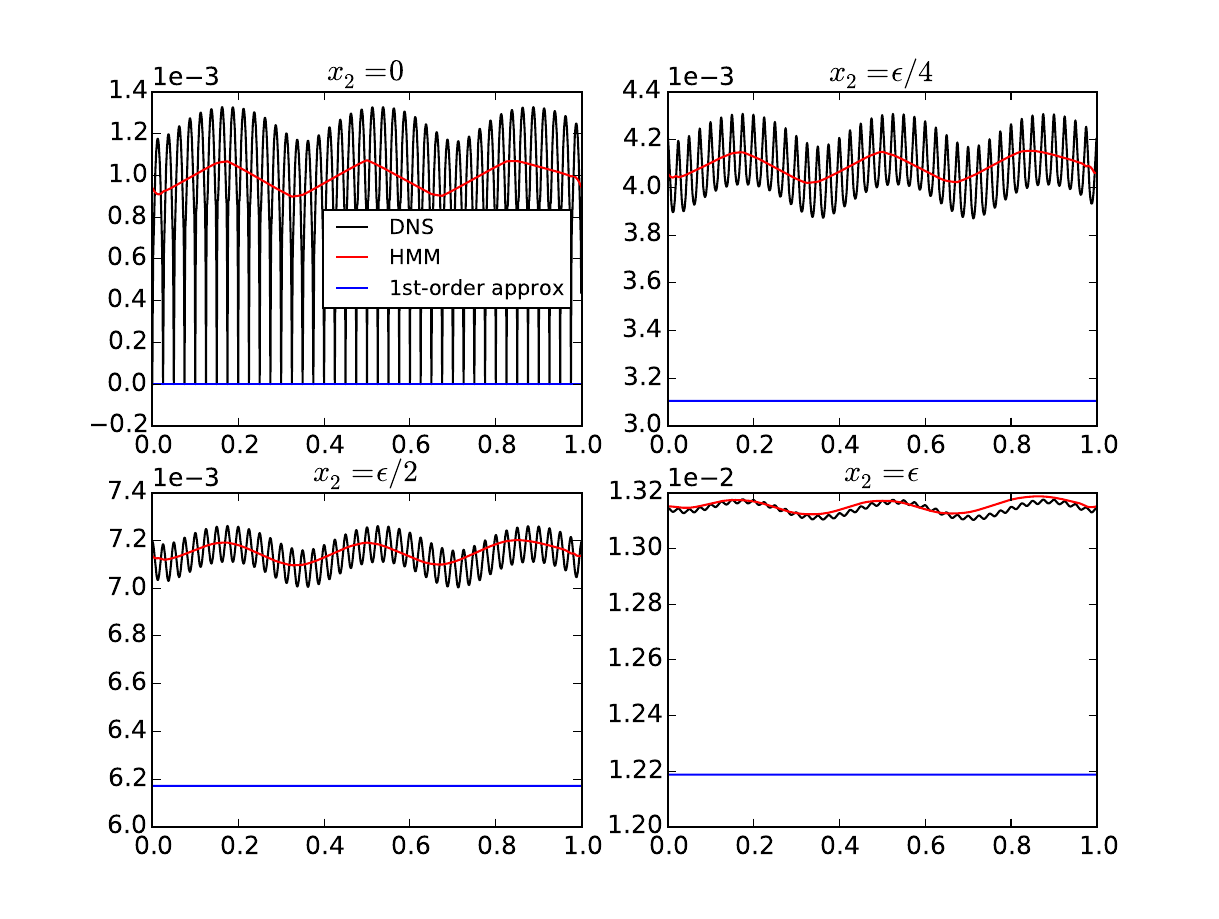}
\caption{Horizontal component of the flow $u_1$ versus $x_1$ plotted at various heights $x_2$
for the domain shown in Figure \ref{fig:domain_c}.}
\label{fig:u1vsx1modulated}
\end{figure}

\begin{figure}[h!]
\centering
\includegraphics[width = 5in, height = 3.33in]{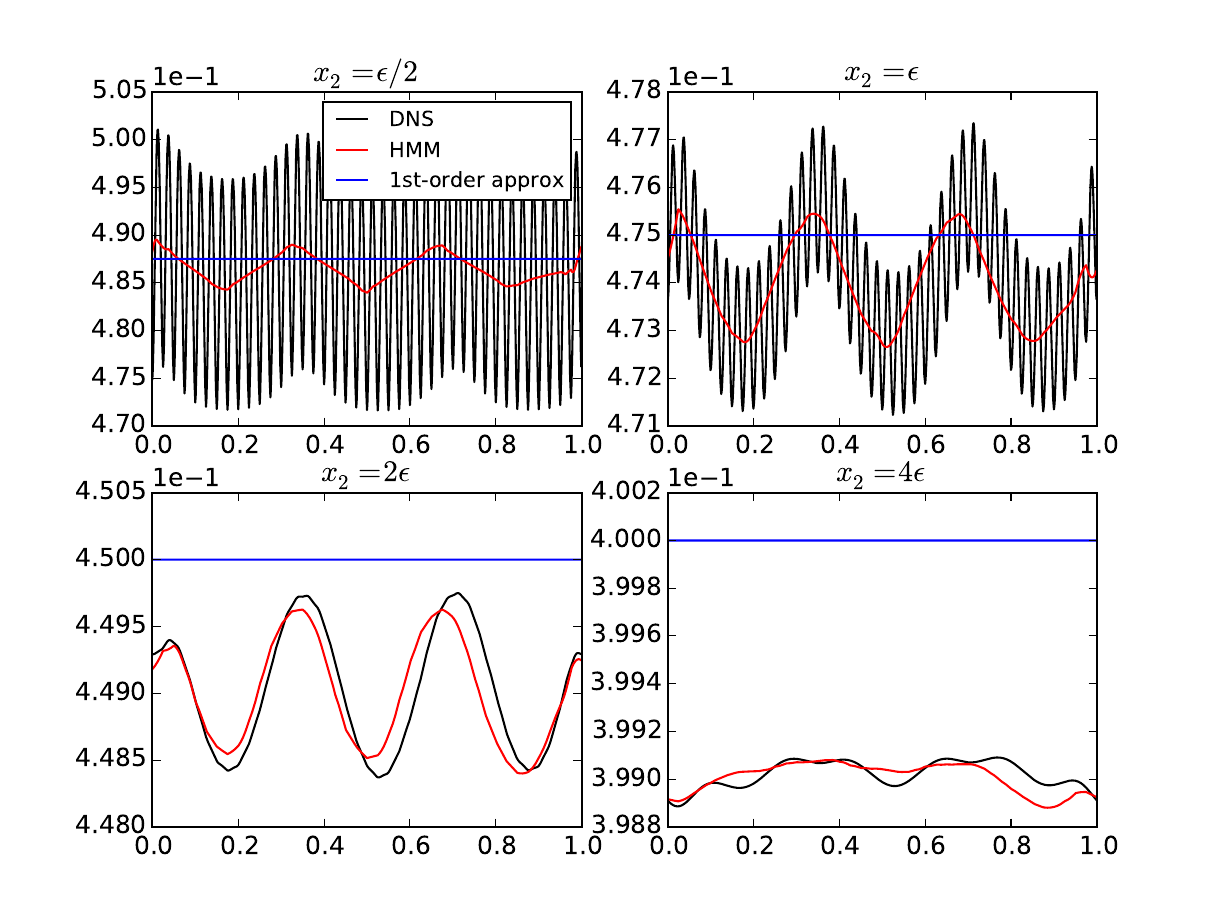}
\caption{Shear $\partial u_1/\partial x_2$ versus $x_1$ plotted at various heights $x_2$
for the domain shown in Figure \ref{fig:domain_c}.}
\label{fig:shearvsx1modulated}
\end{figure}

In this case, the percent difference between the largest and smallest slip amounts 
is a non-negligible 23.4\%. Figures \ref{fig:u1vsx1modulated} and \ref{fig:shearvsx1modulated} illustrate
HMM's successful capturing of the slip amount's horizontal dependence.


%%%%%%%%%%%%%%
%% Channel flow with quasi periodic sinusoidal periodic roughness 
%%%%%%%%%%%%%%
\subsection{Flow in a channel with quasi-periodic roughness.}\label{subsection:channel_quasiper}
\begin{figure}[h!]
\centering
\includegraphics[width = 5in, height = 3.33in]{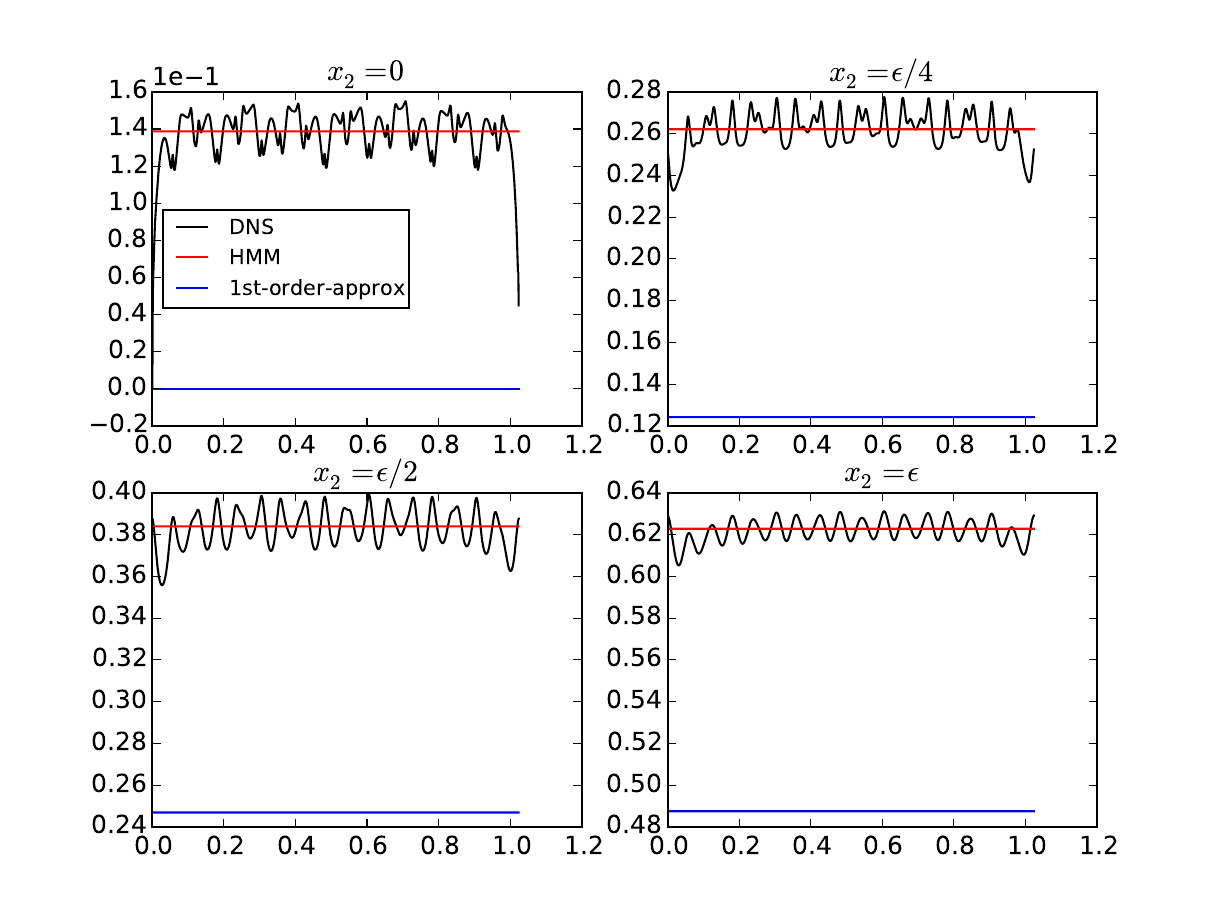}
\caption{Horizontal component of the flow $u_1$ versus $x_1$ plotted at various heights $x_2$
for the domain shown in Figure \ref{fig:domain_d}.}
\label{fig:u1vsx1quasi}
\end{figure}

\begin{figure}[h!]
\centering
\includegraphics[width = 5in, height = 3.33in]{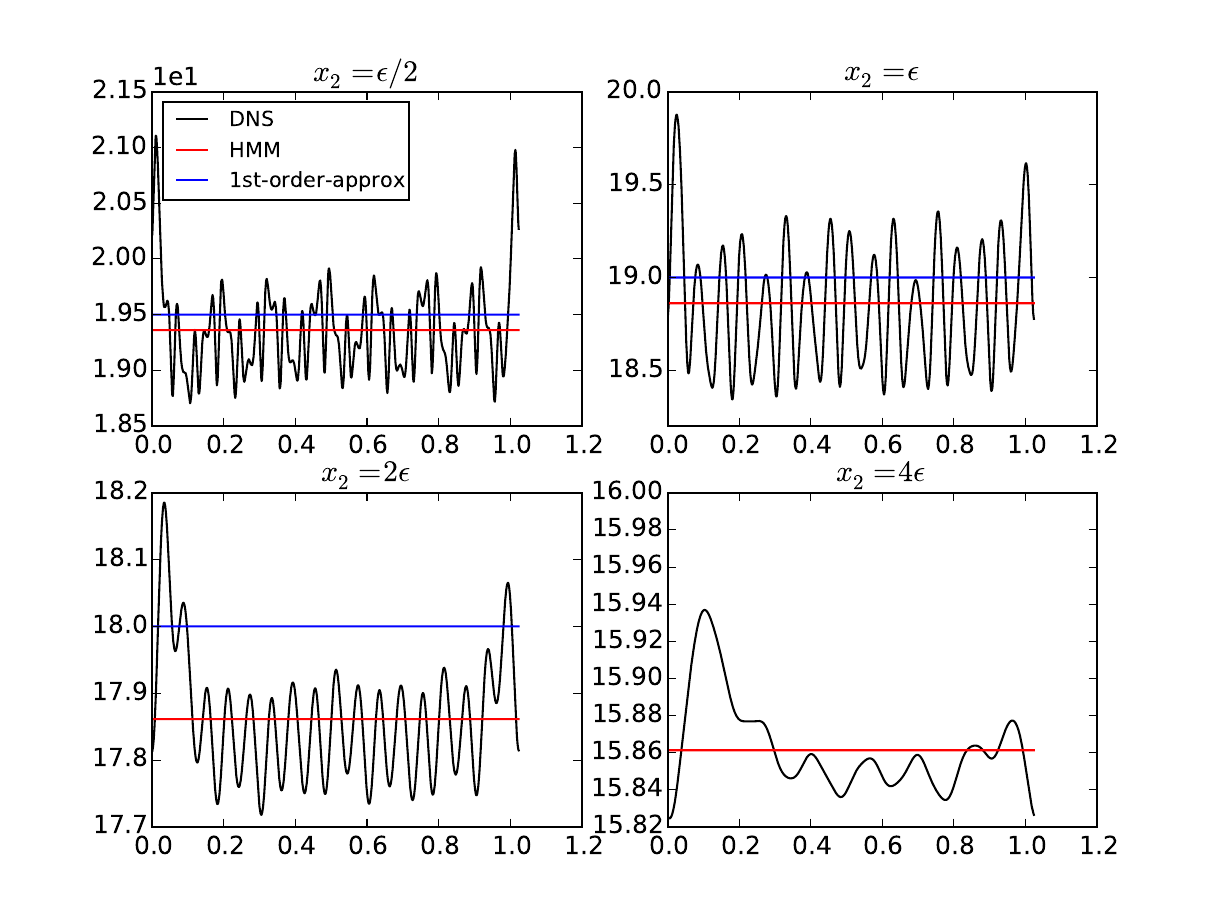}
\caption{Shear $\partial u_1/\partial x_2$ versus $x_1$ plotted at various heights $x_2$
for the domain shown in Figure \ref{fig:domain_d}.}
\label{fig:shearvsx1quasi}
\end{figure}

Consider now a rough boundary parameterized by the quasi-periodic function 
\begin{equation*}
\varphi^{\epsilon}(x_1) = \epsilon/3 \left( \sin(\sqrt{2} \cdot 2\pi x_1/\epsilon) + 
\sin(2\pi x_1/\epsilon) - 2.25\right),
\end{equation*}
like the one displayed in Figure \ref{fig:domain_d}. 
As in Sections \ref{subsection_channel_comp_1} and \ref{subsection_channel_modulated}, 
$\Omegamac = [0,1]^2$, and no-slip is applied at $x_2=0$. 
For both $\Omega^{\epsilon}$ and $\Omegamac$, periodic boundary conditions
are applied at $x_1=0$ and $x_1 = 1$. For $\ueps$ this is only an approximation, 
since $\varphi^{\epsilon}$ is not truly periodic, which explains the
spurious boundary layers in the DNS solution near $x_1=0$ and $x_1=1$, 
as seen in Figures \ref{fig:u1vsx1quasi} and
\ref{fig:shearvsx1quasi}. 

The same problem is encountered in the micro-domain if periodic boundary 
conditions are prescribed. Instead, we use the Dirichlet boundary 
conditions described in Appendix \ref{appendix:micro_bc} 
for the projection operator $\pi_j$. In this case
the no-slip condition is applied at the
locations $(s_j, \varphi^{\epsilon}(s_j))$ and $(s_j + L_j^{\rm mic}, \varphi^{\epsilon}(s_j+L_j^{\rm mic}))$. 
In more general situations such as this where the microscale roughness is no longer
periodic, it is best
to take each $L_j^{\rm mic}>\epsilon$ in order to capture a few ``correlation lengths" 
of $\varphi^{\epsilon}$.
The results shown in Figures \ref{fig:u1vsx1quasi} and \ref{fig:shearvsx1quasi}
are performed with a single micro domain at $s = 0.481561$
 and length $L^{\rm mic}=5\epsilon$. 
%\red{[keep?]}$\Omega^{\epsilon}$ is discretized 
%with 65568 cells, while $\Omegamac$ and $\Omegamic$ are discretized
%with 3280 and 5028 cells, totaling $5\% + 7.7\% = 12.7\%$ the amount
%of rough domain cells. \red{[keep?]}

We note also that if one still wants to use periodic boundary conditions
for the microscale computation
another option is to further increase the horizontal domain length $L$ and 
then replace the smoothing operator from Definition \ref{dfn:smoothing_op} with
\begin{equation*}
\chevron{u}_L(x,y) = \int_{x}^{x+L} K(s) u(s,y) \, ds
\end{equation*}
where $K$ is smooth function that has compact support, integrates to unity, and satisfies
some vanishing moment conditions. Such kernels are well known in the numerical homogenization
community \cite{runborg:2016,tsai:2005,gloria:2011} and likely would be useful in more realistic applications
beyond the academic test cases presented here. 

\subsection{Flow over a backwards facing step.}\label{subsec:bfs}
Consider now flow over a backwards facing step with periodic roughness after the step, as 
in Figure \ref{fig:domain_e}. The roughness is parameterized by the function 
\begin{equation*}
\varphi^{\epsilon}(x) = \epsilon/2\left(\cos(2\pi x_1/\lambda) - 1\right),
\end{equation*}
similar to Section \ref{subsection_channel_comp_1} but with a larger wavelength $\lambda=2.5\epsilon$.
 We are primarily interested in the effect of the roughness
on the flow after the step. Hence for simplicity there is no roughness in the inflow region prior
to the step, and the roughness does not cover the full horizontal extent of the domain. 
In this case, both viscosity $\nu=0.1$ and 
$\epsilon=0.1$, and the horizontal length of the domain is $L=23$. 
A Poiseuille inflow profile drives the flow, and the 
Reynolds number based on the profile is $Re=150$. At this value, some recirculation 
after the step is expected. A zero-stress condition
is applied at the outflow $x_1=23$: $\nu\nabla u - pI = 0$, and for both the full DNS solution and the 1st order
approximation, the no-slip condition is applied on all other domain boundaries. The HMM solution, 
of course however satisfies the slip condition \eqref{eq:mac_wall_law}. 

The algorithm from Section \ref{subsec:algo} is used
with two points $\{s_1, s_2\} = \{7.5, 13.5\}$ (marked with black dots
in Figure \ref{fig:domain_e}) chosen to lie (i) closer to the step, 
and hence within the recirculation bubble,
and (ii) farther away from the step, after the bubble. Given 
$\{\alpha_1, \alpha_2\}$ at these micro domain locations, the slip amount
is given as a piecewise linear interpolant 
\begin{equation}\label{eq:bfs_pw_lin_interp}
\alpha(x_1) = \mathds{1}_{[6,16]}(x_1)
\mathcal{I}_{\rm linear} \left( (6,0), (7.5, \alpha_1), (13.5, \alpha_2), (16,0) \right)(x_1) 
\end{equation}
where $\mathds{1}$ is the indicator function, and $x_1=6$ and $x_1=16$
 are the points at which the roughness begins and ends, and hence
before and after which there should be no-slip. In retrospect a piecewise 
constant interpolant $\mathcal{I}_{\rm constant}$ in the region $6 \le x_1 \le 16$ would
be more appropriate, since \eqref{eq:bfs_pw_lin_interp} does not capture the slip amount
as far out in $x_1$ as it should. Another option would be to simply perform micro simulations
at more points $s_j \in [6,16]$ along the roughness. 

\begin{figure}[h!]
\centering
\includegraphics[width = 5in, height = 3.33in]{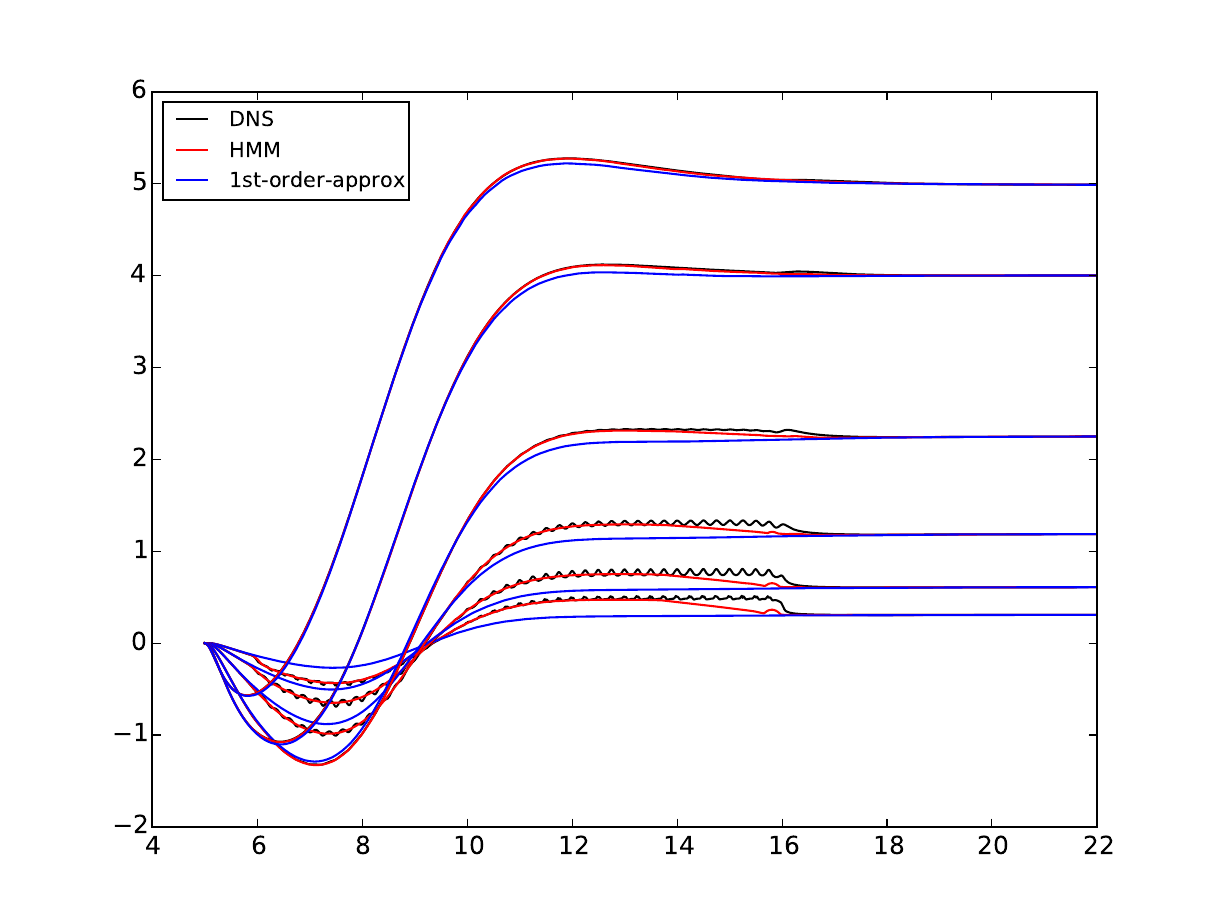}
\caption{Horizontal component of the flow $u_1$ versus $x_1$ plotted at the values 
$x_2 = \epsilon/4, \, \epsilon/2, \, \epsilon, \, 2\epsilon, \, 4\epsilon,\, 0.55$
for the domain shown in Figure \ref{fig:domain_e}.}
\label{fig:u1vsx1bfs}
\end{figure}

\begin{figure}[h!]
\centering
\includegraphics[width = 5in, height = 3.33in]{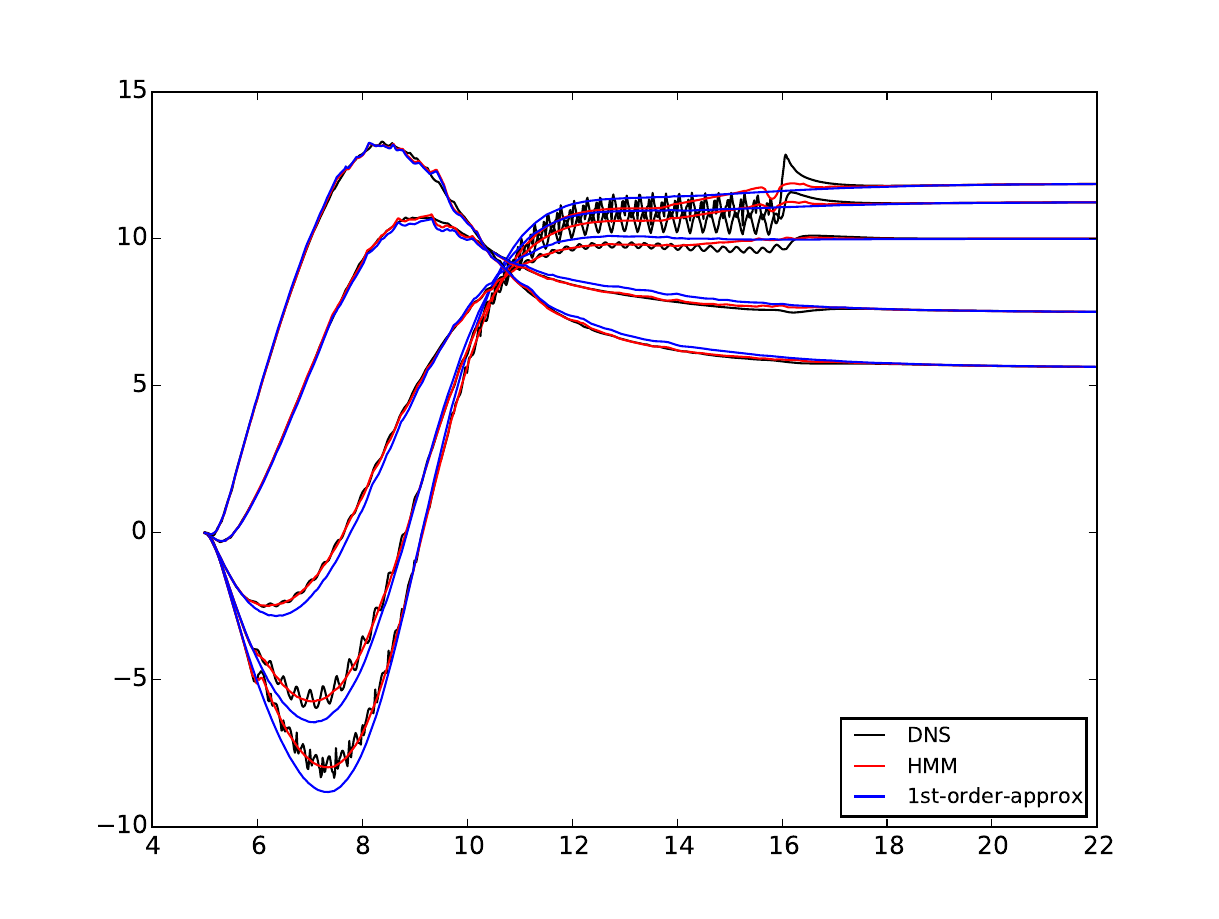}
\caption{Shear $\partial u_1/\partial x_2$ versus $x_1$ plotted at the values 
$x_2= \epsilon/2, \, \epsilon, \, 2\epsilon, \, 4\epsilon,\, 0.55$
for the domain shown in Figure \ref{fig:domain_e}.}
\label{fig:shearvsx1bfs}
\end{figure}

In this case the more general projection operator $\pi_j(U)$ from 
Appendix \ref{appendix:micro_bc} is applied in both micro-domains.
Because of the fluid recirculation, there is a nontrivial mass flux along the upper
computational boundary $x_2 = \gamma = 4\epsilon$ of the micro-domain at $s_1=7.5$. This results
in a 10.2\% difference in the slip amounts computed at $s_1=7.5$ and $s_2=13.5$, even
though the roughness pattern is the same. As a result,
the HMM solution correctly captures the effect of roughness on the size of the recirculation 
bubble, something the 1st-order approximation fails to do. Figures \ref{fig:u1vsx1bfs}
and \ref{fig:shearvsx1bfs} thus illustrate the utility of 
constraining the micro-domains to match the local macroscropic solution.


\subsection{Discussion.}\label{subsec:discussion}
In all cases, the HMM solution computed with the algorithm from Section \ref{subsec:algo} 
is a clear 
improvement over the 1st-order,
no-slip approximation; it captures the average effect of roughness on 
the flow. Although not displayed in the examples above, separate HMM solutions 
were also computed by calculating a slip length at height $x_2 = h>0$
\begin{equation*}
\alpha_{j,h} = \frac{ \langle u^j_1 \rangle_{L_j^{\rm mic}}(s_j,h)}{ 
\langle \partial u^j_1/\partial x_2 \rangle_{L_j^{\rm mic}}(s_j,h)}.
\end{equation*} 
The slip amount fed back to the macro-solver was then defined to be 
\begin{equation}\label{eq:extrap_alpha}
\alpha_j = \alpha_{j, h} - h,
\end{equation}
based on a simple linear extrapolation of $\langle u^j_1\rangle$ from 
$x_2 = h$ to $x_2 = 0$. 
We experimented with both $h = \epsilon$ and $h =2\epsilon$ and found that 
\eqref{eq:extrap_alpha} differed from $\alpha_j$ computed with \eqref{eq:hmm5} 
by at most a few percent. 
%This is consistent with the observation in 
%\cite{pironneau} that the Navier-slip condition \eqref{eq:mac_wall_law}
%is equivalent to a no-slip condition on a flat wall wall at height $\alpha$. 

As mentioned at the beginning of Section \ref{sec:laminar_numerics}, 
the 1st-order approximation satisfying the no-slip condition along $\Gamma$
and the HMM solution are both computed on the same discretization of $\Omegamac$. 
For the first four numerical examples presented 
in Sections \ref{subsection_channel_comp_1}--\ref{subsection:channel_quasiper}, this amounts 
to approximately 4\%--6\% of the mesh cells used to discretize $\Omega^{\epsilon}$. 
After accounting for the mesh cells used to discretize each $\Omegamic_j$, the HMM solution 
is computed on about 8\%--19\% of the mesh cells used to compute the full, oscillatory 
solution $\ueps$ for the parameters described at the 
beginning of Section \ref{sec:laminar_numerics}. Generally speaking, the cost of the micro-solvers
depends on the macroscopic variations of the rough boundary $\Gamma_{\epsilon}$
as well as the correlation length of the oscillations. The former affects the
\emph{number} of micro-domain solutions $J$ while the latter affects the 
\emph{size} of each $\Omegamic_j$. 

In principle, the computational cost of computing each micro-domain solution is 
indepedent of $\epsilon$, since the domain size scales with $\epsilon$. The computational 
savings afforded by solving for the HMM solution
instead of the full, rough-wall solution $\ueps$ of course increase as $\epsilon$ decreases
relative to the measure of $\Omegamac$. For the backwards facing step problem 
described in Section \ref{subsec:bfs}, the system parameters were chosen to highlight the 
effect of the roughness on the size of the recirculation bubble. For the selected
parameters, the total mesh cells for the HMM solution was actually about 20\%  
more than the mesh cells for $\ueps$. For decreasing $\epsilon$, this situation 
would reverse, but the rough-wall's influence on the macroscopic 
flow patterns would be less prominent. 

%Finally, as noted in \cref{remark:one_iter}, the iterative procedure outlined in 
%\cref{alg:hmm} was observed to terminate for all numerical experiments after merely 
%one iteration when the tolerance $\tau = \epsilon^2$ and the initial value $\alpha_{(0)} =0$
%was used. Despite the quick convergence, 
%the need to have some initial $U_{(0)}$ in hand from which to estimate boundary conditions 
%$\Upsilon_0 = \pi(U_{(0)})$ for the micro-solvers is a disadvantage of the proposed HMM scheme. In practice, 
%the overhead cost to produce some $\Upsilon_0$ could be reduced by computing $U_{(0)}$ 
%on an even coarser mesh than that of $U_{(1)}$. For time-dependent problems, the iterative
%procedure could be incorporated into a time-marching scheme, such as a predictor-corrector method,
%quite naturally. 

%Of particular note is the final example of a backwards facing 
%step, for which the HMM solution correctly captures the effect of the 
%roughness on the size of the recirculation bubble. 

\section{Conclusions and future directions}
\label{sec:conclusions}

We develop computational techniques for computing effective boundary conditions
for laminar, viscous flow over a rough surface. The technique is based on coupling 
a `macroscale' solver $M$ in a large domain with smooth boundary to one or more 
`microscale' solvers $m$ localized to the the rough surface. The microscale solvers
estimate the slip amount in the wall-law utilized by the macrosolver. 
 The coupling strategy is described in the framework of the heterogeneous 
multiscale method (HMM). So long as the surface roughness is asymptotically small 
relative to the size of the full domain, the method is designed 
to efficiently capture the average effect of surface roughness on the flow for 
arbitrary roughness patterns. Numerical examples illustrate the method's utility 
for a variety of cases.

We prove the coupled system is consistent for linear Stokes flow in a channel 
with periodic roughness; starting from a no-slip Poiseuille flow, a back-and-forth
iteration will converge to a fixed point. Moreover, we show the convergence happens 
rapidly--the difference between the slip amounts after the first and second iteration
is smaller than $\bigoh{\epsilon^2}$, where $\epsilon$ is the small-scale parameter
characterizing the surface roughness. Finally, we show the HMM slip amount vanishes
as $\epsilon$ vanishes, so that the macroscale HMM flow and the true, rough-wall flow $\ueps$
both converge to the same quantity, namely, Poiseuille flow with no-slip condition 
posed on a smooth boundary. A future challenge remains, however, to rigorously estimate
the difference between the HMM approximation and $\ueps$, 
either by adapting the currently available homogenization techniques 
\cite{jager,jager_couette,mikelic2013} or 
by inventing new ones. 

Laminar incompressible flow over a bed of porous media is a 
setting analogous to the current one; future work could build
off the studies of \cite{Lacis:2017} and \cite{Sudhakar2019HigherOH} 
and analyze a HMM that couples local micro-solvers based on Darcy's 
law to Navier-Stokes macro-solvers for the large-scale flow.
In that case, mathematically rigorous results also exist that suggest
the form of the wall-law at the interface of the fluid and porous 
media \cite{jager_porous}.

%\blue{Asymptotic arguments 
%demonstrate that the model reproduces the results from mathematical homogenization theory 
%whenever the theory may be applied, for instance, when the roughness is periodic or the 
%large domain has a channel geometry. Numerical examples illustrate the model's utility
%for more general situations. }

%In the relatively simple case of laminar flow, mathematically
%rigorous results exist for form of the wall-law that captures 
%the the asymptotic behavior of viscous, rough-wall flow. 
%As such, the current setting is ideal for designing computational techniques to
%estimate effective boundary conditions that homogenize the boundary layer induced 
%by the roughness. 

%An analogous setting is
%laminar incompressible flow over a bed of porous media; 
%mathematical justification of 
%a wall-law, the Beavers-Joseph-Saffman law, at the interface
%of the fluid and porous media is provided in \cite{jager_porous}, for example. Future 
%work could build off the
%studies of \cite{Lacis:2017} and \cite{Sudhakar2019HigherOH} and analyze a HMM that
%couples local micro-solvers based on Darcy's law to Stokes or Navier-Stokes 
%macro-solvers for the large-scale flow.

For numerous other boundary layer problems in fluid mechanics however
there does \emph{not} exist any rigorous theory describing an effective 
boundary condition that captures the system's asymptotic dynamics.
The present work offers some indication that
a computational approach based on coupling macro and microscale simulations could 
be successfully used to numerically derive wall-laws in such cases.
 
%However, there are numerous boundary layer problems in fluid mechanics for which rigorous 
%theory describing a wall-law that captures the system's asymptotic 
%dynamics does not exist. For such situations, similar computational techniques
%based on coupling macro and microscale simulations could be used to numerically derive
%wall-laws. 

One example is wall-bounded electrokinetic flows. In the presence of charged surfaces,
an asymptotically thin charge layer forms \cite{squires_bazant_2004,yariv:2012}. 
This layer is the primary source of both momentum and charge transport for such flows; 
for flows relevant to microfluidics, however, its
physics are poorly represented by continuum models. An HMM approach would  
couple a relatively efficient continuum solver $M$ for the bulk flow to a molecular 
model $m$ for the near-wall charge layer \cite{hmm_microfluidics,bell:2020}. 
One challenge is representing the effect of thermal fluctuations 
in both macro and microscale simulations in a consistent way \cite{donev:2010}. 

Another example is wall-bounded turbulent flow;  
when the friction Reynolds number is 
large enough, there is a clear separation of scales between the near-wall eddies and 
those in the bulk flow \cite{smits:2013}. Based on this observation, Sandham et al.\ 
coupled a large eddy simulation (LES) to a sequence of micro-solvers localized to 
the wall in which direct numerical simulation (DNS) was performed \cite{sandham}. In 
the language of HMM, the coupling was done concurrently, or `on-the-fly'. 
In \cite{carney_engquist_moser_2020} we built on this approach and further developed 
the local DNS solvers for use in a sequential multiscale framework, i.e.\ so that
precomputed information could be utilized by a macroscale LES. 
More work is needed to design a wall-law for an LES that can use
such information \cite{piomelli:2002}. Intermittency remains a formidable obstacle.

%FIGURES
%The use of .eps files is encouraged, in which case you should
%un-comment the \uspackage{graphics} and \uspackage{epstopdf} command above, and use the
%command \include{figure.eps} to insert the figure file.

\section*{Acknowledgments}
The authors thank Yoonsang Lee for discussions regarding the heterogeneous multiscale method, 
as well as Matthew Novack for discussions regarding the proof of the technique's convergence 
to a fixed point. The authors acknowledge support from the Oden Institute for Computational
 Engineering and Sciences and National Science Foundation Grant DMS-1620396.

\appendix
\section{Boundary conditions for micro-scale systems}
\label{appendix:micro_bc}
In more general situations in which the macroscopic flow has nontrivial dependence on $x_1$ and/or
a nontrivial vertical component, a more general approach is needed than the ``free-stream" condition 
mentioned in Remark \ref{rmk:bc_micro}.
We propose to prescribe
quadratic Dirichlet conditions for both the horizontal and vertical components of the velocity
 along each of the three faces of $\partial \Omegamic_{j, \rm D}$ 
(those that intersect $x_1 = s_j$, $x_1 = s_j + \epsilon$ and $x_2 = \gamma$). 
Let 
\begin{align*}
\Gamma_1^{\rm mic} &= \{(x_1,x_2)\vert 0 \le x_2 \le \gamma, \, x_1 = s_j\} \\
\Gamma_2^{\rm mic} &= \{(x_1,x_2)\vert s_j \le x_1 \le s_j + L_j^{\rm mic}, \, x_2 = \gamma\} \\
\Gamma_3^{\rm mic} &= \{(x_1,x_2)\vert 0 \le x_2 \le \gamma, \, x_1 = s_j+L_j^{\rm mic} \} 
\end{align*}
be the left, upper, and right computational boundaries of the micro-domain (the dependence of 
each $\Gamma^{\rm mic}$ on $j$ is implied). Then there are two quadratic profiles
for each boundary, each with three coefficients to be determined, and hence 18 total 
constraints are needed. Let $u_k$ and $v_k$ be the quadratic profile for the horizontal and vertical  
component of the flow at $\Gamma_k^{\rm mic}$, $k=1,2,3$. 
The no-slip requirement \eqref{eq:micro_cont} gives four constraints
\begin{equation*}\label{eq:first_pi_constraint}
0= u_1(s_j, 0) = v_1(s_j, 0) = u_3(s_j+L_j^{\rm mic},0) = v_3(s_j + L_j^{\rm mic} ,0). 
\end{equation*}
Additionally, enforce that the mass flux across each $\Gamma^{\rm mic}_k$ is the same
as the macroscopic mass flux 
\begin{align}
\int_{\Gamma_1^{\rm mic}} u_1\,  ds &= \int_{\Gamma_1^{\rm mic}} U \cdot n\,\, ds \nonumber \\
\int_{\Gamma_2^{\rm mic}} v_2\,  ds &= \int_{\Gamma_2^{\rm mic}} U \cdot n\,\, ds \nonumber \\
\int_{\Gamma_3^{\rm mic}} u_3\,  ds &= \int_{\Gamma_3^{\rm mic}} U \cdot n\,\, ds.\label{eq:flux_mic3} 
\end{align}
Since $U$ is divergence free, these imply conservation of mass along the micro-domain boundaries
\begin{equation*}
\int_{\Gamma_1^{\rm mic}} u_1\,  ds + \int_{\Gamma_2^{\rm mic}} v_2\,  ds + 
\int_{\Gamma_3^{\rm mic}} u_3\,ds= 0, 
\end{equation*}
hence satisfying requirement \eqref{eq:cons_of_mass}.
Equations \eqref{eq:flux_mic3} hence give three more conditions. 
To completely specify the quadratic profiles, one more condition each is needed
for $u_1$ and $u_3$, two more conditions are needed for $v_1$, $v_2$, and $v_3$, and three conditions are
needed for $u_2$. For continuity, enforce the interpolation constraints 
\begin{align*}
u_1(s_j,\gamma) = u_2(s_j,\gamma) =  U(s_j, \gamma) \cdot e_1 &\\
v_1(s_j,\gamma) = v_2(s_j,\gamma) =  U(s_j, \gamma) \cdot e_2 &\\
u_2(s_j+L_j^{\rm mic},\gamma) = u_3(s_j+L_j^{\rm mic}, \gamma) = U(s_j+L_j^{\rm mic},\gamma) \cdot e_1  &\\
v_2(s_j+L_j^{\rm mic},\gamma) = v_3(s_j+L_j^{\rm mic}, \gamma) = U(s_j+L_j^{\rm mic}, \gamma) \cdot e_2 &, 
\end{align*} 
which leaves one more condition each for $v_1$, $u_2$, and $v_3$. Adding one more interpolation 
point for each 
\begin{align*}
v_1(s_j, \gamma/2) = U(s_j, \gamma/2)\cdot e_2& \\ 
u_2(s_j + L_j^{\rm mic}/2, \gamma) = U(s_j + L_j^{\rm mic}/2, \gamma)\cdot e_1& \\ 
v_3(s_j+L_j^{\rm mic}, \gamma/2) = U(s_j+L_j^{\rm mic},\gamma/2)\cdot e_2&. \label{eq:last_pi_constraint}
\end{align*} 
ensures the Dirichlet conditions are thus uniquely determined along each $\Gamma^{\rm mic}_k$, 
which then completely defines a projection operator $\pi_j$.

\section{Estimates for Stokes flow}
\label{appendix:estimates}
By definition 
\begin{equation*}
\overline B_2 = \int_{0}^{\epsilon} u_2(x_1,0)\, dx
\end{equation*}
meaning 
\begin{equation*}
|\overline{B}_2  | \le \int_0^{\epsilon} |u_2(x_1,0)| \, dx \le \norm{u_2}_{L^1} 
\le |\Omega_{\rm mic}|^{1/2} \norm{u_2}_{L^2} \le [(M+\gamma_2)\epsilon]^{1/2} \norm{u_2}_{H^1}
\end{equation*}
By the Lax-Milgram theorem \cite{arbogast:2008}, we know 
\begin{equation*}
\norm{u_2}_{H^1} \le (2+C_p^2) \norm{g_2}_{H^1}, 
\end{equation*}
where $C_p$ is the Poincar\'{e} constant (which is estimated in Appendix \ref{appendix:poincare}) and $g_2$ is the (nonunique) pre-image of the trace-map $\textrm{Tr}$ of the boundary values for $u_2$: 
\begin{equation*}
\textrm{Tr}(u_2) = 
\begin{cases*} 
0, \qquad\qquad x_2=\gamma_2 \\
-u_1\big\vert_{\Gamma_{\epsilon}}, \qquad (x_1,x_2) \in \Gamma_{\epsilon}.
\end{cases*}
\end{equation*}
Let 
\begin{align*}
g_2(x_1,x_2) = m(x_1) (x_2-\gamma_2) &= \left(\frac{u_1\big(-\varphi^{\epsilon}(x_1)\big)}
{\gamma_2 + \varphi^{\epsilon}(x_1)}\right) (x_2-\gamma_2) \nonumber \\
&= \left(\frac{-(M-\varphi^{\epsilon}(x_1))(\varphi^{\epsilon}(x_1)+\gamma_2)}{\gamma_2 + \varphi^{\epsilon}(x_1)}\right) (x_2-\gamma_2),
\end{align*}
then $g_2(x_1,\gamma_2) = 0$ and $g_2(x_1,-\varphi^{\epsilon}(x_1)) = -u_1(-\varphi^{\epsilon}(x_1))$, i.e.\ $\rm{Tr}(g_2) = \rm{Tr}(u_2)$. 
The gradient is 
\begin{equation*}
\pderiv{g_2}{x_1} = \pderiv{m}{x_1} (x_2-\gamma_2) = \varphi'(x_1/\epsilon) (x_2-\gamma_2) 
\end{equation*}
and 
\begin{equation*}
\pderiv{g_2}{x_2} = m(x_1) = \frac{u_1\big(-\varphi^{\epsilon}(x_1)\big)}{\gamma_2 + \varphi^{\epsilon}(x_1)}.
\end{equation*}
Measuring each term in the $H^1$ norm: 
\begin{align*}
\norm{g_2}_{L^2}^2 &= \int_0^{\epsilon} \int_{-\varphi^{\epsilon}(x_1)}^{\gamma_2} \big(m(x_1)\big)^2 (x_2-\gamma_2)^2 \, dx_2 \, dx_1 \nonumber \\
&= \frac13\int_0^{\epsilon} \big(m(x_1)\big)^2 (\varphi^{\epsilon}(x_1)+\gamma_2)^3 \, dx_1  \nonumber \\
&= \frac13\int_0^{\epsilon} \left[u_1\big(-\varphi^{\epsilon}(x_1)\big)\right]^2 (\varphi^{\epsilon}(x_1)+\gamma_2)\, dx_1 \nonumber \\
&\le \frac13 \,\epsilon\,  (M\gamma_2)^2 (M+\gamma_2)
\end{align*}
since $u_1(\varphi^{\epsilon})$ is largest at $y=0$ and $\varphi^{\epsilon}(x) \le M$. 
\begin{align*}
\norm{\partial g_2/\partial x_1}_{L^2}^2 &= \int_0^{\epsilon} \int_{-\varphi^{\epsilon}(x_1)}^{\gamma_2} 
\left(\varphi'(x_1/\epsilon)\right)^2 (x_2-\gamma_2)^2 \, dx_2 \, dx_1  \nonumber \\
&= \frac13\int_0^{\epsilon} \left(\varphi'(x_1/\epsilon)\right)^2 (\varphi^{\epsilon}(x_1)+\gamma_2)^3 \, dx_1 \nonumber \\
&\le \frac13 \,\epsilon\,  (M+\gamma_2)^3 \norm{\varphi'}_{\infty}^2 .
\end{align*}
Finally, 
\begin{align*}
\norm{\partial g_2/\partial x_2}_{L^2}^2 &= \int_0^{\epsilon} \int_{-\varphi^{\epsilon}(x_1)}^{\gamma_2} \left(m(x_1)\right)^2 \, dx_2 \, dx_1 \nonumber \\
 &= \int_0^{\epsilon} \frac{\left(u_1\big(-\varphi^{\epsilon}(x_1)\big)\right)^2}{\gamma_2 + \varphi^{\epsilon}(x_1)} \, dx_2 \, dx_1 \nonumber \\
&\le \frac{1}{\gamma_2} \int_0^{\epsilon} \left(u_1\big(-\varphi^{\epsilon}(x_1)\big)\right)^2  \, dx_1 \nonumber \\
&\le \frac{\epsilon}{\gamma_2} \, (M\gamma_2)^2 = \epsilon \gamma_2 M^2  .
\end{align*}
Putting it all together:  %%%%%%%%%%%%%%%%%%%%%%%%%%%%%%
\begin{equation} \label{eq:sigma_defn}
\norm{g_2}_{H^1} \le \left[ \frac13 \,\epsilon\,  (M\gamma_2)^2 (M+\gamma_2) + 
\frac13 \,\epsilon\,  (M+\gamma_2)^3 \norm{\varphi'}_{\infty}^2 + \epsilon \gamma_2 M^2  \right]^{1/2} =: \sigma_{\epsilon},
\end{equation}
so that in total          %%%%%%%%%%%%%%%%%%%%%%%%%%%%%%
\begin{equation} \label{eq:B2_bound}
|\overline{B}_2  | \le (2 + C_p^2) [(M+\gamma_2)\epsilon]^{1/2} \sigma_{\epsilon}.
\end{equation}

Next we bound $\tilde B_3$. By definition 
\begin{equation*}
\tilde B_3 = \int_{0}^{\epsilon} \tilde u_3(x_1,0)\, dx_1
\end{equation*}
meaning 
\begin{equation*}
|\tilde B_3  | \le \int_0^{\epsilon} |\tilde u_3(x_1,0)| \, dx_1 \le \norm{\tilde u_3}_{L^1} 
\le |\Omega_{\rm mic}|^{1/2} \norm{\tilde u_3}_{L^2} \le [(M+\gamma_2)\epsilon]^{1/2} \norm{\tilde u_3}_{H^1}
\end{equation*}
Again we know by the Lax-Milgram theorem \cite{arbogast:2008} that
\begin{equation*}
\norm{\tilde u_3}_{H^1} \le (2+C_p^2) \norm{g_3}_{H^1}, 
\end{equation*}
where $g_3$ is the pre-image of the trace-map $\textrm{Tr}$ of the boundary values for $\tilde u_3$: 
\begin{equation*}
\textrm{Tr}(\tilde u_3) = 
\begin{cases} 
1, \qquad x_2=\gamma_2 \\
0, \qquad (x_1,x_2) \in \Gamma_{\epsilon} 
\end{cases}
\end{equation*}
Let 
\begin{equation*}
g_3(x_1,x_2) = 
\begin{cases} 
x_2/\gamma_2, \qquad (x_1,x_2) \in [0,\epsilon]\times [0, \gamma_2] \\
0, \qquad (x_1,x_2) \in \Omega_{\rm mic}\setminus[0,\epsilon]\times [0, \gamma_2] .
\end{cases}
\end{equation*}
Then 
\begin{equation*}
\norm{g_3}_{H^1}^2 = \int_{0}^{\epsilon} \int_{0}^{\gamma_2} \left(x_2^2+1\right)/\gamma_2^2 \, dx_2 \, dx_1
= \epsilon \left(\frac{\gamma_2}{3} + \frac{1}{\gamma_2}\right)
\end{equation*}
meaning
\begin{equation}\label{eq:B3_bound}
|\tilde B_3| \le (2+C_p^2) [(M+\gamma_2)\epsilon]^{1/2}  \norm{g_3}_{H^1} 
= \epsilon [(M+\gamma_2)]^{1/2} (2+C_p^2)\left(\frac{\gamma_2}{3} + \frac{1}{\gamma_2}\right)^{1/2} .
\end{equation}

\section{Estimate for Poincar\'{e} constant}
\label{appendix:poincare}
The bounds \eqref{eq:B2_bound} and \eqref{eq:B3_bound} both depend on the Poincar\'{e} constant 
$C_p$, i.e.\ the constant such that 
\begin{equation*}
\norm{u}_{L^2(\Omegamic)}^2 \le C_p^2 \norm{\nabla u}_{L^2(\Omegamic)}^2 \qquad \forall u \in \mathcal{H}
\end{equation*}
where 
\begin{equation*}
\mathcal{H} = \{\psi \in H^1(\Omegamic) \,\, : \,\, 
\psi\vert_{\Gamma_{\epsilon}} = \psi\vert_{x_2 = \gamma_2} = 0, \, \,  \psi \text{ is periodic at } x_1 = 0 \text{ and }  x_1 = \epsilon \}.
\end{equation*}
We now show $C_p \to 0$ as $\epsilon \searrow 0$, which is necessary to determine the
asymptotic character of $\overline{B}_2$ and $\tilde B_3$. 

Let $u \in \mathcal{H}$. For $x_1 \in [0, \epsilon)$, we know from direct computation 
of the spectrum of the operator $-\partial^2/\partial x_2^2$ acting on the interval
$[-\varphi^{\epsilon}(x_1), \gamma_2]$ with homogeneous Dirichlet boundary conditions 
that 
\begin{equation*}
\int_{-\varphi^{\epsilon}(x_1)}^{\gamma_2} |u(x_1,x_2)|^2 \, dx_2 \le 
\left(\frac{L_y(x_1)}{\pi}\right)^2 
\int_{-\varphi^{\epsilon}(x_1)}^{\gamma_2} \left| \pderiv{u}{x_2} (x_1,x_2)\right|^2 \, dx_2
\end{equation*}
where $L_y(x_1) = \gamma_2 + \varphi^{\epsilon}(x_1) \le \gamma_2 + M$. Integrating
both sides of the inequality from $x_1 = 0$ to $x_1 = \epsilon$ then gives 
\begin{align}
\norm{u}_{L^2(\Omegamic)}^2 &\le \left(\frac{\gamma_2 + M}{\pi}\right)^2 \norm{\pderiv{u}{x_2}}_{L^2(\Omegamic)}^2 \nonumber \\
&\le \left(\frac{\gamma_2 + M}{\pi}\right)^2 \norm{\nabla u}_{L^2(\Omegamic)}^2
\end{align} 
so that indeed $C_p \to 0$ as $\epsilon \searrow 0 $ as desired.

\bibliographystyle{siamplain}
\bibliography{references}

%\par{\bf References.}\, Always use $\backslash$cite$\{biblabelname\}$ (eg. \cite{taubes1}) to cite
%          references which have been named in the bibliography via
%          $\backslash$bibitem$\{biblabelname\}$.
%
%
%          % Non-BibTeX users please use
%
%
%          % and use \bibitem to create references.
%
%           \bibitem{mnemonic}Author, {\em Article's title}, Journal Name,
%         Volume Number(Issue Number):page numbers, year.
%
%          % Format for Journal Reference. For example
%
%          \bibitem{taubes1} C. Taubes, {\em The Seiberg-Witten invariants
%          and
%          symplectic forms}, Math. Res. Letters, 1:809--822, 1994.

\end{document}